\documentclass{article}
\usepackage{amsmath,amsthm,amssymb}

\title{Elementary Equivalence of Categories\\
of Modules over Rings, Endomorphism Rings,\\
and Automorphism Groups of Modules}

\author{E.~I.~Bunina, A.~V.~Mikhalev}

\newtheorem{theorem}{Theorem}
\newtheorem{corollary}{Corollary}
\newtheorem*{corol}{Corollary}
\newtheorem{lemma}{Lemma}
\newtheorem{proposition}{Proposition}

\makeatletter
\@addtoreset{theorem}{section}
\@addtoreset{lemma}{section}
\@addtoreset{proposition}{section}
\@addtoreset{corollary}{section}
\makeatother

\newcommand{\prisv}{\mathbin{{:}\!=}}

\newcommand{\GL}{\mathrm{GL}}
\newcommand{\PGL}{\mathrm{PGL}}
\newcommand{\PSL}{\mathrm{PSL}}
\newcommand{\SL}{\mathrm{SL}}
\newcommand{\mU}{\mathrm{U}}
\newcommand{\SO}{\mathrm{SO}}
\newcommand{\PSO}{\mathrm{PSO}}
\newcommand{\Spin}{\mathrm{Spin}}
\newcommand{\mSp}{\mathrm{Sp}}
\newcommand{\PSp}{\mathrm{PSp}}
\newcommand{\Mor}{\mathop{\mathrm{Mor}}\nolimits}
\newcommand{\Obj}{\mathop{\mathrm{Obj}}\nolimits}
\newcommand{\ring}{\mathsf{ring}}
\newcommand{\Cn}{\mathrm{Cn}}
\newcommand{\Th}{\mathrm{Th}}
\newcommand{\Endom}{\mathop{\mathrm{End}}\nolimits}
\newcommand{\Aut}{\mathop{\mathrm{Aut}}\nolimits}
\newcommand{\Dom}{\mathrm{Dom}}
\newcommand{\Un}{\mathrm{Un}}
\newcommand{\On}{\mathrm{On}}
\newcommand{\card}{\mathop{\mathrm{card}}\nolimits}
\newcommand{\cf}{\mathrm{cf}}
\newcommand{\rng}{\mathop{\mathrm{rng}}\nolimits}
\newcommand{\ucup}{\mathop{{\cup}}}
\newcommand{\Simp}{\mathop{\mathrm{Simp}}\nolimits}
\newcommand{\tSum}{\mathop{\mathrm{Sum}}\nolimits}
\newcommand{\Fin}{\mathop{\mathrm{Fin}}\nolimits}
\newcommand{\Pret}{\mathop{\mathrm{Pret}}\nolimits}
\newcommand{\Under}{\mathop{\mathrm{Under}}\nolimits}
\newcommand{\Und}{\mathop{\mathrm{Und}}\nolimits}
\newcommand{\Finite}{\mathop{\mathrm{Finite}}\nolimits}
\newcommand{\Proobr}{\mathop{\mathrm{Proobr}}\nolimits}
\newcommand{\Gr}{\mathrm{Gr}}
\newcommand{\Rng}{\mathop{\mathrm{Rng}}\nolimits}
\newcommand{\Cl}{\mathop{\mathrm{Cl}}\nolimits}
\newcommand{\Gen}{\mathop{\mathrm{Gen}}\nolimits}
\newcommand{\Gener}{\mathrm{Gener}}
\newcommand{\Ring}{\mathrm{Ring}}
\newcommand{\Sets}{\mathrm{Sets}}
\newcommand{\Select}{\mathrm{Select}}
\newcommand{\Surj}{\mathrm{Surj}}
\newcommand{\Inj}{\mathrm{Inj}}
\newcommand{\Bij}{\mathrm{Bij}}
\newcommand{\Inf}{\mathrm{Inf}}
\newcommand{\Next}{\mathrm{Next}}
\newcommand{\Matrix}{\mathrm{Matrix}}
\newcommand{\Func}{\mathrm{Func}}
\newcommand{\Mod}{\mathrm{Mod}}
\newcommand{\Undir}{\mathrm{Undir}}
\newcommand{\Ker}{\mathop{\mathrm{Ker}}\nolimits}
\newcommand{\Image}{\mathop{\mathrm{Im}}\nolimits}

\newenvironment{mVentry}[1]%
{\begin{list}{}{%
\setlength{\itemsep}{-\parsep}%
\settowidth{\labelwidth}{#1}%
\setlength{\leftmargin}{\labelwidth}%
\addtolength{\leftmargin}{\labelsep}}}%
{\end{list}}

\newcommand{\logic}[1]{\mathrel{#1}}
\newcommand{\losp}{\,}
\newcommand{\ttm}{\mskip-1mu}

\allowdisplaybreaks

\begin{document}

\maketitle

\tableofcontents

\section*{Introduction}
\addcontentsline{toc}{section}{\hspace*{2.3em}Introduction}
The \emph{first order language} (see Sec.~\ref{ss1.1})
of some algebraic theory (for example the
group theory or the ring theory) is the language, where in formulas
we use quantifiers $\forall$~and~$\exists$, logical symbols
$\neg$, $\land$, $\lor$, $\Rightarrow$,
parentheses and variables, and also
predicate and function symbols, and constant symbols
of this theory. For example, in the group theory we use
the subformulas $x\cdot y$, $x^{-1}$, $1$,
in the ring theory we use
the subformulas $x\cdot y$, $x^{-1}$, $1$, $x+y$, $-x$, $0$, and so on.

Two models $\mathcal U$~and~$\mathcal V$ of the language~$\mathcal L$
(for example, two groups or two rings) are called \emph{elementarily
equivalent} if for every sentence~$\varphi$ of the language~$\mathcal L$
we have that it is true in~$\mathcal U$ if and only if
it is true in~$\mathcal V$.
We denote this relation between models by
$\mathcal U \equiv \mathcal V$.

The first result in elementary equivalence
of linear groups was proved by A.~I.~Maltsev
in~1961 (see~\cite{3}). He proved the following theorem.

\begin{theorem}
The group $G_m(K_1)$ is elementarily equivalent to the group
$G_n(K_2)$ \textup{(}$G=\GL,\PGL,\SL,\PSL$, $m\ge n\ge 3$, $K_1$ and~$K_2$
are fields of characteristic~$0$\textup{)} if and only if $m=n$
and $K_1\equiv K_2$.
\end{theorem}

In his proof of this theorem A.~I.~Maltsev used the Jordan normal
form of matrices and explained how to write for each matrix~$M$
a~formula $\varphi(A)$ which is true in the given group
if and only if the matrix $A$ has the same Jordan form
as the matrix~$M$.

If we consider linear groups over skewfields or rings
we still do not have any adequate analogue of the theory
of Jordan normal forms.

But recent progress in the model theory (the construction of ultraproducts
and ultrapowers) (see~\cite{4} and also Sec.~\ref{ss1.4}) has
helped us to continue investigations in this field.
Using this construction in 1992 C.~I.~Beidar and
A.~V.~Mikhalev formulated a~general approach to
problems of elementary equivalence of different algebraic structures
(see~\cite{2}).
Taking into account some results in the theory of linear groups over
rings, they obtained easy proofs
of Maltsev-type theorems in rather general
situations (for linear groups over prime rings,
for multiplicative semigroups, lattices of submodules, and so on).

We give some of their results which extend the Maltsev theorem.

\begin{theorem}
Let $R$ and $S$ be prime associative rings with~$1$
\textup{(}$1/2$\textup{)} and
$m,n\ge 3$ \textup{(}$m,n\ge 2$\textup{)}.
Then $\GL_m(R)\equiv \GL_n(S)$ if and only
if either $M_m(R)\equiv M_n(S)$ or $M_m(R)\equiv M_n(S)^{\mathrm{op}}$.
\end{theorem}

\begin{theorem}
Let $R$ and $S$ be skewfields and $m,n\ge 3$. Then
$\GL_m(R)\equiv \GL_n(S)$ if and only if either $m=n$ and
$R\equiv S$ or $m=n$ and $R\equiv S^{\mathrm{op}}$.
\end{theorem}

In 1998--2001 E.~I.~Bunina continued to study elementary properties
of linear groups (see \cite{12,13,14,15}).
In 1998 (see \cite{12,15}) the results of A.~I.~Maltsev were
generalized to unitary linear groups over fields with involution.
The proof, as in the paper~\cite{3} of A.~I.~Maltsev,
was based on the Jordan normal form of matrices.

Let $K$ be an infinite field with characteristic not equal to~$2$ and
with an involution~$j$ (an involution is an antiautomorphism of order~$2$),
$M_n(K)$ be the total $(n\times n)$-matrix ring over~$K$,
and $\GL_n(K)$ be the linear group over~$K$.
Let $Q_{2n}$ be the following matrix from $\GL_{2n}(K)$:
$$
\left.
\begin{pmatrix}
0 & 1& \dots& 0& 0\\
-1 & 0& \dots& 0& 0\\
& & \ddots\\
0& 0& \dots & 0 & 1\\
0& 0&\dots &-1 & 0
\end{pmatrix}
\right\} {\scriptstyle 2n}.
$$
Let $\mU_{2n}(K,j,Q)$ be the \emph{unitary group} of all matrices
$A\in \GL_{2n}(K)$ such that
$AQ_{2n} A^*=Q_{2n}$,
where
$$
A^*=(A^j)^{\mathrm T}=
\begin{pmatrix}
a_{11}^j & \ldots & a_{1n}^j\\
\hdotsfor{3}\\
a_{n1}^j & \ldots & a_{nn}^j
\end{pmatrix}^{\mathrm T}=
\begin{pmatrix}
a_{11}^j & \ldots & a_{n1}^j\\
\hdotsfor{3}\\
a_{1n}^j & \ldots & a_{nn}^j
\end{pmatrix}.
$$

The following theorem was proved by E.~I.~Bunina.

\begin{theorem}
If $K_1$ and $K_2$ are infinite fields
of characteristic not equal to~$2$ with involutions $j_1$~and~$j_2$,
respectively, and $n,m\ge 2$,
then
the groups $\mU_{2n}(K_1,j_1,Q_{2n})$ and $\mU_{2m}(K_2,j_2,Q_{2m})$
are elementarily equivalent
if and only if $m=n$ and the fields $K_1$~and~$K_2$
are elementarily equivalent as fields with involution.
\end{theorem}

Elementary equivalence of fields with involution means that
in sentences together with the ring operations we use
the operation of involution.

As it was done for linear groups over rings,
using the construction of ultraproducts, E.~I.~Bunina in 1998
(see \cite{13,15}) considered elementary equivalence of unitary
linear groups over rings and skewfields with involution.

Involution in a~ring~$K$ is an antiautomorphism of order~$2$,
i.e., it is a~bijective mapping~$j$ from the ring~$K$ onto itself
such that
\begin{enumerate}
\item
$j(a+b)=j(a)+j(b)$ for all $a,b\in K$;
\item
$j(a\cdot b)=j(b)\cdot j(a)$ for all $a,b\in K$;
\item
$j^2(a)=j(j(a))=a$ for all $a\in K$.
\end{enumerate}

If $K$ is a~ring with involution~$j$, then
by~$\tau$ we shall denote the involution of the ring $M_{2n}(K)$
of matrices over~$K$ having the form
$$
\tau\colon
A=
\begin{pmatrix}
a_{11}& \ldots& a_{1\,2n}\\
\hdotsfor{3}\\
a_{2n\,1}& \ldots& a_{2n\,2n}
\end{pmatrix}
\mapsto
Q_{2n}\circ
\begin{pmatrix}
a^j_{1\,1}& \ldots& a^j_{2n\,1}\\
\hdotsfor{3}\\
a^j_{1\,2n}& \ldots& a^j_{2n\,2n}
\end{pmatrix}
\circ Q_{2n}^{-1},
$$
where the matrix $Q_{2n}$ has been defined above.

The unitary linear group $\mU_{2n}(K,j,Q_{2n})$ over a~ring~$K$
with an involution~$j$ is the group of matrices $A\in M_{2n}(K)$ such that
$AA^\tau =E$.

Now we formulate two theorems which were proved by E.~I.~Bunina.

\begin{theorem}
If $K_1$~and~$K_2$ are associative \textup{(}commutative\textup{)}
rings with $1/2$ and $1/3$, $j_1$~and~$j_2$ are involutions
in the rings $K_1$~and~$K_2$, respectively,
and $n,m>2$ \textup{(}$n,m>1$\textup{)},
then the unitary linear groups $\mU_{2n}(K_1,j_1,Q_{2n})$ and
$\mU_{2m}(K_2,j_2,Q_{2m})$
are elementarily equivalent if and only if
the rings $M_{2n}(K_1)$ and $M_{2m}(K_2)$ are elementarily equivalent
as rings with involutions $\tau_1$~and~$\tau_2$, respectively.
\end{theorem}

\begin{theorem}
If skewfields \textup{(}fields\textup{)}
$F_1$~and~$F_2$ have characteristic which is not equal to~$2$,
$j_1$~and~$j_2$ are involutions in skewfields \textup{(}fields\textup{)}
$F_1$~and~$F_2$, respectively, and
$n,m>2$ \textup{(}$n,m>1$\textup{)}, then the unitary linear groups
$\mU_{2n}(F_1,j_1,Q_{2n})$ and $\mU_{2m}(F_2,j_2,Q_{2m})$
are elementarily equivalent if and only if
the skewfields \textup{(}fields\textup{)}
$F_1$~and~$F_2$ are elementarily equivalent
as the skewfields \textup{(}fields\textup{)}
with involutions $j_1$~and~$j_2$, respectively.
\end{theorem}

In 2001 E.~I.~Bunina (see \cite{14,15}) studied
elementary properties of Chevalley groups over algebraically closed
fields. The class of all Chevalley groups contains many classical
groups like $\SL_n(K)$, $\PSL_n(K)$, $\SO_n(K)$, $\Spin_n(K)$,
$\PSO_n(K)$, $\mSp_{2n}(K)$, $\PSp_{2n}(K)$. Therefore,
the studied groups intersect with the groups which were considered by
A.~I.~Maltsev, but there are many other algebraic groups in this class.

The main result is the following theorem.

\begin{theorem}
Suppose that Chevalley groups $\mathcal G_1$ and $\mathcal G_2$
are constructed respectively by algebraically closed fields
$K_1$~and~$K_2$ of characteristic not equal to~$2$,
simple Lie algebras $\mathcal L_1$~and~$\mathcal L_2$, and lattices
$M\prisv L_{V_1}$ and $N\prisv L_{V_2}$.
Let $M/M_0\cong \varphi_1$ and $N/N_0\cong \varphi_2$,
where $\varphi_1$~and~$\varphi_2$ are finite groups.
Then $\mathcal G_1\equiv \mathcal G_2$ if and only if $K_1\equiv K_2$,
$\mathcal L_1\cong \mathcal L_2$, and $\varphi_1\cong \varphi_2$,
except the case where $\mathcal L_1$ and $\mathcal L_2$
have the same type~$D_{2l}$, $l\ge 3$, and
$\varphi_1\cong \varphi_2\cong \mathbb Z_2$. In this case
there exist two nonequivalent
groups such that the corresponding fields are elementarily equivalent.
\end{theorem}

In this paper we consider elementary properties of categories
of modules over rings, endomorphism rings of almost free modules
of infinite ranks over rings, and automorphism groups of almost
free modules of infinite ranks over rings.

The first section includes some basic notions from the set theory and
the model theory: definitions of first order languages,
models of a~language, deducibility, interpretability,
axioms and basic notions of the theory NBG (Neumann--Bernays--G\"odel),
which is used for all later constructions, and also some basic
notions from category theory (see~\cite{5}),
which we need in the following sections.

The second section is devoted to elementary properties and elementary
equivalence of categories of modules over rings.

In Sec.~\ref{ss2.1}, we give some additional notions about
the category $\textup{mod-}R$.

In Sec.~\ref{ss2.2}, we prove that in the category
$\textup{mod-}R$ the notion of a~progenerator object is elementary,
i.e., there exists a~formula of the first order language
of category theory with one free object variable
such that the formula
is true in the category $\textup{mod-}R$
for progenerators and only for them.

In Sec.~\ref{ss2.3}, we show that for a~given progenerator~$P$
on the semigroup $\Mor(P,P)$ we can introduce the operations of
addition and multiplication to make this semigroup isomorphic
to the ring $\Endom_R(P)$.

In Sec.~\ref{ss2.4}, we consider the case where the rings are finite
and prove the theorem that
the categories $\textup{mod-}R$ and $\textup{mod-}S$,
where the ring~$R$ is finite, are elementarily equivalent if and only
if they are Morita-equivalent.

In Sec.~\ref{ss2.5}, we remind the results of S.~Shelah from~\cite{6}
on interpretation of the set theory in a~category.

In Sec.~\ref{ss2.6}, we use the results from Sec.~\ref{ss2.5} to select
in the category $\textup{mod-}R$ for some fixed modules $X$~and~$Y$
a~set of linearly independent projectors from~$X$ on~$Y$.

In Sec.~\ref{ss2.7}, we describe
the structure $\langle \Cn, \ring\rangle$, consisting
of the class $\Cn$ of all cardinal numbers
and the ring $\ring$
with usual ring relations $+$~and~$\circ$, and
we also describe
the \emph{second-order logic} of this
structure ($L_2(\langle \Cn,\ring\rangle)$) which
allows us to use in formulas arbitrary predicate
symbols of the form
$$
P_{\lambda_1,\dots,\lambda_k}(c_1,\dots,c_k;v_1,\dots, v_n),
$$
where $\lambda_1,\dots,\lambda_k$ are fixed cardinal numbers,
$c_1,\dots,c_k$ are variables for elements from
$\lambda_1 ,\dots ,\lambda_k$,
respectively, and $v_1,\dots,v_n$ are variables for ring elements.

Further, in this section the following theorem is proved.

\begin{theorem}
Let $R$~and~$S$ be rings.
Suppose that there exists a~sentence~$\psi$
of the language $L_2(\langle \Cn,\ring\rangle)$
which is true in the ring~$R$, false in any ring
similar to~$R$, and not equivalent to it in the language
$L_2(\langle \Cn,\ring \rangle )$.
If the categories $\textup{mod-}R$ and $\textup{mod-}S$ are
elementarily equivalent, then there exists a~ring~$S'$ which is similar
to~$S$ and such that the structures
$\langle \Cn ,R\rangle $ and $\langle \Cn,S'\rangle$ are equivalent in the
logic~$L_2$.
\end{theorem}

Section~\ref{ss2.8} is devoted to the proof of the ``opposite'' theorem.

\begin{theorem}
Let $R$~and~$S$ be
arbitrary rings with unit.
If the structures $\langle \Cn,R\rangle$ and $\langle \Cn, S\rangle$
are equivalent in the second-order logic~$L_2$, then the
categories $\textup{mod-}R$ and $\textup{mod-}S$ are elementarily
equivalent.
\end{theorem}

Finally, in Sec.~\ref{ss2.9} two previous theorems imply a~theorem
which is an analogue of the Morita theorem for elementary equivalence,
as well as some useful corollaries from it.

\begin{theorem}
Let $R$~and~$S$ be rings.
Suppose that there exists
a~sentence~$\psi$ of the language $L_2(\langle \Cn,\ring\rangle)$
which is true in the ring~$R$ and is false in any ring similar
to~$R$ and not equivalent to it in the language
$L_2(\langle \Cn,\ring\rangle )$.
Then the categories $\textup{mod-}R$ and $\textup{mod-}S$ are
elementarily equivalent if and only if
there exists a~ring~$S'$ similar to the ring~$S$ and such that
the structures $\langle \Cn ,R\rangle $ and $\langle \Cn,S'\rangle$
are equivalent in the logic~$L_2$.
\end{theorem}

\begin{corollary}
For any skewfields $F_1$~and~$F_2$ the categories $\textup{mod-}F_1$ and
$\textup{mod-}F_2$ are elementarily equivalent if and only if
the structures $\langle \Cn,F_1\rangle$ and $\langle \Cn,F_2\rangle$
are equivalent in the second-order logic~$L_2$.
\end{corollary}

\begin{corollary}
For any commutative rings $R_1$~and~$R_2$
the categories $\textup{mod-}R_1$ and
$\textup{mod-}R_2$ are elementarily equivalent if and only if
the structures $\langle \Cn,R_1\rangle$ and
$\langle \Cn,R_2\rangle$ are equivalent in the second-order
logic~$L_2$.
\end{corollary}

\begin{corollary}
For arbitrary local rings $R_1$~and~$R_2$ the categories
$\textup{mod-}R_1$ and $\textup{mod-}R_2$
are elementarily equivalent if and only if
the structures $\langle \Cn,R_1\rangle$ and $\langle \Cn,R_2\rangle$
are equivalent in the second-order logic~$L_2$.
\end{corollary}

\begin{corollary}
For arbitrary integral domains $R_1$~and~$R_2$ the categories
$\textup{mod-}R_1$ and $\textup{mod-}R_2$ are elementarily
equivalent if and only if the structures
$\langle \Cn,R_1\rangle$ and $\langle \Cn,R_2\rangle$
are equivalent in the logic~$L_2$.
\end{corollary}

\begin{corollary}
For any Artinian rings $R_1$~and~$R_2$
the categories $\textup{mod-}R_1$ and
$\textup{mod-}R_2$ are elementarily equivalent if and only if there exist
rings $S_1$~and~$S_2$ such that the ring~$S_1$ is similar to
the ring~$R_1$, the ring~$S_2$ is similar to the ring~$R_2$,
and the structures $\langle \Cn, S_1\rangle$
and $\langle \Cn,S_2\rangle$ are equivalent in the
second-order logic~$L_2$.
\end{corollary}

Section~\ref{s3} is devoted to the same question for endomorphism rings
of modules of infinite ranks.

In this section, we suppose that a~ring~$R$ and an infinite
cardinal number~$\varkappa$ are such that in the ring~$R$ there exists
a~maximal ideal generated by ${\le}\,\varkappa$ elements
(for example, it is true when $\varkappa\ge |R|$ or when
the ring~$R$ is semisimple or is an integral domain).

In Sec.~\ref{ss3.1}, for every free module~$V$ of infinite rank
over a~ring we introduce some special category~$C_{M(V)}$
such that elementary equivalence of endomorphism rings
of two free modules of infinite ranks
over rings is equivalent to elementary equivalence
of the corresponding categories.

Section~\ref{ss3.2} is devoted to elementary
equivalence of categories $C_{M(V)}$. In Sec.~\ref{ss3.3}, we prove
the following main theorem and the corollaries from it.

\begin{theorem}
Let $V_1$~and~$V_2$ be free modules of infinite ranks
$\varkappa_1$~and~$\varkappa_2$ over rings $R_1$~and~$R_2$,
respectively.
Suppose that there exists a~sentence
$\psi\in \Th_2^{\varkappa_1}(\langle \varkappa_1,R_1\rangle)$
such that
$\psi\notin \Th_2^{\varkappa_1}(\langle \varkappa_1,R'\rangle)$
for every ring~$R'$
such that $R_1$ is similar to~$R'$ and
$\Th_2^{\varkappa_1}(\langle \varkappa_1,R_1\rangle)\ne
\Th_2^{\varkappa_1} (\varkappa_1, R'\rangle)$.
Then the categories $C_{M(V_1)}$ and $C_{M(V_1)}$
are elementarily equivalent if and only if
there exists a~ring~$S$ similar to the ring~$R_2$ and such that
the theories $\Th_2^{\varkappa_1}(\langle \varkappa_1 ,R_1\rangle) $
and $\Th_2^{\varkappa_2}(\langle \varkappa_2,S\rangle)$ coincide.
\end{theorem}

\begin{corollary}
Let $V_1$~and~$V_2$ be two spaces of infinite ranks
$\varkappa_1$~and~$\varkappa_2$ over arbitrary skewfields
\textup{(}integral domains\textup{)} $F_1$~and~$F_2$.
Then the rings
$\Endom_{F_1} V_1$ and $\Endom_{F_2}V_2$
are elementarily equivalent if and only if the theories
$\Th_2^{\varkappa_1}(\langle \varkappa_1,F_1\rangle)$ and
$\Th_2^{\varkappa_2}(\langle \varkappa_2,F_2\rangle)$ coincide.
\end{corollary}

\begin{corollary}
Suppose that $\varkappa_1$~and~$\varkappa_2$ are infinite
cardinal numbers, $R_1$~and~$R_2$ are commutative \textup{(}local\textup{)}
rings, and every maximal ideal of the ring~$R_1$
is generated by at most~$\varkappa_1$ elements of the ring.
Then for free modules $V_1$~and~$V_2$ of ranks
$\varkappa_1$~and~$\varkappa_2$ over the rings $R_1$~and~$R_2$, respectively,
the rings $\Endom_{R_1}V_1$ and $\Endom_{R_2}V_2$ are elementarily
equivalent if and only if the theories
$\Th_2^{\varkappa_1}(\langle \varkappa_1,R_1\rangle)$ and
$\Th_2^{\varkappa_2}(\langle \varkappa_2,R_2\rangle)$ coincide.
\end{corollary}

\begin{corollary}
Suppose that $\varkappa_1$~and~$\varkappa_2$ are infinite cardinal numbers,
$R_1$~and~$R_2$ are Artinian rings, and every maximal
ideal of the ring~$R_1$ is generated by at most~$\varkappa_1$
elements of the ring. Then for free modules $V_1$~and~$V_2$ of ranks
$\varkappa_1$~and~$\varkappa_2$ over the rings $R_1$~and~$R_2$,
respectively, the rings $\Endom_{R_1}V_1$ and $\Endom_{R_2}V_2$
are elementarily equivalent if and only if
there exist rings $S_1$~and~$S_2$ similar to the rings $R_1$~and~$R_2$,
respectively, such that the theories
$\Th_2^{\varkappa_1}(\langle \varkappa_1,S_1\rangle)$ and
$\Th_2^{\varkappa_2}(\langle \varkappa_2,S_2\rangle)$ coincide.
\end{corollary}

\begin{corollary}
For free modules $V_1$~and~$V_2$ of infinite ranks
$\varkappa_1$~and~$\varkappa_2$ over semisimple rings
$R_1$~and~$R_2$, respectively, the rings
$\Endom_{R_1}(V_1)$ and $\Endom_{R_2}(V_2)$
are elementarily equivalent if and only if
there exist rings $S_1$~and~$S_2$ similar to the rings
$R_1$~and~$R_2$, respectively,
such that the theories
$\Th_2^{\varkappa_1}(\langle \varkappa_1,S_1\rangle)$ and
$\Th_2^{\varkappa_2}(\langle \varkappa_2,S_2\rangle)$ coincide.
\end{corollary}

In Sec.~\ref{s4}, we consider projective spaces of modules of
infinite ranks.

In Sec.~\ref{ss4.1}, we describe the language of projective spaces
and basic notions which can be expressed in this language.

In Sec.~\ref{ss4.2}, we show
how in a~projective space of a~module of infinite rank
one can interpret
a~ring that is isomorphic to the ring $\Endom_RP$ for some
progenerator~$P$.

In Sec.~\ref{ss4.3}, we show how
to interpret the ring $\Endom_RV$
in a~projective space of the module~$V$.

Finally, in this section, we prove the following theorem.

\begin{theorem}
For free modules $V_1$~and~$V_2$ of infinite ranks
over arbitrary rings $R_1$~and~$R_2$, respectively,
elementary equivalence of the lattices of submodules $P(V_1)$ and
$P(V_2)$ implies elementary equivalence of the endomorphism rings
$\Endom_{R_1}(V_1)$ and $\Endom_{R_2}(V_2)$.
\end{theorem}

In Sec.~\ref{ss4.4}, we prove the ``inverse'' theorem.

\begin{theorem}
Suppose that $V_1$~and~$V_2$ are free modules of infinite ranks
$\varkappa_1$~and~$\varkappa_2$ over rings $R_1$~and~$R_2$, respectively,
and every submodule of the module $V_1$~\textup{(}$V_2$\textup{)}
has at most~$\varkappa_1$~\textup{(}$\varkappa_2$\textup{)}
generating elements
\textup{(}for example,
this is true if $\varkappa_1\ge |R_1|$
and $\varkappa_2\ge |R_2|$, or if $R_1$ and~$R_2$ are semisimple rings
or integral domains\textup{)}.
Then
$\Endom_{R_1}(V_1)\equiv \Endom_{R_2}(V_2)$
implies $P(V_1)\equiv P(V_2)$.
\end{theorem}

In Sec.~\ref{s5}, we consider automorphism groups of modules
of infinite ranks over rings.

In Sec.~\ref{ss5.1}, as in~\cite{1}, we prove that if rings
$R$~and~$S$ with $1/2$ do not contain any central idempotents that are
not equal to $0$~and~$1$,
$V$~and~$V'$ are free modules of infinite ranks over the rings $R$~and~$S$,
respectively, then the groups $\Aut_R(V)$ and $\Aut_S(V')$
are isomorphic if and only if $\Endom_R(V)\cong \Endom_S(V')$.

In Sec.~\ref{ss5.2}, all results of Sec.~\ref{ss5.1} are proved
for elementary equivalences.
We do this with the help of ultrapowers, like in the paper~\cite{2}
of C.~I.~Beidar and A.~V.~Mikhalev.
We prove the following theorem.

\begin{theorem}
Suppose that rings $R$ and $S$ contain $1/2$ and do not contain any
central idempotents which are not equal to $1$~and~$0$.
Suppose that $V$ and~$V'$ are
free modules of infinite ranks over the rings $R$~and~$S$, respectively.
Then the groups $\Aut_R(V)$ and $\Aut_S(V')$
are elementarily equivalent if and only if
the rings $\Endom_R(V)$ and $\Endom_S(V')$ are elementarily equivalent.
\end{theorem}

In Sec.~\ref{ss5.3}, we assume that the cardinal number~$\varkappa_1$
is such that the ring~$R_1$ has a~maximal ideal generated by
at most $\varkappa_1$ elements.

\begin{theorem}
Suppose that rings $R_1$~and~$R_2$ contain~$1/2$ and do not contain
any central idempotents which are not equal to $1$~or~$0$. Let
$V_1$~and~$V_2$ be free modules of infinite ranks
$\varkappa_1$~and~$\varkappa_2$ over the rings $R_1$~and~$R_2$,
respectively, and let
$\psi \in \Th_2^{\varkappa_1}(\langle \varkappa_1,R_1\rangle)$
be such that
$\psi\notin \Th_2^{\varkappa_1} (\langle \varkappa_1, R'\rangle )$
for any ring~$R'$ such that $R'$ is similar to~$R_1$ and
$\Th_2^{\varkappa_1}(\langle \varkappa_1,R_1\rangle )\ne
\Th_2^{\varkappa_1}(\varkappa_1, R'\rangle)$.
Then the groups $\Aut_{R_1}(V_1)$ and $\Aut_{R_2}(V_2)$ are elementarily
equivalent if and only if there exists a~ring~$S$
similar to the ring~$R_2$ and such that the theories
$\Th_2^{\varkappa_1}(\langle \varkappa_1,R_1\rangle)$ and
$\Th_2^{\varkappa_2}(\langle \varkappa_2,S\rangle)$ coincide.
\end{theorem}

\begin{corollary}
For free modules $V_1$~and~$V_2$ of infinite ranks
$\varkappa_1$~and~$\varkappa_2$ over skewfields
\textup{(}integral domains, commutative or local
rings without central idempotents not equal to $1$~or~$0$\textup{)}
$F_1$~and~$F_2$ with~$1/2$, respectively,
the groups $\Aut_{F_1}(V_1)$ and $\Aut_{F_2}(V_2)$
are elementarily equivalent if and only if the theories
$\Th_2^{\varkappa_1}(\langle \varkappa_1, F_1\rangle)$
and $\Th_2^{\varkappa_2}(\langle \varkappa_2, F_2\rangle)$
coincide.
\end{corollary}

\begin{corollary}
For free modules $V_1$~and~$V_2$ of infinite ranks
$\varkappa_1$~and~$\varkappa_2$ over Artinian
rings $R_1$~and~$R_2$
with~$1/2$ without central idempotents not equal
to $1$~or~$0$, respectively,
the groups $\Aut_{R_1}(V_1)$ and $\Aut_{R_2}(V_2)$
are elementarily equivalent
if and only if there exist rings $S_1$~and~$S_2$ such that
the ring~$R_1$ is similar to the ring~$S_1$, the ring~$R_2$
is similar to the ring~$S_2$, and the theories
$\Th_2^{\varkappa_1}(\langle \varkappa_1, S_1\rangle)$
and $\Th_2^{\varkappa_2}(\langle \varkappa_2, S_2\rangle)$
coincide.
\end{corollary}

\section{Basic Notions from the Set Theory, Model Theory, and
Category Theory}
\subsection{First Order Languages}\label{ss1.1}
The \emph{first order language}~$\mathcal L$ is some set of symbols.
This set consists of
\begin{itemize}
\item[]
the blank symbol;
\item[]
the parentheses $(,)$;
\item[]
the connectives $\Rightarrow$ (``implies'') and $\neg$ (``not'');
\item[]
the quantifier $\forall$ (for all);
\item[]
the equality symbol $=$;
\item[]
a~countable set of variables~$v_i$ ($i\ge 0$);
\item[]
a~nonempty countable set of predicate symbols~$P_i^n$ ($n\ge 1$, $i\ge 0$);
\item[]
a~countable set of function symbols~$F_i^n$ ($n\ge 1$, $i\ge 0$);
\item[]
a~countable set of constant symbols~$c_i$ ($i\ge 0$).
\end{itemize}

Some symbol-strings constructed from these symbols of the first order
language~$\mathcal L$ are called
\emph{terms} and \emph{formulas} of this language.

\emph{Terms} are defined in the following way:
\begin{enumerate}
\item
a~variable is a~term;
\item
a~constant symbol is a~term;
\item
if $F_i^n$ is some function symbol,
$t_0,\dots,t_{n-1}$ are terms, then $F_i^n(t_0,\dots,t_{n-1})$ is
a~term;
\item
a~symbol-string is a~term if and only if this follows from the rules
(1)--(3).
\end{enumerate}

If $P_i^n$ is some predicate symbol and
$t_0,\dots,t_{n-1}$ are terms,
then the symbol-string $(P_i^n(t_0,\dots,t_{n-1}))$ is called an
\emph{elementary formula}.

\emph{Formulas} of the language~$\mathcal L$
are defined in the following way:
\begin{enumerate}
\item
every elementary formula is a~formula;
\item
if $\varphi$~and~$\psi$ are formulas and
$v$~is a~variable, then each of the symbol-strings
$(\neg \varphi)$, $(\varphi\Rightarrow \psi)$, and
$(\forall v \losp \varphi)$ is a~formula;
\item
a~symbol-string is a~formula
if and only if this follows from the rules (1)~and~(2).
\end{enumerate}

Let us introduce the following abbreviations:
\begin{itemize}
\item[]
$(\varphi\land \psi)$ stands for $(\neg (\varphi\Rightarrow (\neg \psi)))$;
\item[]
$(\varphi\vee \psi)$ stands for $((\neg \varphi)\Rightarrow \psi)$;
\item[]
$(\varphi\equiv \psi)$ stands for
$((\varphi\Rightarrow \psi)\land (\psi\Rightarrow \varphi))$;
\item[]
$(\exists v\losp \varphi)$ is an abbreviation for
$(\neg (\forall v\losp (\neg \varphi)))$;
\item[]
$\varphi_1\vee \varphi_2\vee\dots \vee\varphi_n$ stands for
$(\varphi_1\vee (\varphi_2\vee \dots \vee \varphi_n))$;
\item[]
$\varphi_1\wedge \varphi_2\wedge \dots \wedge \varphi_n$ stands for
$(\varphi_1\wedge (\varphi_2\wedge \dots \wedge \varphi_n))$;
\item[]
$(\forall x_1x_2\dots x_n)\varphi$ stands for
$(\forall x_1)(\forall x_2)\dots (\forall x_n)\varphi$;
\item[]
$(\exists x_1x_2\dots x_n)\varphi)$ stands for
$(\exists x_1)(\exists x_2)\dots (\exists x_n)\varphi$.
\end{itemize}

We introduce the notion of \emph{free} and \emph{bound} occurrences
of a~variable in a~formula. An occurrence of a~variable~$v$ in a~given
formula is called \emph{bound}
if $v$ is either the variable of
a~quantifier prefix~$\forall v$ occurring in this formula
or is under the action of
a~quantifier prefix~$\forall v$ occurring in this formula;
otherwise an occurrence of a~variable in
a~given formula is called \emph{free}.
Thus, one variable can have free and bound occurrences
in the same formula.
A~variable is called \emph{free}
(\emph{bound}) \emph{variable} in a~given formula if there exist free
(bound) occurrences of this variable in this formula, i.e., a~variable
can at the same time be free and bound in one formula.

A~\emph{sentence} is a~formula with no free variables.

If $\zeta$ is a~term or a~formula, $\theta$ is a~term, and
$v$ is a~variable, then $\zeta(v\|\theta)$ denotes
the symbol-string obtained by replacing every free
occurrence of the variable~$v$ in the symbol-string~$\zeta$
by the symbol-string~$\theta$.

A substitution $v\|\theta$ in $\zeta$ is called \emph{admissible} if
for every free occurrence of a~variable~$w$ in the symbol-string~$\theta$
every free occurrence~$v$ in~$\zeta$ is not a~free occurrence
in some formula~$\psi$ occurring in some formula
$\forall w\losp \psi(w)$ or $\exists w\losp \psi(w)$
that occurs in the
symbol-string~$\zeta$.

In the sequel, if a~substitution $v\|\theta$ in~$\zeta$ is admissible,
then along with $\zeta(v\|\theta)$ we shall write $\zeta(\theta)$.

If $\zeta$ is a~term or a~formula, $\theta$~is a~term,
and $v$~is a~variable
such that the substitution $v\|\theta$ in~$\zeta$ is admissible,
then the substitution $\zeta(v\|\theta)$ is a~term or a~formula,
respectively.

Every free occurrence of some variable~$u$
(except~$v$) in a~symbol-string~$\zeta$ and every free occurrence of some
variable~$w$ in a~symbol-string~$\theta$
are free occurrences of these variables
in a~symbol-string $\zeta(v\|\theta)$
(provided that the variable~$v$ is free in~$\zeta$).

A~symbol-string~$\gamma$, equipped with some rule, is called
a~\emph{formula scheme of a~language~$\mathcal L$} if
\begin{enumerate}
\item
this rule marks some letters (in particular, free and bound
variables) occurring in~$\gamma$;
\item
this rule determines the necessary substitution of these marked letters
in~$\gamma$ by some terms (in particular, variables);
\item
after every such substitution in~$\gamma$ some propositional
formula~$\varphi$ of the language~$\mathcal L$ is obtained.
\end{enumerate}
Each such propositional formula~$\varphi$ is called a~\emph{formula,
obtained by the application of the formula scheme~$\gamma$}.

A~text~$\Gamma$ consisting of symbol-strings separated by blank-symbols
is called an \emph{axiom text} if every symbol-string~$\gamma$ occurring
in~$\Gamma$ is either a~formula or a~formula scheme
of the language~$\mathcal L$.
If $\gamma$ is a~formula, then $\gamma$ is called an
\emph{explicit axiom of the language~$\mathcal L$}.
If $\gamma$ is a~formula scheme, then it is called an
\emph{axiom scheme of the language~$\mathcal L$}.
Every formula obtained by the application
of the axiom scheme~$\gamma$ is called an \emph{implicit axiom of the
language~$\mathcal L$}.

We need \emph{logical axioms} and \emph{rules of deduction} to
construct a~formal system.

Logical axiom schemes of any first order language are cited below.
\begin{mVentry}{\textbf{LAS12.}}
\item[\textbf{LAS1.}]
$\varphi\Rightarrow (\psi\Rightarrow \varphi)$.
\item[\textbf{LAS2.}]
$(\varphi\Rightarrow (\psi\Rightarrow \chi))\Rightarrow
((\varphi \Rightarrow \psi)\Rightarrow (\varphi\Rightarrow \chi))$.
\item[\textbf{LAS3.}]
$(\varphi\land \psi)\Rightarrow \varphi$.
\item[\textbf{LAS4.}]
$(\varphi\land \psi)\Rightarrow \psi$.
\item[\textbf{LAS5.}]
$\varphi \Rightarrow (\psi\Rightarrow (\varphi \land \psi))$.
\item[\textbf{LAS6.}]
$\varphi \Rightarrow (\varphi\lor \psi)$.
\item[\textbf{LAS7.}]
$\psi \Rightarrow (\varphi\lor \psi)$.
\item[\textbf{LAS8.}]
$(\varphi\Rightarrow \chi)\Rightarrow
((\psi \Rightarrow \chi)\Rightarrow ((\varphi \lor \psi)\Rightarrow \chi))$.
\item[\textbf{LAS9.}]
$(\varphi \Rightarrow \psi )\Rightarrow
((\varphi \Rightarrow (\neg \psi))\Rightarrow (\neg \varphi))$.
\item[\textbf{LAS10.}]
$(\neg (\neg \varphi)) \Rightarrow\varphi$.
\item[\textbf{LAS11.}]
$(\forall v \varphi)\Rightarrow \varphi(v\|\theta)$
if $v$ is a~variable and
$\theta$ is a~term such that the substitution $v\|\theta$ in~$\varphi$
is admissible.
\item[\textbf{LAS12.}]
$\varphi(v\| \theta)\Rightarrow (\exists v \varphi)$ in the same
conditions as in \textbf{LAS11}.
\item[\textbf{LAS13.}]
$(\forall v (\psi\Rightarrow \varphi(v)))\Rightarrow
(\psi\Rightarrow (\forall v \varphi))$
if the variable~$v$ is not free in $\psi$.
\item[\textbf{LAS14.}]
$(\forall v (\varphi(v)\Rightarrow \psi))\Rightarrow
((\exists v \varphi)\Rightarrow \psi)$
if the variable~$v$ is not free in $\psi$.
\end{mVentry}

Rules of deduction are the following.
\begin{description}
\item[\hspace*{-\parindent}the rule of implication (modus ponens or MP):]
from $\varphi$ and $\varphi\Rightarrow \psi$ it follows that $\psi$;
\item[\hspace*{-\parindent}the rule of generalization (Gen):]
from $\varphi$ it follows that $(\forall x)(\varphi)$.
\end{description}

Let $\Sigma$ be a~totality of formulas and $\psi$ be a~formula of
the language~$\mathcal L$.
A~sequence
$f\equiv (\varphi_i|i\in n+1)\equiv (\varphi_0,\dots,\varphi_n)$
of formulas of the language~$\mathcal L$ is called
a~\emph{deduction of the formula~$\psi$ from the totality~$\Sigma$}
if $\varphi_n=\psi$ and for any $0\le i\le n$
one of the following conditions is fulfilled:
\begin{enumerate}
\item
$\varphi_i$ belongs to~$\Sigma$ or is a~logical axiom;
\item
there exist $0\le k< j< i$ such that $\varphi_j$ is
$(\varphi_k \Rightarrow \varphi_i)$, i.e.,
$\varphi_i$~is obtained from $\varphi_k$ and
$\varphi_k\Rightarrow \varphi_i$ by the rule of implication MP;
\item
there exists $0\le j< i$ such that $\varphi_i$ is $\forall x\losp \varphi_j$,
where $x$ is not a~free variable of any formula from~$\Sigma$, i.e.,
$\varphi_i$~is obtained from~$\varphi_j$ by the rule of generalization
Gen with the given \emph{structural requirement}.
\end{enumerate}

Denote this deduction either by
$f\equiv (\varphi_0,\dots,\varphi_n)\colon \Sigma\vdash \psi$,
or by $(\varphi_0,\dots,\varphi_n)\colon \Sigma\vdash \psi$,
or by $f\colon \Sigma\vdash \psi$.

If there exists a~deduction $f\colon \Sigma\vdash \psi$,
then the formula~$\psi$ is called
\emph{deducible in the language~$\mathcal L$ from the set~$\Sigma$},
and the deduction~$f$ is called a~\emph{proof of the formula~$\psi$}.

A~(\emph{first order}) \emph{theory}~$T$ in the language~$\mathcal L$ is
a~set of sentences of the language~$\mathcal L$.
A~\emph{set of axioms} of a~theory~$T$ is any set of sentences,
which has the same corollaries as~$T$.

Now we introduce axioms and basic notions of the set theory
NBG (von Neumann--Bernays--G\"odel) (see~\cite{7}),
which is a~first order theory. We shall use it for all our
constructions.

\subsection{Axioms and Basic Notions of the Theory NBG}\label{ss1.2}
The set theory NBG (see~\cite{7})
has one predicate symbol~$P$, which denotes a~2-place relation,
no function symbols, and no constant symbols.
We shall use Latin letters $X$, $Y$, and~$Z$
with subscripts and apostrophes as variables of this system.
We also introduce the abbreviations $X\in Y$ for $P(X,Y)$
and $X\notin Y$ for $\neg P(X,Y)$. The sign~$\in$ can be interpreted
as the symbol of belonging.

The formula $X=Y$ ($X$~is equal to~$Y$)
is an abbreviation for the formula
$\forall Z\losp (Z\in X\Leftrightarrow Z\in Y)$, i.e.,
two objects are equal if they consist of the same elements.

The formula $X\subseteq Y$ is an abbreviation
for the formula $\forall Z\losp (Z\in X\Rightarrow Z\in Y)$
(\emph{inclusion}),
$X\subset Y$ is an abbreviation for $X\subseteq Y \logic\wedge X\ne Y$
(\emph{proper inclusion}).

From these definitions we can easily get the following proposition.

\begin{proposition}
\begin{enumerate}
\renewcommand{\theenumi}{\alph{enumi}}
\item
$\vdash X=Y\Leftrightarrow (X\subseteq Y \logic\wedge Y\subseteq X)$\textup;
\item
$\vdash X=X$\textup;
\item
$\vdash X=Y\Rightarrow Y=X$\textup;
\item
$\vdash X=Y\Rightarrow (Y=Z\Rightarrow X=Z)$\textup;
\item
$\vdash X=Y\Rightarrow (Z\in X\Rightarrow Z\in Y)$.
\end{enumerate}
\end{proposition}

Objects of the theory NBG are called \emph{classes}.
A class is called a~\emph{set} if it is an element of some class.
A class which is not a~set is called a~\emph{proper class}.
We introduce small Latin letters $x$, $y$, and~$z$ with subscripts
as special variables bounded by sets.
This means that the formula $\forall x\losp A(x)$ is an abbreviation
for $\forall X\losp (\text{$X$ is a~set}\Rightarrow A(X))$, and it has
the sense ``$A$~is true for all sets'',
and $\exists x\losp A(x)$ is an abbreviation for
$\exists X\losp (\text{$X$~is a~set} \logic\wedge A(X))$,
and it has the sense ``$A$~is true for some set.''

\begin{description}
\item[\hspace*{-\parindent}A1
{\mdseries (\emph{the extensionality axiom}).}]
$X=Y\Rightarrow (X\in Z\Leftrightarrow Y\in Z)$.
\item[\hspace*{-\parindent}A2
{\mdseries (\emph{the pair axiom}).}]
$\forall x\losp \forall y\losp \exists z\losp \forall u\losp
(u\in z\Leftrightarrow u=x \logic\vee u=y)$,
i.e., for all sets $x$~and~$y$
there exists a~set~$z$ such that $x$~and~$y$ are the only
elements of~$z$.
\item[\hspace*{-\parindent}A3
{\mdseries (\emph{the empty set axiom}).}]
$\exists x\losp \forall y\losp \neg (y\in x)$,
i.e., there exists a~set which does not contain any elements.
\end{description}

Axioms \textbf{A1}~and~\textbf{A3} imply that this set is unique,
i.e., we can introduce a~constant symbol~$\varnothing$
(or~$0$), with the condition $\forall y\losp (y\notin \varnothing)$.

Also we can introduce a~new function symbol
$f(x,y)$ for the pair, and write it in the form
$\{ x,y\}$. We can even define a~pair $\{ X,Y\}$
for arbitrary classes $X$~and~$Y$, setting $\{ X,Y\}=0$ if one of
the classes $X$,~$Y$ is not a~set.
Further, set $\{ X\}=\{ X,X\}$.
The class $\langle X,Y\rangle \equiv\{ \{X\},\{ X,Y\}\}$ is called
the \emph{ordered pair} of classes $X$~and~$Y$.
Similarly we can introduce \emph{ordered triplets, quadruplets}
and so on.

\begin{description}
\item[\hspace*{-\parindent}AS4
{\mdseries (\emph{the axiom scheme of existence of classes}).}]
Let
$$
\varphi(X_1,\dots,X_n,Y_1,\dots,Y_m)
$$
be a~formula. We shall call this formula
\emph{predicative} if only variables for sets are bound in it
(i.e.,
if it can be transferred to this form with the help of abbreviations).
For every predicative formula $\varphi(X_1,\dots,X_n,Y_1,\dots,Y_m)$
$$
\exists Z\losp \forall x_1\dots \forall x_n\losp
(\langle x_1,\dots,x_n\rangle \in Z
\Leftrightarrow
\varphi(x_1,\dots,x_n,Y_1,\dots,Y_m)).
$$
\end{description}

The class~$Z$ which exists by the axiom scheme~\textbf{AS4}
will be denoted by
$$
\{ x_1,\dots,x_n\mid \varphi(x_1,\dots,x_n,Y_1,\dots,Y_m)\}.
$$

Now, by the axiom scheme \textbf{AS4}, we can define for arbitrary
classes $X$~and~$Y$ the following derivative classes:
\begin{itemize}
\item[]
$X\cap Y\equiv \{ u\mid u\in X \logic\land u\in Y\}$
(\emph{the intersection of classes $X$~and~$Y$});
\item[]
$X\cup Y\equiv \{ u\mid u\in X \logic\lor u\in Y\}$
(\emph{the union of classes $X$~and~$Y$});
\item[]
$\bar X\equiv \{ u\mid u\notin X\}$
(\emph{the addition to a~class~$X$});
\item[]
$V\equiv \{ u\mid u=u\}$ (\emph{the universal class});
\item[]
$X\setminus Y\equiv \{ u\mid u\in X \logic\land u\notin Y\}$
(\emph{the difference of classes $X$~and~$Y$});
\item[]
$\Dom(X)\equiv \{ u\mid \exists v\losp (\langle u,v\rangle \in X)\}$
(\emph{the domain of a~class~$X$});
\item[]
$X\times Y\equiv \{ u\mid \exists x\losp \exists y\losp
(u=\langle x,y\rangle \logic\land x\in X \logic\land y\in Y)\}$
(\emph{the Cartesian product of classes $X$~and~$Y$});
\item[]
$\mathcal P(X)\equiv \{ u\mid u\subseteq X\}$
(\emph{the class of all subsets of a~class~$X$});
\item[]
$\ucup X\equiv \{ u\mid \exists v\losp (u\in v \logic\land v\in X)\}$
(\emph{the union of all elements of a~class~$X$}).
\end{itemize}

Introduce now other axioms.

\begin{description}
\item[\hspace*{-\parindent}A5 {\mdseries (\emph{the union axiom}).}]
$\forall x\losp \exists y\losp \forall u\losp
(u\in y\Leftrightarrow \exists v\losp (u\in v \logic\land v\in x))$.
\item[\hspace*{-\parindent}A6 {\mdseries (\emph{the power set axiom}).}]
$\forall x\losp \exists y\losp \forall u\losp
(u\in y\Leftrightarrow u\subseteq x)$.
\item[\hspace*{-\parindent}A7 {\mdseries (\emph{the separation axiom}).}]
$\forall x\losp \forall Y\losp \exists z\losp \forall u\losp
(u\in z\Leftrightarrow u\in x \logic\wedge u\in Y)$.
\end{description}

Denote the class $X\times X$ by $X^2$,
the class $X\times X\times X$ by~$X^3$
and so on.
Denote the formula
$\forall x\losp \exists y\losp \forall z\losp
(\langle x,y\rangle \in X \logic\wedge \langle x,z\rangle \in X
\Rightarrow y=z)$
by
$\Un(X)$.

\begin{description}
\item[\hspace*{-\parindent}A8 {\mdseries (\emph{the replacement axiom}).}]
$\forall X\losp \forall x\losp
(\Un(X)\Rightarrow \exists y\losp \forall u\losp
(u\in y \Leftrightarrow \exists v\losp
(\langle v,u\rangle \in X \logic\wedge v\in x)))$.
\item[\hspace*{-\parindent}A9 {\mdseries (\emph{the infinity axiom}).}]
$\exists x\losp (0\in x \logic\wedge
\forall u\losp (u\in x\Rightarrow u\cup \{ u\}\in x))$.
It is clear that for such a~set~$x$
we have
$\{ 0\}\in x$, $\{ 0,\{ 0\}\}\in x$,
$\{ 0,\{ 0\},\{ 0,\{ 0\}\}\}\in x$,\ldots{}
If we now set $1\prisv \{0\}$, $2\prisv \{ 0,1\}$,\ldots,
$n\prisv \{ 0,1,\dots,n-1\}$,
then for every integer $n\ge 0$ the condition $n\in x$ is fulfilled
and $0\ne 1$, $0\ne 2$, $1\ne 2$,\ldots
\item[\hspace*{-\parindent}A10 {\mdseries (\emph{the regularity axiom}).}]
$\forall X\losp (X\ne \varnothing\Rightarrow
\exists x\in X\losp (x\cap X = \varnothing))$.
\item[\hspace*{-\parindent}A11 {\mdseries (\emph{the axiom of choice AC}).}]
For every set~$x$ there exists
a~mapping~$f$ such that for every nonempty subset
$y\subseteq x$ we have $f(y)\in y$
(this mapping is called a~\emph{choice} mapping for~$x$).
\end{description}

The list of axioms of the theory NBG is finished.

A~class~$P$ is called \emph{ordered by a~binary relation~$\le$ on~$P$}
if the following conditions hold
\begin{enumerate}
\item
$\forall p\in P\losp (p\le p) $;
\item
$\forall p,q\in P\losp (p\le q \logic\land q\le p\Rightarrow p=q)$;
\item
$\forall p,q,r\in P\losp (p\le q \logic\land q\le r\Rightarrow p\le r)$.
\end{enumerate}

If, in addition,
\begin{enumerate}
\setcounter{enumi}{3}
\item
$\forall p,q\in P\losp (p\le q \logic\lor q\le p)$,
\end{enumerate}
then the relation~$\le$ is called a~\emph{linear order}
on the class~$P$.

An ordered class~$P$ is called \emph{well-ordered} if
\begin{enumerate}
\setcounter{enumi}{4}
\item
$\forall q\losp (\varnothing \ne q\subseteq P\Rightarrow
\exists x\in q\losp (\forall y\in q\losp (x\le y)))$,
i.e., every nonempty subset of the class~$P$
has the smallest element.
\end{enumerate}

If a~class~$P$ is ordered by a~relation~$\le$ and $A$ is
a~nonempty subclass of the class~$P$, then an element $p\in P$
is called the \emph{least upper bound} or the
\emph{supremum of the subclass}~$A$ if
$$
\forall x\in A\losp (x\le p) \logic\land
\forall y\in P\losp ((\forall x'\in A\losp (x'\le y))\Rightarrow p\le y).
$$
This formula is denoted by $p=\sup A$.

A~class~$S$ is called \emph{transitive} if
$\forall x\losp (x\in S\Rightarrow x\subseteq S)$.

A~class (a~set)~$S$ is called an \emph{ordinal}
(an \emph{ordinal number}) if $S$
is transitive and well-ordered by the relation ${\in}\cup {=}$ on~$S$.
The property of a~class~$S$ to be an ordinal will be denoted by
$\On(S)$.

Ordinal numbers are usually denoted by Greek letters
$\alpha$, $\beta$, $\gamma$, and so on. The class of all
ordinal numbers is denoted by~$\On$.
The natural ordering of the class of ordinal numbers is the relation
${\alpha\le \beta} \prisv \allowbreak
\alpha=\beta \logic\lor \alpha\in \beta$.
The class $\On$ is transitive and linearly ordered by the relation
${\in} \cup {=}$.

There are some simple assertions about ordinal numbers:
\begin{enumerate}
\item
if $\alpha$ is an ordinal number, $a$ is a~set, and $a\in \alpha$,
then $a$ is an ordinal number;
\item
$\alpha= \{ \beta\mid \beta\in \alpha\}$ for every ordinal
number~$\alpha$;
\item
$\alpha+1\equiv \alpha\cup \{ \alpha\}$
is the smallest ordinal number that is greater
than~$\alpha$;
\item
every nonempty set of ordinal numbers has the smallest element.
\end{enumerate}

Therefore the ordered class $\On$ is well-ordered.
Thus $\On$ is an ordinal.

\begin{lemma}\label{4}\label{l9}
Let $A$ be a~nonempty subclass of the class~$\On$.
Then $A$ has the smallest element.
\end{lemma}

\begin{lemma}\label{l10}
If $a$ is a~nonempty set of ordinal numbers, then the following
statements hold\textup:
\begin{enumerate}
\item
the class $\ucup a$ is an ordinal number\textup;
\item
$\ucup a=\sup a$ in the ordered class~$\On$.
\end{enumerate}
\end{lemma}

An ordinal number $\alpha$ is called a~\emph{successor} if
$\alpha=\beta+1$ for some ordinal number~$\beta$.
This unique number~$\beta$ will be denoted by $\alpha-1$.
In the opposite case $\alpha$ is called a~\emph{limit ordinal number}.

\begin{lemma}\label{l11}
An ordinal number~$\alpha$ is a~limit ordinal number if and only if
$\alpha =\sup \alpha$.
\end{lemma}

The smallest (in the class $\On$)
nonzero limit ordinal is denoted by~$\omega$.
Ordinals which are smaller than~$\omega$
are called \emph{natural numbers}.

The classes~$F$ which are mappings with domains equal to~$\omega$
are called \emph{infinite sequences}.
Mappings
with domains equal to $n\in \omega$ are called
\emph{finite sequences}.

\begin{theorem}[the principle of transfinite induction]\label{t1}
Let $C$ be a~class of ordinal numbers such that the following
statements hold\textup:
\begin{enumerate}
\item
$\varnothing \in C$\textup;
\item
$\alpha\in C \Rightarrow \alpha+1\in C$\textup;
\item
$\text{\textup($\alpha$ is a~limit ordinal number} \logic\land
\alpha\subset C) \Rightarrow \alpha\in C$.
\end{enumerate}
Then $C=\On$.
\end{theorem}

Sets $a$~and~$b$ are called
\emph{equivalent} (notation: $|a|=|b|$ or $a\sim b$) if there
exists a~bijective mapping $u\colon a\to b$.

An ordinal number~$\alpha$ is called a~\emph{cardinal} if
for every ordinal number~$\beta $ the conditions
$\beta \leq \alpha$ and $|\beta|=|\alpha|$
imply $\beta =\alpha$. The class of all cardinal numbers
will be denoted by~$\Cn$. The class $\Cn$
with the order induced from the class $\On$
is well-ordered.

The axiom of choice implies the following lemma.

\begin{lemma}\label{l13}
For every set~$a$ there exists an ordinal number~$\alpha$
such that $|a|=|\alpha|$.
\end{lemma}

Now for a~set~$a$ consider the class
$\{ x\mid x\in \On \logic\land x\sim a\}$. By Lemma~\ref{l13},
this class is not empty and therefore it contains
the smallest element~$\alpha$.
It is clear that $\alpha$ is a~cardinal number.
Further, this class contains only one cardinal number~$\alpha$.
This number~$\alpha$ is called the \emph{cardinality of the set}~$a$
(it is denoted by $|a|$ or $\card a$).
Two sets having the same cardinality are equivalent.
A~set of cardinality~$\omega$ is called \emph{denumerable}.
Sets of cardinality $n\in \omega$ are called
\emph{finite}. A~set is called \emph{countable}
if it is finite or denumerable. A~set is called
\emph{infinite} if it is not finite.
A~set is called \emph{uncountable} if it is not countable.

Note that if $\varkappa$ is an infinite cardinal number, then
$\varkappa$ is a~limit ordinal number.

As for ordinal numbers, we use Greek letters for cardinal numbers:
$\xi$th infinite cardinal number will be denoted by~$\omega_\xi$
(i.e., the cardinal number~$\omega$ will be denoted also by~$\omega_0$).

Let $\alpha$ be an ordinal. A~\emph{confinality} of~$\alpha$
is the ordinal number $\cf(\alpha)$ which is equal to the smallest
ordinal number~$\beta$ for which there exists a~function~$f$
from~$\beta$ into~$\alpha$ such that $\sup f[\beta]=\alpha$.

A~cardinal~$\varkappa$ is called \emph{regular} if
$\cf (\varkappa)=\varkappa$, i.e., for every ordinal number~$\beta$
for which there exists a~function $f\colon \beta\to \varkappa$ such that
$\ucup \rng  f=\varkappa$ the inequality $\varkappa \le \beta$ holds,
where $\ucup \rng  f=\varkappa$ means that for every $y\in \varkappa$
there exists $x\in \beta$ such that $y< f(x)$.

A~cardinal $\varkappa >\omega$ is called
(\emph{strongly}) \emph{inaccessible}
if $\varkappa$ is regular and $\card \mathcal P(\lambda) < \varkappa$
for all ordinal numbers $\lambda< \varkappa$.

\subsection{Models, Deducibility, and Elementary Equivalence}\label{ss1.3}
Since we now suppose that all our constructions are made in the theory
NBG,
it follows that in the definition of deduction of a~formula~$\psi$ from
a~totality~$\Sigma$
we can change condition~(1) (``$\varphi_i$~belongs
to~$\Sigma$ or is a~logical axiom'')
to ``$\varphi_i$~belongs to~$\Sigma$ or is a~logical axiom
or is a~proper axiom of the theory NBG.''

Let in the theory NBG some object~$A$ be selected.
This selected object~$A$
is called a~\emph{universe}
if in the theory NBG for all $n\ge 1$ the notions
of $n$-finite sequence
$(x_i\in A\mid i\in n)$ of elements of the object~$A$,
$n$-place relation
$R\subset A^n$, $n$-place operation $O\colon A^n\to A$,
and infinite sequence $x_0,\dots,x_q,\dots$
of elements of the object~$A$ are defined.

A~\emph{model of a~first order language~$\mathcal L$
equipped with the universe~$A$}
is a~pair~$\mathcal U$
consisting of the object~$A$ and some correspondence~$I$
that assigns to every predicate symbol~$P_i^n$ some
$n$-place relation in~$A$, to every function symbol~$F_i^n$ some
$n$-place operation in~$A$, and to every constant symbol~$c_i$
some element of~$A$.

Let $s$ be an infinite sequence $x_0,\dots,x_q,\dots$
of elements of the object~$A$.

Define the \emph{value of a~term~$t$ of the language~$\mathcal L$
on the sequence~$s$ in the model~$\mathcal U$}
(notation: $t_{\mathcal U}[s]$) by induction in the following way:
\begin{itemize}
\item[--]
if $t\equiv v_i$, then $t_{\mathcal U}[s]\equiv x_i$;
\item[--]
if $t\equiv c_i$, then $t_{\mathcal U}[s]\equiv I(c_i)$;
\item[--]
if $t\equiv F^n_i(t_1,\dots,t_n)$, where $F$ is
a~function symbol and $t_1,\dots,t_n$ are terms, then
$t_{\mathcal U}[s]\equiv
I(F^n_i)({t_1}_{\mathcal U}[s],\dots,{t_{n}}_{\mathcal U} [s])$.
\end{itemize}

Define the \emph{translation of a~formula $\varphi$ on the sequence~$s$
in the model~$\mathcal U$}
(notation: $\mathcal U\vDash \varphi[s]$) by induction in the following way:
\begin{itemize}
\item[--]
if $\varphi\equiv (P^n_i(t_1,\dots,t_{n}))$, where $P^n_i$ is
a~predicate symbol and $t_1,\dots,t_{n}$ are terms, then
$\mathcal U\vDash \varphi[s]\equiv
(({t_1}_{\mathcal U}[s],\dots,{t_{n}}_{\mathcal U}[s])\in I(P^n_i))$;
\item[--]
if $\varphi\equiv (\neg \theta)$, then
$\mathcal U\vDash \varphi[s]\equiv (\neg \mathcal U\vDash \theta[s])$;
\item[--]
if $\varphi\equiv (\theta_1\Rightarrow \theta_2)$, then
$\mathcal U\vDash \varphi[s]\equiv
(\mathcal U\vDash \theta_1[s]\Rightarrow \mathcal U\vDash \theta_2[s])$;
\item[--]
if $\varphi\equiv (\forall v_i\losp \theta)$, then
$\mathcal U\vDash \varphi[s]\equiv
(\forall x\losp (x\in A\Rightarrow
\mathcal U\vDash \theta[x_1,\dots,x_{i-1},x,x_{i+1},\dots,x_q,\dots]))$.
\end{itemize}

Using the abbreviations cited above, we have also the following:
\begin{itemize}
\item[--]
if $\varphi\equiv (\theta_1\land \theta_2)$,
then $\mathcal U\vDash \varphi[s]\equiv
(\mathcal U\vDash \theta_1[s] \logic\land \mathcal U\vDash \theta_2[s])$;
\item[--]
if $\varphi\equiv (\theta_1\lor \theta_2)$, then
$\mathcal U\vDash \varphi[s]\equiv
(\mathcal U\vDash \theta_1 \logic\lor \mathcal U\vDash \theta_2[s])$;
\item[--]
if $\varphi\equiv (\exists v_i\losp \theta)$, then
$\mathcal U\vDash \varphi[s]\equiv
(\exists x\losp (x\in A \logic\land
\mathcal U\vDash \theta[x_1,\dots,x_{i-1},x,x_{i+1},\dots,x_q,\dots]))$;
\item[--]
if $\varphi\equiv (\theta_1\Leftrightarrow \theta_2)$, then
$\mathcal U\vDash \varphi[s]\equiv
(\mathcal U\vDash \theta_1[s]\Leftrightarrow
\mathcal U\vDash \theta_2[s])$.
\end{itemize}

Models $\mathcal U$~and~$\mathcal U'$ of a~language~$\mathcal L$
are called \emph{isomorphic}
if there exists a~bijective mapping~$f$ of the set (universe)~$A$
onto the set~$A'$
and the following conditions are satisfied:
\begin{enumerate}
\item
for every $n$-place relation~$R$ of the model~$\mathcal U$
and the corresponding
relation~$R'$ of the model~$\mathcal U'$
$R(x_1,\dots,x_n)$ if and only if
$R'(f(x_1),\dots,f(x_n))$ for all $x_1,\dots,x_n$ from~$A$;
\item
for every $m$-place operation~$G$ of the model~$\mathcal U$
and the corresponding operation~$G'$ of the model~$\mathcal U'$
$$
f(G(x_1,\dots,x_m))=G'(f(x_1),\dots,f(x_m))
$$
for all $x_1,\dots,x_m$ from~$A$;
\item
for every constant~$x$ of the model~$\mathcal U$ and the corresponding
constant~$x'$ of the model~$\mathcal U'$
$$
f(x)=x'.
$$
\end{enumerate}

Every mapping~$f$ satisfying these conditions is called
an \emph{isomorphism of the model~$\mathcal U$ onto the model~$\mathcal U'$}
or an \emph{isomorphism between the models $\mathcal U$~and~$\mathcal U'$}.
The fact that $f$ is an
isomorphism of the model~$\mathcal U$ onto the model~$\mathcal U'$
will be denoted by $f\colon \mathcal U\cong \mathcal U'$, and the formula
$\mathcal U\cong \mathcal U'$ means that
the models $\mathcal U$~and~$\mathcal U'$ are isomorphic.

A~model~$\mathcal U'$ is called a~\emph{submodel}
of a~model~$\mathcal U$ if $A'\subset A$ and
\begin{enumerate}
\item
every $n$-place relation~$R'$ of the model~$\mathcal U'$ is
the restriction on the set~$A'$ of the corresponding relation~$R$
of the model~$\mathcal U$, i.e., $R'=R\cap (A')^n$;
\item
every $m$-place operation~$G'$ of the model~$\mathcal U'$ is
the restriction on the set~$A'$ of the corresponding operation~$G$
of the model~$\mathcal U$, i.e., $G'=G|(A')^m$;
\item
every constant of the model~$\mathcal U'$ coincides with
the corresponding constant of the model~$\mathcal U$.
\end{enumerate}

We shall use the notation $\mathcal U'\subset\mathcal U$
to express the fact that
$\mathcal U'$ is a~submodel of the model~$\mathcal U$.
If $\mathcal U$ is a~submodel of a~model~$\mathcal V$,
then $\mathcal V$ is called an \emph{extension} of the model~$\mathcal U$.

Now we shall give a~formal definition of satisfiability.
Let $\varphi$ be an arbitrary
formula of a~language~$\mathcal L$,
let all its variables, free and bound,
be contained in the set $v_0,\dots,v_q$, and let
$x_0,\dots,x_q$ be an arbitrary sequence of elements
of the set~$A$. We define the predicate
$$
\text{\emph{$\varphi$ is true on the sequence $x_0,\dots,x_q$
in the model~$\mathcal U$}, or \emph{$x_0,\dots,x_q$ satisfies the
formula~$\varphi$ in~$\mathcal U$}.}
$$

Let $\mathcal U$ be some fixed model of a~language~$\mathcal L$.
The following sentence shows that the assertion
$$
\mathcal U\models \varphi(v_0,\dots,v_p)[x_0,\dots,x_q]
$$
depends only on the values $x_0,\dots,x_p$, where $p< q$.

\begin{proposition}
\begin{enumerate}
\item
Let $t(v_0,\dots,v_p)$ be a~term, and let $x_0,\dots,x_q$ and
$y_0,\dots,y_r$ be two sequences of elements such that
$p\le q$, $p\le r$, and $x_i=y_i$ whenever $v_i$ is a~free
variable of the term~$t$. Then
$$
t[x_0,\dots,x_q]=t[y_0,\dots,y_r].
$$
\item
Let $\varphi$ be a~formula, let all its variables,
free and bound, belong to the set $v_0,\dots,v_p$, and let
$x_0,\dots,x_q$ and $y_0,\dots,y_r$
be two sequences of elements such that
$p\le q$, $p\le r$, and $x_i=y_i$ whenever $v_i$ is a~free variable
in the formula~$\varphi$. Then
$$
\mathcal U\models \varphi [x_0,\dots,x_q]\quad
\text{if and only if}\quad
\mathcal U\models \varphi [y_0,\dots,y_r].
$$
\end{enumerate}
\end{proposition}

This proposition allows us to give the following definition.
Let $\varphi(v_0,\dots,v_p)$ be a~formula, and let all its variables,
free and bound, be contained in the set $v_0,\dots,v_q$, where $p\le q$.
Let $x_0,\dots, x_p$ be a~sequence of elements of the set~$A$.
We shall say that \emph{$\varphi$~is true in~$\mathcal U$ on
$x_0,\dots,x_p$},
$$
\mathcal U\models \varphi[x_0,\dots,x_p]
$$
if $\varphi$ is true in~$\mathcal U$ on $x_0,\dots,x_p,\dots,x_q$
with some (or, equivalently, any) sequence $x_{p+1},\dots,x_q$.

Let $\varphi$ be a~sentence, and let all its bound variables
be contained in the set $v_0,\dots,v_q$. We shall say that
\emph{$\varphi$~is true in the model~$\mathcal U$}
(notation: $\mathcal U\models \varphi$) if $\varphi$~is true
in~$\mathcal U$ on some (equivalently, any) sequence
$x_0,\dots,x_q$.

Now we say that
$$
\text{a~sentence~$\sigma$ \emph{is true} in~$\mathcal U$,}
$$
if
\begin{gather*}
\text{$\mathcal U\models \sigma [x_0,\dots,x_q]$
\emph{for some} (\emph{or, equivalently, for any})
\emph{sequence $x_0,\dots,x_q$ of elements from~$A$.}}
\end{gather*}

We use special notation
$\mathcal U\models \sigma$ to express this fact.

In the case where $\sigma$ is not true in~$\mathcal U$,
we say that $\sigma$ is \emph{false} in~$\mathcal U$,
or that $\sigma$ \emph{does not hold} in~$\mathcal U$,
or that $\mathcal U$~\emph{is a~model of the sentence} $\neg \sigma$.
If we have a~set~$\Sigma$
of sentences, we say that $\mathcal U$ is a~model of this set
if $\mathcal U$ is a~model of every sentence $\sigma\in \Sigma$.
It is useful to denote this concept by $\mathcal U\models \Sigma$.
A~sentence~$\sigma$ which holds in every model of a~language~$\mathcal L$
is called true. A~sentence (or a~set of sentences)
is called \emph{satisfiable} if it has at least one model.
A~sentence~$\sigma$ is called \emph{refutable}
if $\neg \sigma$ is satisfiable.

A~sentence~$\varphi$ is called a~\emph{corollary} from a~sentence~$\sigma$
(notation: $\sigma\models \varphi$) if every model of the
sentence~$\sigma$ is also a~model of~$\varphi$.
A~sentence~$\varphi$ is called a~corollary of a~set of sentences~$\Sigma$
(notation: $\Sigma \models \varphi$) if every model
of~$\Sigma$ is also a~model of~$\varphi$. Therefore
$$
\Sigma \cup \{ \sigma\}\models \varphi\quad
\text{if and only if}\quad
\Sigma\models \sigma\Rightarrow \varphi.
$$

Models $\mathcal U$~and~$\mathcal V$ of a~language~$\mathcal L$ are called
\emph{elementarily equivalent} if every sentence holds
in~$\mathcal U$ if and only if it holds in~$\mathcal V$.
We express this relation between models
by the notation~$\equiv$. It is clear that the relation~$\equiv$
is an equivalence relation.

Any two isomorphic models of the same language are elementarily equivalent.
If two models of the same language are elementarily equivalent
and one of them is
finite, then they are also isomorphic. If models are infinite
and elementarily equivalent, they are not necessarily isomorphic.
For example, the field~$\mathbb C$ of complex numbers
and the field~$\bar{\mathbb Q}$ of algebraic numbers
are elementarily equivalent,
but not isomorphic, because they have different cardinalities.

Besides first order languages described above, we shall need
to consider second-order languages, in which we can
also quantify predicates, i.e., use predicate symbols as variables.
Such languages will be described in the following sections. We shall
say that two models of the same language
(for example, a~second-order language)~$\mathcal L$ are equivalent
in this language
if for every sentence of the language~$\mathcal L$
it holds in the first model if and only if it holds in the second one.

\subsection{Ultrafilters, Ultraproducts, and Ultrapowers}\label{ss1.4}
The construction of ultraproduct became a~strong instrument in the model
theory. We shall describe it in this section (see~\cite{4}).
{\sloppy

}

Let $I$ be any nonempty set.
By $\mathcal P (I)$ we denote
the set of all subsets of the set~$I$.
A~\emph{filter}~$D$ over the set~$I$ is a~set
$D\subset \mathcal P (I)$ which satisfies the following conditions:
\begin{enumerate}
\item
$I\in D$,
\item
if $X,Y\in D$, then $X\cap Y\in D$,
\item
if $X\in D$ and $X\subset Z\subset I$, then $Z\in D$.
\end{enumerate}

Since $I\in D$, every filter~$D$ is nonempty.
Now we give some examples:
the \emph{trivial filter} $D=\{ I\}$;
the \emph{improper filter} $D=\mathcal P(I)$;
the filter $D=\{ X\subset I\colon Y\subset X\}$ for any
set $Y\subset I$ (this filter is called the \emph{principal filter},
generated by the set~$Y$).

A~filter~$D$ over a~set~$I$
is called an \emph{ultrafilter} over~$I$ if for any
$X\in \mathcal P(I)$
$$
X\in D\quad
\text{if and only if}\quad
(I\setminus X)\notin D.
$$

Let $I$ be any nonempty set, $D$~be a~proper filter over~$I$,
and let $A_i$ be a~nonempty set for each $i\in I$. Consider
$$
C=\prod_{i\in I}A_i,
$$
the Cartesian product of these sets. In other words,
$C$~is the set of all mappings~$f$
which are defined on~$I$ and are such that
$f(i)\in A_i$ for each $i\in I$. The mappings $f,g\in C$
are said to be \emph{$D$-equivalent} (notation: $f=_D g$) if
$$
\{ i\in I\colon f(i)=g(i)\}\in D.
$$

\begin{proposition}
The relation $=_D$ is an equivalence relation on the set~$C$.
\end{proposition}

Now let $f_D$ be the equivalence class which contains the mapping~$f$:
$$
f_D =\{ d\in C\colon f=_D g\}.
$$
We now define the \emph{filter product over sets~$A_i$
by the filter~$D$} as the set of all equivalence classes
of the relation~$=_D$. We denote it by $\smash[b]{\prod\limits_D A_i}$.
Therefore,
$$
\prod_D A_i=\biggl\{ f_D\colon f\in \prod_{i\in I}A_i\biggr\}.
$$
The set~$I$ is called the \emph{set of indices} of $\prod\limits_D A_i$.
If $D$ is an ultrafilter over~$I$,
the filter product $\prod\limits_D A_i$ is called
the \emph{ultraproduct}.
If all $A_i$ coincide (i.e., $A_i=A$),
the filter product is denoted by $\prod\limits_D A$
and called the
\emph{filter power of the set~$A$ by the filter~$D$}. In particular, if
$D$ is an ultrafilter, then $\smash[b]{\prod\limits_D A}$
is called the \emph{ultrapower of the set~$A$ by the filter~$D$}.

Now we give the definition of the filter product of models.
Suppose that $I$ is any nonempty set, $D$~is a~proper filter over~$I$,
and $\mathcal U_i$ is a~model of the language~$\mathcal L$
for every $i\in I$.
We suppose that the predicate symbols~$P$ are interpreted in the
model~$\mathcal U_i$ as~$P_i$, the function symbols~$F$ as~$F_i$,
and the constant symbols~$c$ as~$c_i$.

By definition, \emph{the filter product}
$\prod\limits_D \mathcal U_i$ is
the model of the language~$\mathcal L$ which is defined by the following:
\begin{enumerate}
\renewcommand{\theenumi}{\roman{enumi}}
\item
the universe of the model is the set $\prod\limits_DA_i$;
\item
let $P$ be some $n$-place predicate symbol
of the language~$\mathcal L$. This symbol~$P$ is interpreted in
the model $\prod\limits_D \mathcal U_i$
as the relation~$\bar P$, satisfying the condition
$$
\bar P(f_D^1,\dots,f_D^n)\quad
\text{if and only if}\quad
\{ i\in I\colon P_i(f^1(i),\dots,f^n(i))\}\in D;
$$
\item
let $F$ be some $n$-place function symbol
of the language~$\mathcal L$. The symbol~$F$ is interpreted in
$\smash[b]{\prod\limits_D \mathcal U_i}$
by means of the following mapping~$\bar F$:
$$
\bar F(f_D^1,\dots,f_D^n)=(F_i(f^1(i),\dots,f^n(i))\colon i\in I)_D;
$$
\item
let $c$ be a~constant symbol of the language~$\mathcal L$.
This symbol is interpreted as the element
$$
\bar c=(c_i\colon i\in I)_D
$$
of the set $\prod\limits_D A_i$.
\end{enumerate}

\begin{proposition}\label{ult-t-1}
Let $\prod\limits_D \mathcal U$ be an ultrapower
of a~model~$\mathcal U$.
Then $\prod\limits_D \mathcal U\equiv \mathcal U$.
\end{proposition}

The following important theorem was proved by Keisler and Shelah
(the proof can be found in~\cite{4}).

\begin{theorem}[the isomorphism theorem]\label{ult-t-2}
Let $\mathcal U$~and~$\mathcal V$ be
models of the language~$\mathcal L$. Then
$\mathcal U$~and~$\mathcal V$ are elementarily equivalent
if and only if they have isomorphic ultrapowers.
\end{theorem}

\subsection{Basic Notions from the Category Theory.
Category of Modules over Rings}\label{ss1.5}
We took the basic definitions and notions of this section
from~\cite{5}.


We shall consider an algebraic system~$C$, consisting of two classes
$\Obj$ and $\Mor$, and three operations: \emph{collection},
\emph{composition} (denoted by~$\circ$)
and \emph{identification}, satisfying the following conditions.
\begin{enumerate}
\item
Collection maps every element of the class $\Mor$ to an ordered
pair of elements of the class $\Obj$ (if $f$ is an element of the class
$\Mor$ and $A,B\in \Obj$ are the corresponding elements, then we write
$f\in \Mor(A,B)$).
\item
Composition maps some pairs of elements from $\Mor$ to elements
from $\Mor$ (if $f$,~$g$ are elements from $\Mor$
and $h$~is the corresponding
element from $\Mor$, then we write $h=f\circ g$).
\item
Identification maps every element~$A$ from the class $\Obj$
to some element $f\in \Mor(A,A)$ (we write $f=1_A$).
\item
For every $A\in \Obj $ we have $1_A \in \Mor(A,A)$.
\item
For every $A,B,C\in \Obj$, $f\in \Mor(A,B)$, $g\in \Mor(B,C)$
there exists $h\in \Mor(A,C)$ such that $h=g\circ f$.
\item
For every $A,B,C,D\in \Obj$, $w\in \Mor(A,B)$, $v\in \Mor(B,C)$,
$u\in \Mor(C,D)$ we have $(u\circ v)\circ w=u\circ (v\circ w)$.
\item
For every $A,B\in \Obj$, $u\in \Mor(B,A)$, $v\in \Mor(A,B)$
we have $1_A\circ u=u$ and $v\circ 1_A=v$.
\end{enumerate}


Elements $u\in \Mor(A,B)$ are called \emph{morphisms} from the
object~$A$ into the object~$B$.
The formula $f\in \Mor(A,B)$ will be denoted also by $f\colon A\to B$.

The category $\textup{mod-}R$ of left modules over a~fixed ring~$R$
consists of all left modules over the ring~$R$
and all homomorphisms between them.

If $C$~and~$D$ are categories then a~\emph{covariant functor}
$T\colon C\to D$ is a~pair of mappings
$$
T\begin{cases}
\Obj C\to \Obj D,\\
X\mapsto T X,
\end{cases}\quad
T\begin{cases}
\Mor C\to \Mor D,\\
f\mapsto T f,
\end{cases}
$$
which preserve composition of morphisms and identity morphisms:
\begin{alignat*}{2}
T (f\circ g)&= T f\circ T g &\quad &\forall f,g\in \Mor C,\\
T 1_A     &= 1_{T A} && \forall A\in \Obj C.
\end{alignat*}

A~functor $ T\colon C\to D$ is called \emph{univalent}
if for all objects $X$,~$Y$ of the category~$ C$ the induced
mapping
$$
\begin{cases}
\Mor_C(X,Y)\to \Mor_D ( T X, T Y),\\
f\mapsto T f
\end{cases}
$$
is injective.

The \emph{category of sets} SETS is the category~$C$
in which $\Obj C$ is the class of all sets and $\Mor C$
is the class of all mappings of sets.

A~morphism $f\in \Mor (A,B)$ of a~category~$ C$ is called
an \emph{equivalence}
if there exists a~morphism $g\in \Mor(B,A)$ such that $g\circ f=1_A$ and
$f\circ g=1_B$. A~morphism~$g$ with this property is denoted by~$f^{-1}$.
An object~$A$ is \emph{equivalent} to an object~$B$ (notation:
$A\sim B$) if
there exists an equivalence $f\in \Mor (A,B)$.
It is clear that all these notions can be
expressed in the first order language:
\begin{align*}
\text{$f\in \Mor (A,B)$ is an equivalence}
& \Leftrightarrow \exists g\in \Mor (B,A)\losp
(f\circ g=1_B \logic\land g\circ f=1_A);\\
A\sim B & \Leftrightarrow \exists f\in \Mor (A,B)\losp
(\text{$F$ is an equivalence}).
\end{align*}

In the category $\textup{mod-}R$ an equivalence $f\in \Mor(A,B)$ is called
an \emph{isomorphism} of the modules $A$~and~$B$, and equivalent modules
are called \emph{isomorphic} ($A\cong B$).
An equivalence $f\in \Mor (A,A)$ is called an \emph{automorphism}
of the module~$A$.

Let $S\colon C\to D$ and $T\colon C\to D$ be two covariant functors.
A~\emph{natural transformation} $S\to T$ is a~function~$h$
which maps every object $A\in C$ to a~morphism
$h(A)\colon S(A)\to T(A)$ such that
for every morphism $f\colon A\to A'$ of the category~$C$ we have
$$
T(f)h(A)=h(A')S(f).
$$
A~natural transformation $h\colon S\to T$ between functors $S$~and~$T$
is called
a~\emph{natural equivalence} of $S$~and~$T$ if $h(A)$
is an equivalence for all $A\in \Obj(C)$. In this case, we use the notation
$S\approx T$.

An \emph{equivalence} $C\to D$ between two categories consists of an
ordered pair $(T,S)$ of covariant functors $T\colon C\to D$ and
$S\colon D\to C$ and a~pair of natural equivalences
$$
ST\approx 1_C\ \ \text{and}\ \ TS\approx 1_D
$$
of functors. In this case, we say that $C$~and~$D$ are
\emph{equivalent categories} (notation: $C\approx D$).

An object $T\in \Obj$ of a~category~$ C$ is called a~\emph{left zero}
(an \emph{initial} object) of the category~$C$ if for every object
$X \in \Obj$ there exists a~unique morphism $f\in \Mor (T,X)$.
In the first order language this property can be expressed as
$$
\text{$T$ is a~left zero} \Leftrightarrow
\forall X\in \Obj\losp
\exists f\in \Mor(T,X)\losp \forall g\in \Mor(T,X)\losp (g=f).
$$

An object~$F$ is called a~\emph{right zero}
of a~category~$C$ if for every object $X\in \Obj$
there exists a~unique morphism $f\in \Mor(X,F)$.
An object of a~category~$C$
is called a~\emph{zero} object if it is a~left and right zero simultaneously.
This object is definable in the first order language.
In the category $\textup{mod-}R$ a~zero object is the zero module.

{\sloppy
We say that a~morphism $f\in \Mor(X,Y)$ can be \emph{let trough}
an object~$Q$ if
$\exists g\in \Mor (X,Q)\losp\allowbreak
{\exists h\in \Mor (Q,Y)}\losp\allowbreak (f=h\circ g)$.
A~morphism is called a~\emph{zero} morphism if it can be let trough a~zero
object:
$$
\text{$f\in \Mor(A,B)$ is a~zero morphism}
\Leftrightarrow
\exists g\in \Mor(A,0)\losp \exists h \in \Mor (0,B)\losp (f=h\circ g).
$$

}

In the category $\textup{mod-}R$ zero morphisms between modules
$A$~and~$B$
are morphisms having the form $f(a)=0\in B$ for all $a\in A$.

A~morphism $f\in \Mor (A,B)$ is called a~\emph{retraction} if
$\exists g\in \Mor (B,A)\losp (f\circ g=1_B)$.
A~morphism $f\in \Mor(A,B)$ is called a~\emph{coretraction}
if $\exists g\in \Mor (B,A)\losp (g\circ f=1_A)$. If the category
$\textup{mod-}R$ every retraction
$f\in \Mor (A,B)$ is an epimorphic homomorphism of the module~$A$
onto the module~$B$,
i.e., a~homomorphism $f\colon A\to B$ such that
$\forall b\in B\losp \exists a\in A\losp (f(a)=b)$.
If $f$ is a~retraction $f\colon A\to B$ in the category
$\textup{mod-}R$, then let us consider the set $A'\equiv g[B]$.
It is clear that
$A'$ is a~submodule in~$A$. It is clear that $f|_{A'}\circ g=1_B$.
We show that $g\circ f|_{A'} =1_{A'}$. Let $a\in A'$.
Then $\exists b\in B\losp (g(b)=a)$. In this case, $g(f(a))=g(f(g(b)))=g(b)=a$.
Therefore $A\cong B$. Further, consider $A''\equiv \Ker f$, i.e.,
$a\in A''\Leftrightarrow f(a)=0\in B$. It is clear that
$A=A'\oplus A''$. Thus a~retraction in the category $\textup{mod-}R$
is an isomorphism of some direct
summand of the module~$A$ onto the module~$B$. Similarly, a~coretraction
is an isomorphic embedding of the module~$A$ in the module~$B$ such that
the image of the module~$A$ is a~direct summand in~$B$.

An object~$A$ of a~category~$C$ is called a~\emph{generator} in~$C$ if
$\forall X,Y \in \Obj\losp \forall f,f'\in \Mor(X,Y)\losp\allowbreak
{\exists g \in \Mor(A,X)}\losp (f\circ g \!\ne\! f'\circ g)$.
An object~$A$ is called a~\emph{cogenerator} in~$C$
if $\forall X,Y\in \Obj\losp\allowbreak {\forall f,f'\in \Mor(X,Y)}\losp
\allowbreak
\exists g\in \Mor(Y,A)\losp (g\circ f\ne g\circ f')$.

A~morphism $f\in \Mor(A,B)$ is called a~\emph{monomorphism} if
$\forall C\in \Obj\losp \forall g_1,g_2\in \Mor(C,A)\losp
(f\circ g_1=f\circ g_2\Rightarrow g_1=g_2)$.
A~morphism $f\in \Mor(A,B)$ is called an \emph{epimorphism}
if $\forall C\in \Obj\losp \forall g_1,g_2\in \Mor(B,C)\losp\allowbreak
(g_1\circ f=g_2\circ f \Rightarrow g_1=g_2)$.
A~morphism $f\in \Mor$ is called a~\emph{proper monomorphism} if it is
a~monomorphism and is not an equivalence. We shall say that $f\le g$
for some $f,g\in \Mor$ if $f$~and~$g$ are monomorphisms and
$\exists h\in \Mor\losp (f=g\circ h)$.

An object $A\in \Obj$ is called \emph{projective} if
$$
\forall X,Y\in \Obj\losp \forall f\in \Mor (X,Y)\losp
(\text{$f$ is an epimorphism} \Rightarrow
\forall \tilde g\in \Mor(A,Y)\losp \exists g\in \Mor (A,X)\losp
(\tilde g=g\circ f)).
$$

An object $A\in \Obj$ is called \emph{injective} if
$$
\forall X,Y\in \Obj\losp \forall f\in \Mor(X,Y)\losp
(\text{$f$ is a~monomorphism} \Rightarrow
\forall \tilde g\in \Mor (Y,A)\losp (\tilde g=g\circ f)).
$$

All these properties are elementary, i.e., they can be expressed
in the first order language.

Let $I$ be some subset of the universe and
$\{ A_i\}_{i\in I}$ be a~set of left $R$-modules.
Consider the set~$F$ of mappings from the set~$I$ such that
$\forall i\in I\losp f(i)\in A_i$. On the set~$F$
we can introduce the structure of $R$-module
in the following way: if $f,g\in F$,
then $(f+g)(i)\prisv f(i)+g(i)\in A_i$;
if $f\in F$, $r\in R$, then $(rf)(i)\prisv rf(i)\in A_i$. This module~$F$
is called the \emph{product of the set of modules} $\{ A_i\}_{i\in I}$
and is denoted by $\prod\limits_{i\in I} A_i$.
If $A_i=A$ for all $i\in I$, then the product
$\prod\limits_{i\in I}A_i$ is denoted by~$A^I$. For every $k\in I$ the set
of mappings satisfying the condition $f(i)=0$ for $i\ne k$ is
a~module which is isomorphic to the module~$A_k$. Such a~module
will be considered as a~natural embedding of the module~$A_k$ in the module
$\prod\limits_{i\in I} A_i$.

Further, consider the set $S$ of mappings~$f$ from the set~$I$
such that $\forall i\in I\losp f(i)\in A_i$ and
$f(i)\ne 0$ only for a~finite set of elements from~$I$.
On the set~$S$ we can similarly introduce the structure of $R$-module.
The module~$S$ is called the \emph{direct sum
of a~set of modules} $\{ A_i\}_{i\in I}$
and is denoted by $\smash[b]{\bigcup\limits_{i\in I} A_i}$.
If $A_i=A$ for all $i\in I$, then the direct sum
$\bigcup\limits_{i\in I}A_i$ is denoted by~$A^{(I)}$.

The product of a~finite set $A_1,\dots, A_n$ is denoted by
$$
\prod_{i=1}^n A_i\quad \text{or}\quad A_1\times \dots \times A_n,
$$
and the direct sum is denoted by
$$
\bigcup_{i=1}^n A_i\quad \text{or}\quad A_1\oplus \dots \oplus A_n.
$$

An $I\times J$-\emph{matrix} over a~set $S$ is a~mapping
$f\colon I\times J\to S$. Therefore a~matrix is an element of the set
$S^{I\times J}$. If $S$ contains only two different elements $0$~and~$1$,
then the \emph{Kronecker delta} is a~matrix $\delta\colon I\times I\to S$
such that $\delta_{ii}=1$ and $\delta_{ij}=0$ if $i\ne j$.

\begin{proposition}\label{kateg1}
Let $C$ be a~category with zero. If $A=\prod\limits_{i\in I}A_i$
is a~product in the category~$C$,
then there exists a~family of retractions
$p_i\colon A\to A_i$ and a~family of coretractions $u_i\colon A_i\to A$
which can be uniquely defined by the condition
$$
p_i u_j=\delta_{ij} 1_{A_i}
$$
for all $i,j\in I$.

Dually, if $\bigcup\limits_{i\in I} A_i$
is a~direct sum, then there exist coretractions
$u_i\colon A_i\to A$ and there exist uniquely defined by~$u_i$
retractions $p_i\colon A \to A_i$ such that
$$
p_iu_j=\delta_{ij} 1_{A_i}
$$
for all $i,j\in I$.
\end{proposition}

Product of two objects is definable in the first order language.
The same is fulfilled for a~product and direct sum of any given
finite number of objects.

Let $C$ be some concrete category. For an arbitrary set~$S$ consider the
category $(S, C)$, with objects which are mappings $f\colon S\to A$,
where $A$ is an object of the category~$ C$.
Morphisms of the category $(S, C)$
are defined as morphisms $A\to B$ of the category~$ C$ such that
for the given objects $S\to A$ and $S\to B$ of the category $(S, C)$
the diagram
$$
\begin{array}{ccc}
S& \to& A\\
\searrow&&\swarrow\\
&B&
\end{array}
$$
is commutative.

A~left zero $f_S\colon S\to F(S)$ of a~category $(S, C)$ is called
a~\emph{free object of the category~$ C$ over the set~$S$}.

In other words, for every mapping $f\colon S\to A$ there exists
a~unique morphism $h\colon F(S)\to A$ such that
$$
f=h\circ f_S.
$$

For objects of the category $\textup{mod-}R$ the following
notions are definable in the first order language:
\begin{itemize}
\item[]
modules $X$~and~$Y$ are isomorphic;
\item[]
a module $X$ is embeddable in a~module $Y$;
\item[]
there exists a~surjection from a~module~$X$ onto a~module~$Y$;
\item[]
a~module~$X$ is isomorphic to a~direct summand of a~module $Y$;
\item[]
a~module~$X$ is isomorphic to the direct sum of modules
$Y$~and~$Z$;
\item[]
a~module~$X$ is projective;
\item[]
a~module~$X$ is injective;
\item[]
a~module~$X$ is a~generator;
\item[]
a~module~$X$ is a~cogenerator.
\end{itemize}

In the general case the following
properties of modules of the category
$\textup{mod-}R$ are not definable:
\begin{itemize}
\item[]
a module $X$ is free;
\item[]
a module $X$ is equal to $A^I$ for some set $A$;
\item[]
a module $X$ is equal to $A^{(I)}$ for some set $A$.
\end{itemize}

\section[An Analogue of the Morita Theorem
for Elementary Equivalence of Categories of Modules]{An Analogue
of the Morita Theorem\\
for Elementary Equivalence of Categories of Modules}
In 2003 we studied elementary properties of categories of modules
over rings, endomorphism rings of modules, and automorphism groups
of modules over rings. Our interest to these questions was attracted
by the paper~\cite{16} of V.~Tolstykh.

\subsection[Some Facts
about the Category $\textup{mod-}R$]{Some Facts
about the Category $\boldsymbol{\textup{mod-}R}$}\label{ss2.1}
A~\emph{factormodule} of a~module~$M$ by its submodule~$N$ is the module
consisting of all equivalence classes $a\sim b\Leftrightarrow a-b\in N$
and such that $(a+N)r=ar+N$. The property of a~module~$L$ to be isomorphic
to a~factormodule of the module~$M$ is a~first order property:
$\exists f\in \Mor(M,L)\losp (\text{$f$ is an epimorphism})$.

Let $\mathbf C$ be some concrete category. If $B$~and~$A$ are its objects
and $B\subseteq A$, then $B$ is a~\emph{subobject} in~$A$. If $A$ is
a~subset and $N$ is a~subobject in~$A$, then $S$ \emph{generates}~$N$
if $N$ is an intersection of all subobjects of the object~$A$
containing~$S$. In this case, we use the notation $N=(S)$.
A~subobject~$M$ of an object~$A$ is called \emph{finitely generated},
\emph{countably generated}, or \emph{generated by $a$~elements}
if $M=(T)$, where $|T|< \omega_0$, $|T|\le \omega_0$, or $|T|\le a$,
respectively. In the general case these properties are not elementary.

A~family $\{ x_i\}_{i\in I}$ which generates a~submodule~$N$
of a~module~$M$ is called
a~\emph{system of generators} of the submodule~$N$.
If every element of the module can be uniquely represented
as a~linear combination of generators, then $\{ x_i\}_{i\in I}$
is called a~\emph{basis} of the module~$N$,
and the cardinality of the set~$I$ is called a~\emph{basic
number} of this module. A~family $\{ x_i\}_{i\in I}$ is called
\emph{linearly independent over~$R$}.

The module $R^{(I)}$ is a~free module over the set~$I$.

\begin{proposition}
\begin{enumerate}
\item
If $R$ is a~ring and $X$~is an object of the category $\textup{mod-}R$,
then there exist a~set of indices~$I$ and some epimorphism
$f\colon R^{(I)}\to X$, i.e., every $R$-module is isomorphic
to a~factormodule of a~free $R$-module.
\item
If $\{ u_i\colon R\to R^{(I)}\}_{i\in I}$ are injections
into the direct sum, then
$\{ u_i(1)\}_{i\in I}$ is a~basis of the free module $R^{(I)}$.
\item
The object~$R$ is a~generator in the category $\textup{mod-}R$.
\end{enumerate}
\end{proposition}

The basic number in a~general situation depends on the choice of the basis
and
therefore can not be an invariant of the module $F=R^{(I)}$. But
it does not depend on the choice of the basis if $F$ is a~free module under
an infinite set~$I$.

\begin{proposition}
An $R$-module~$P$ is projective if and only if
it is isomorphic to a~direct summand of a~free module.
\end{proposition}

\begin{corol}
A~module~$P$ is finitely generated and projective if and only if
$R^{(n)}\cong P\oplus X$ for some integer $n> 0$ and some
module~$X$.
\end{corol}

\begin{proof}
If $R^{(n)}\cong P\oplus X$ for some integer $n> 0$, then
it is clear that $P$ is projective and finitely generated.

Conversely, let $P$ be finitely generated and projective. Since $P$ is
projective, we have $P\oplus Q\cong R^{(I)}$ for some
set~$I$. Let the set~$I$ be infinite. Consider the set
$\{ p_1,\dots,p_k\}$ of elements which generate~$P$
and the basis $\{ e_i\}_{i\in I}$ of the module $R^{(I)}$. Every~$p_j$
is a~linear combination
of a~finite number of elements of the basis, therefore
only a~finite subset of $\{ e_i\}_{i\in I}$
belongs to all linear combinations
for all~$p_i$. Thus
$P\subset R^{(n)}\subset R^{(I)}$, and $R^{(n)}$ is a~direct summand
in~$R^{(I)}$. Consequently, $P\oplus (Q\cap R^{(n)})\cong R^{(n)}$.
\end{proof}

\begin{proposition}
A~module $M\in \textup{mod-}R$ is a~generator if and only if
every $R$-module~$X$ is a~factormodule of the module $M^{(I)}$
for some set~$I$.
\end{proposition}

\begin{proposition}
An object~$G$ of the category $\textup{mod-}R$ is a~generator if and only if
there exist an integer $n> 0$ and an isomorphism
$G^{(n)}\cong R\oplus X$ for some object $X\in \textup{mod-}R$.
\end{proposition}

A~module~$M$ is called \emph{simple} if it has only two submodules
$0$~and~$M$. If $M$ is some module and $S$ is its submodule, then
$M/S$ is simple if and only if $S$ is a~maximal submodule. Every
finitely generated module~$M$ has maximal submodules.
Therefore for every ring~$R$ in the category $\textup{mod-}R$ there exist
some simple modules (they can be isomorphic to each other).
It is clear that the property of a~module to be simple is definable in the
first order language.

\begin{proposition}\label{sim}
For every simple module~$M$ every submodule~$P$ of the module $M^{(I)}$
is isomorphic to $M^{(J)}$ for some set~$J$ with cardinality
not greater than the cardinality of~$I$.
\end{proposition}


A~module $P\in \textup{mod-}R$ is
called a~\emph{progenerator} if it is finitely generated and projective and
$P$ is a~generator in $\textup{mod-}R$.

Two rings $R$~and~$S$ are called \emph{similar} (denoted by $R\sim S$)
if there exist a~progenerator $P\in \textup{mod-}R$ and
a~ring isomorphism $S\cong \Endom_R P$.

The following famous theorem is cited here without proof
(its proof can be found in~\cite[Theorem~4.29]{5}).

\begin{theorem}[Morita theorem]\label{morita}
The following conditions are equivalent\textup:
\begin{enumerate}
\item
$\textup{mod-}R\approx \textup{mod-}S$\textup;
\item
$R\sim S$.
\end{enumerate}
\end{theorem}

In the sequel, we shall also need the following theorem
from~\cite{5} (see~4.35).

\begin{theorem}\label{kommut}
If $A$ is a~commutative ring and a~ring~$B$ is similar to the ring~$A$,
then $A$ is isomorphic to the center of~$B$.
Therefore two commutative rings
are similar if and only if they are isomorphic.
\end{theorem}

\subsection[Progenerators in the Category
$\textup{mod-}R$]{Progenerators in the Category
$\boldsymbol{\textup{mod-}R}$}\label{ss2.2}
Let a~formula $\Simp(M)$ be true in the category $\textup{mod-}R$
for simple modules and only for them.
Consider an object~$X$ satisfying the formula
$$
\tSum^\omega (X,M) \prisv \Simp(M) \logic\land
(X\oplus M\cong X) \logic\land
\forall Y\in \Obj\losp
(Y\oplus M\cong Y \Rightarrow \exists Q\in \Obj \losp (Y\cong X\oplus Q)).
$$

The property $Y\oplus M\cong Y$ means that $Y\cong M^{(\omega)}\oplus Z$
for some object $Z\in \textup{mod-}R$. Therefore $X$ is a~module
which contains $M^{(\omega)}$ as its direct summand and it itself
is a~direct summand in~$M^{(\omega)}$. It follows from Proposition~\ref{sim}
that in this case $X\cong M^{(\omega)}$. Thus
for every simple module~$M$
the formula $\tSum^\omega_M(X) \prisv \tSum^\omega(X,M)$
defines the module~$M^{(\omega)}$.

The formula
$$
\tSum^{\Fin}(X,M) \prisv \tSum_M^{\Fin}(X) \prisv
\Simp(M) \logic\land \exists Y\in \Obj\losp
(\tSum^\omega_M(Y) \logic\land \exists Q\in \Obj\losp
(Y\cong X\oplus Q) \logic\land X\not\cong Y)
$$
holds for all finite direct sums of the simple module~$M$ and only for them.

The formula
$$
\tSum(X,M) \prisv \tSum_M(X) \prisv
\Simp(M) \logic\land
\forall Y\losp
(Y\subset X \logic\land Y\ne 0\Rightarrow \exists P\losp (Y\cong P\oplus M))
$$
defines the class $\tSum_M$ of all direct sums of the module~$M$.
Introduce the relation
$$
(X\le Y) \prisv \exists f\in \Mor(X,Y)\losp (\text{$f$ is a~monomorphism})
$$
on this class.

The class $\tSum_M$ is well-ordered with respect to the order~$\le$
and there exists a~natural bijection from the class $\tSum_M$
onto the class $\Cn$ of all cardinal numbers.

The formula
\begin{multline*}
\Pret(P) \prisv (\text{$P$ is projective}) \logic\land
(\text{$P$ is a~generator})\\
\logic\land
\exists M\in \Obj\losp \exists f\in \Mor(P,M)\losp
(\Simp(M) \logic\land (\text{$f$ is an epimorphism}))
\end{multline*}
holds for all projective modules having maximal
submodules, in particular it must hold for projective finitely generated
generators (\emph{progenerators}).

By $\langle M,f \rangle^P$ (or $\langle M^P, f^P\rangle$)
we shall denote a~pair (a~simple module~$M$, an epimorphism~$f$
from~$P$ onto~$M$) for a~module~$P$ such that $\Pret(P)$.

Consider a~module~$N$ satisfying the formula
$\tSum^{\Fin}_{M^P}(N)$. Such a~module~$N$ has the form $M^{(n)}$
for some natural~$n$. We shall denote this module by
$N_{\mathrm{fd}}(M)$.

Consider now the formula
\begin{align*}
& \Under (P,M,N,X) \prisv \Under_{M^P,N}(X) \prisv
N\cong N_{\mathrm{fd}}(M) \logic\land
\exists g\in \Mor(X,N)\losp\\
& \quad (\text{$g$~is an epimorphism} \logic\land
\forall i_M\in \Mor (M,N)\losp \forall p_M \in \Mor (N,M)\losp\\
& \quad (p_M\circ i_M=1_M\Rightarrow
\exists i\in \Mor(P,X)\losp \exists p\in \Mor(X,P)\losp
(p\circ i=1_P \logic\land g\circ i=i_M\circ f \logic\land
f\circ p=p_M\circ g))\\
& \quad \logic\land
\forall i_M,i_M'\in \Mor (M,N)\losp \forall p_M,p_M'\in \Mor(M,M)\losp
\forall i,i'\in \Mor(P,X)\losp \forall p,p'\in \Mor(X,P)\losp\\
& \quad (p_M\circ i_M=p_M'\circ i_M'=1_M \logic\land
p\circ i=p'\circ i'=1_P \logic\land
g\circ i=i_M\circ f \logic\land f\circ p=p_M\circ g\\
& \quad \logic\land
g\circ i'=i_M'\circ f \logic\land f\circ p'=p_M'\circ g \logic\land
p_M\circ i_M'=p_M'\circ i_M=0 \Rightarrow p\circ i'=p'\circ i=0)).
\end{align*}

This formula means that
\begin{enumerate}
\item
for the module~$X$ there exists an epimorphism
$g\colon X\to N$
such that for every pair $(i_M,p_M)$ consisting of an embedding
of the module~$M$ into the module~$N$ and an inverse projection of
the module~$N$ onto the module~$M$
there exists a~pair $(i,p)$ consisting of an embedding of the module~$P$
into the module~$X$ and an inverse projection of the module~$X$
onto the module~$P$ such that the diagrams
$$
\begin{array}{ccc}
P&\stackrel{i}{\rightarrow}& X\\
\downarrow\lefteqn{f}&&\downarrow\lefteqn{g}\\
M&\stackrel{i_M}{\rightarrow}& N
\end{array}\quad \text{ and }\quad
\begin{array}{ccc}
P& \stackrel{p}{\leftarrow}& X\\
\downarrow\lefteqn{f}&&\downarrow\lefteqn{g}\\
M& \stackrel{p_M}{\leftarrow}& N
\end{array}
$$
are commutative;
\item
if embeddings $i_M$ and $i_M'$ of the module~$M$ into the module~$N$
are such that their images in~$N$ do not intersect, then the images
of the corresponding embeddings $i,i'\colon P\to X$ also do not intersect.
\end{enumerate}

Look at the module~$X$ in this case.

Suppose that $N\cong M^{(n)}\cong M_1\oplus \dots \oplus M_n$,
where $M_i\cong M$ for every $1\le i\le n$.
Let $i_l^M\colon M\to N$ and $p_l^M\colon N\to M$
be such that $\rng i_l^M =M_l$ and $p_l^M\circ i_l^M=1_M$.
To these pairs of embeddings correspond the pairs
$(i_l,p_l)$ such that $i_l\colon P\to X$, $p_l\colon X\to P$,
$p_l\circ i_l=1_P$,
and the images of embeddings $i_l$~and~$i_m$ for different $l$~and~$m$
do not intersect and are independent. Therefore the module
$P^{(n)}$ is a~direct summand in~$X$.
Now we only need to consider a~module~$X'$ satisfying the formula
$$
\Und(P,M,N,X') \prisv \Und_{N,M^P}(X') \prisv
\forall X\losp
(\Under_{M^P,N}(X) \Rightarrow \exists Q\losp (X\cong X'\oplus X')).
$$
We shall get the module~$X'$ which is a~direct summand
in the module~$P^{(n)}$ and has the module $P^{(n)}$ as its direct summand.

Now consider the following formula:
$$
\Finite(P,X)\prisv \Finite_P(X)\prisv
\exists (M^P,f^P)\losp \exists Y\in \Obj\losp
(\tSum^{\Fin}_M(Y) \logic\land \Und_{Y,M^P}(X)).
$$

This formula defines modules~$X$ with the property
$$
\exists n\in\omega\losp \exists Q,Q'\losp
(X\oplus Q\cong P^{(n)} \logic\land X\cong P^{(n)}\oplus Q'),
$$
i.e., all modules of the form $P^{(n)}$ and some other \emph{finitely
generated} modules.

Every projective finitely generated module
is a~direct summand of the module $R^{(n)}$
for some~$n\in \omega$, and, respectively,
if $P$ is finitely generated and projective, then for every
generator~$S$
$$
P\oplus Q\cong S^{(m)}
$$
for some $m\in \omega$ and some module~$Q$. If a~module~$P$
is not finitely generated,
then there exists a~progenerator~$S$ such that $P$ can not
be embedded in $S^{(n)}$ for any $n\in \omega$.

Therefore the formula
$$
\Proobr(P)\prisv \Pret(P) \logic\land
\forall S \in \Obj\losp
(\Pret(S)\Rightarrow \exists X \in \Obj\losp
(\Finite_S(X) \logic\land \exists Q\in \Obj\losp (P\oplus Q\cong X)))
$$
defines all progenerators in the category $\textup{mod-}R$.

Thus having the category $\textup{mod-}R$ we automatically have
(with the help of the formula $\Proobr()$) the class of all
progenerators in this category.

Note also that having some fixed progenerator~$P$,
we have the class of all modules which are direct summands in $P^{(I)}$
and simultaneously have $P^{(I)}$ as a~direct summand.
It is clear that each such module has the form $P^{(I)}\oplus X$, where $X$
is some projective module which can be embedded in $P^{(I)}$.
Every such module can be represented as $R^{(I)}\oplus Y$, where $Y$ is
a~projective module of rank~${\le}\,|I|$. We shall call
such modules \emph{almost free modules of rank~$|I|$ over the ring~$R$}.

\subsection[The Ring $\protect\Endom_RP$]{The Ring
$\boldsymbol{\Endom_RP}$}\label{ss2.3}
Consider now some progenerator~$P$ and the set $\Mor(P,P)$.
The operation of multiplication on this set can be introduced as
$$
(f=g\times h) \prisv (f=g\circ h).
$$
Introduce now the operation of addition. For this purpose we consider
the module $P\oplus P$
with two embeddings $i_1,i_2\in \Mor(P,P\oplus P)$ and two projections
$p_1,p_2\in \Mor(P\oplus P,P)$ with the conditions
$p_1\circ i_1=p_2\circ i_2=1_P$, $p_1\circ i_2=p_2\circ i_1=0$.

For a~given $f\in \Mor(P,P)$ consider the morphism
$\Gr_f\in \Mor(P\oplus P,P\oplus P)$ which is defined by the conditions
\begin{align*}
p_1\circ \Gr_f \circ i_1&=1_P,\\
p_2\circ \Gr_f \circ i_2&=1_P,\\
p_2\circ \Gr_f\circ i_1&=0,\\*
p_1\circ \Gr_f\circ i_2&=f.
\end{align*}

It is clear that the mapping
$$
\Gr\colon \Mor(P,P)\to \Mor(P\oplus P,P\oplus P),\quad
f\mapsto \Gr_f,
$$
is injective and that for every morphism $g\in \Mor(P\oplus P,P\oplus P)$
satisfying the conditions
$p_1\circ g\circ i_1=p_2\circ g\circ i_2=1_P$
and $p_2\circ g\circ i_1=0$ there
exists a~morphism $f\in \Mor (P,P)$ such that $\Gr_f=g$.

Define
$$
(f=g+h) \prisv (\Gr_f=\Gr_g\circ \Gr_h).
$$
Thus we have introduced on the set $\Mor(P,P)$ the structure of
ring which is isomorphic to the ring $\Endom_R(P)$.

Indeed, let us show that
for any three endomorphisms $f,g,h\in \Endom_R P=\Mor(P,P)$ the relation
$f=g+h$ is true if and only if $\Gr_f=\Gr_g\circ \Gr_h$.
Consider the morphisms $\Gr_g$ and $\Gr_h$
and the morphism $G=\Gr_g\circ \Gr_h$.
The mappings $k_1\equiv i_1\circ p_1$ and $k_2\equiv i_2\circ p_2$
from $\Mor(P\oplus P,P\oplus P)$ are such that
$\forall x\in P\oplus P\losp (x=k_1(x)+k_2(x))$, i.e.,
$k_1+k_2=1_{P\oplus P}$. Thus,
\begin{align*}
& p_1\circ G\circ i_1=p_1\circ \Gr_g\circ \Gr_h\circ i_1=
p_1\circ \Gr_g\circ 1_{P\oplus P} \circ \Gr_h \circ i_1\\
& \quad = p_1\circ \Gr_g\circ
(i_1\circ p_1+i_2\circ p_2)\circ \Gr_h\circ i_1=
p_1\circ \Gr_g\circ (i_1\circ p_1\circ \Gr_h\circ i_1+i_2\circ p_2\circ
\Gr_h\circ i_1)\\
& \quad = p_1\circ \Gr_g\circ i_1\circ 1_P+0=1_P,
\end{align*}
and, similarly,
$p_2\circ G\circ i_2=1_P$, $p_2\circ G\circ i_1=0$ and, finally,
\begin{multline*}
p_1\circ G\circ i_2=p_1\circ \Gr_g\circ \Gr_h\circ i_2=
p_1\circ \Gr_g\circ (i_1\circ p_1+i_2\circ p_2)\circ \Gr_h\circ i_2\\
=(p_1\circ \Gr_g \circ i_1)\circ (p_1\circ \Gr_h\circ i_2)+
(p_1\circ \Gr_g\circ i_2)\circ (p_2\circ \Gr_h\circ i_2)=
g\circ 1_P+ 1_P\circ h=g+h.
\end{multline*}
Thus we get the required equivalence.

\subsection{The Case of Finite Rings}\label{ss2.4}
\begin{lemma}\label{mod-l-1}
The endomorphism ring $\Endom_R P$ of any progenerator~$P$
of the category $\textup{mod-}R$ with a~finite ring~$R$ is finite.
\end{lemma}

\begin{proof}
The module~$P$ is a~submodule of the module $R^{(n)}$ for some~$n$.
Since the ring~$R$ is finite, also the module $R^{(n)}$ is finite and
therefore the module~$P$ is also finite. It is clear that
the endomorphism ring of a~finite module is finite.
\end{proof}

\begin{lemma}\label{mod-l-2}
For every finite ring~$R$ there exists a~sentence~$\varphi_R$
of the first order language of the ring theory
which is true in a~ring~$X$ if and only if $X\cong R$.
\end{lemma}

\begin{proof}
Consider a~finite ring~$R$. Suppose that it contains exactly
$m$ different elements $a_1,\dots,a_m$, and $a_i+a_j=a_{s(i,j)}$,
$a_i\cdot a_j=a_{p(i,j)}$. Then the required sentence~$\varphi_R$
has the form
$$
\exists x_1 \dots \exists x_m\losp
\biggl(\,\bigwedge_{i,j\in m, i\ne j} x_i\ne x_j \biggr) \logic\land
\biggl(\forall x \bigvee_{i\in m} x=x_i\biggr) \logic\land
\biggl(\,\bigwedge_{i,j\in m} a_i+a_j=a_{s(i,j)}\biggr) \logic\land
\biggl(\,\bigwedge_{i,j\in m} a_i\cdot a_j=a_{p(i,j)}\biggr).
$$
\end{proof}

\begin{theorem}\label{mod-t-1}
If categories $\textup{mod-}R$ and $\textup{mod-}S$
are elementarily equivalent
and the ring~$R$ is finite, then $R\cong \Endom_S P$ for some
progenerator module~$P$ of the category $\textup{mod-}S$.
\end{theorem}

\begin{proof}
In the category $\textup{mod-}R$ the sentence
$$
\xi \prisv \exists P\in \Obj\losp
(\Proobr(P) \logic\land \varphi_{\Mor(P,P)})
$$
is true.
Therefore, the sentence~$\xi$ is true in the category $\textup{mod-}S$,
i.e., the endomorphism ring of some progenerator
is isomorphic to the ring~$R$.
\end{proof}

\begin{corol}
The categories $\textup{mod-}R$ and $\textup{mod-}S$,
where $R$ is a~finite ring,
are elementarily equivalent if and only if they are Morita-equivalent.
\end{corol}

\begin{proof}
If categories $\textup{mod-}R$ and $\textup{mod-}S$ are Morita-equivalent,
they are clearly elementarily equivalent.

If categories $\textup{mod-}R$ and $\textup{mod-}S$
are elementarily equivalent
and the ring~$R$ is finite, then, by Theorem~\ref{mod-t-1},
$R\cong \Endom_S P$ for some progenerator~$P$ of the category
$\textup{mod-}S$, i.e., the rings $R$~and~$S$
are similar. By the Morita Theorem (Sec.~2.1, Theorem~\ref{morita}),
in this case the categories $\textup{mod-}R$ and
$\textup{mod-}S$ are Morita-equivalent.
\end{proof}

Now we assume that the rings $R$~and~$S$ are infinite.

\subsection{Beautiful Linear Combinations}\label{ss2.5}
We apply the results of S.~Shelah (1976) (see~\cite{6})
on interpretation of the set theory in a~category.

Suppose that we have some fixed ring~$R$,
the category $\textup{mod-}R$, and in the category $\textup{mod-}R$ we have
some simple module~$M$ which corresponds to the fixed progenerator~$P$,
$V=M^{(I)}$, $|I|=\mu$, where $\mu$ is an infinite cardinal number.
Let a~set $A=\{ a_i\mid i\in I\}$ be such that
$\forall i\in I\losp (a_i\in M_i \logic\land a_i\ne 0)$.

For every $f\in \Mor(A,B)$ let $\Rng f$ be the image of~$f$ in~$B$,
$\Cl B$ be the closure of the set $B\subset V$ in~$V$,
i.e., the smallest submodule in~$V$ containing the set~$B$.
Further, let $\tilde b = \Cl \{ b\}$.

As usual, $\vec x$ denotes a~finite sequence
of variables $\vec x=\langle x_1,\dots, x_n\rangle$.
A~linear combination $\alpha_1 x_1+\dots +\alpha_n x_n$,
where $\alpha_i\in R$, will be also denoted by $\tau(x_1,\dots,x_n)$
or $\tau(\vec x)$. We shall call such a~linear combination
\emph{reduced} if all~$\alpha_i$ are nonzero.

A~linear combination $\tau(x_1,\dots,x_n)=\alpha_1 x_1+\dots+\alpha_nx_n$
is called \emph{beautiful} (see~\cite{6}) if
\begin{enumerate}
\renewcommand{\theenumi}{\alph{enumi}}
\item
for every linear combination
$\sigma(x_1,\dots,x_m)=\beta_1 x_1+\dots+\beta_m x_m$ we have
the equality
\begin{multline*}
\tau(\sigma(x_1^1,\dots,x_m^1),\sigma(x_1^2,\dots,x_m^2),\dots,
\sigma(x_1^n,\dots, x_m^n))\\
=\sigma(\tau(x_1^1,x_1^2,\dots,x_1^n),\tau(x_2^1,x_2^2,\dots,x_2^n),\dots,
\tau(x_m^1,x_m^2,\dots,x_m^n));
\end{multline*}
\item
we have the equality
$$
\tau(\tau(x_1^1,\dots,x_n^1),\tau(x_1^2,\dots, x_n^2),\dots
\tau(x_1^n,\dots, x_n^n))
=\tau(x_1^1,\dots,x_n^n);
$$
\item
we have the equality
$$
\tau(x,\dots,x)=x.
$$
\end{enumerate}

It is easy to show that all beautiful linear combinations
have the form
$$
\alpha_1 x_1+\dots+\alpha_n x_n, \quad
\alpha_i\in Z(R),\ \ \alpha_i \alpha_j=\delta_{ij}\alpha_i,\ \
\sum_{i=1}^n \alpha_i =1.
$$

\begin{theorem}\label{shel_t4_1}
There exists a~formula~$\varphi_m$ satisfying the following condition.
Let $\bar f_i$ be an $m$-tuple of elements of $\Mor(V,M)$
for every $i< i_0< \mu^+$.
Then we can find a~vector $\bar g$ such that the formula
$\varphi_m(\bar f, \bar g)$ holds in $\textup{mod-}R$ if and only if
$\bar f= \tau(\bar f_{i_1},\dots, \bar f_{i_n})$ for some
beautiful linear combination~$\tau$ and some
$i_1< \dots< i_n< i_0< \mu^+$.
\end{theorem}

\subsection[A Generating Set of the Module~$V$]{A Generating Set
of the Module~$\boldsymbol{V}$}\label{ss2.6}
Recall that by $V$ we denote a~module $M^{(\mu)}$ for some
infinite cardinal number~$\mu$ and a~fixed simple
module~$M$.

Let $V=\bigcup\limits_{t\in \mu} M_t$, where $M_t\cong M$ for every
$t\in \mu$, and suppose that in the module~$M$ some generating
(i.e., nonzero) element~$a$ is fixed and in every~$M_t$ the corresponding
element $a_t$ is fixed.

We shall use Theorem~\ref{shel_t4_1} for
$m=1$ and $f_i\in \Mor(V,M)$ such that $f_i(a_t)=\delta_{it}a$.
Then there exist $\bar g^*$ and a~formula $\varphi(f,\bar g^*)$
such that the formula $\varphi(f,\bar g^*)$ holds if and only if
$f=\tau(f_{i_1},\dots,f_{i_n})$, where $i_1< \dots< i_n< \mu$
and the linear combination~$\tau$ is beautiful.

We know that in this case
$\tau(x_1,\dots,x_n)=r_1 x_1+\dots+r_n x_n$,
where $r_ir_j=\delta_{ij} r_i$ for all $i,j=1,\dots,n$ and
$r_1+\dots+r_n=1$.

Consider a~mapping $r_k f_{i_k}\colon V\to M$. We know that
$r_k f_{i_k}(a_{i_k})=r_k\cdot a$ and $r_k f_{i_k}(a_t)=0$ for $t\ne i_k$.
In the module~$M$ consider a~set $N\subseteq M$ such that
$n\in N\Leftrightarrow r_k\cdot n=0$.
If $n_1,n_2\in N$, then $r_k (n_1+n_2)=0$, whence $n_1+n_2\in N$.
If $r\in R$, $n\in N$, then $r_k(rn)=r(r_k n)=0$, whence $rn\in N$.
Thus $N$ is an ideal in~$M$, i.e., $N=\{ 0\}$ or $N=M$.
Let $r_ka\ne 0$ and $r_l a\ne 0$ for some different $k$ and~$l$.
Then $r_kb\ne 0$ and $r_l b\ne 0$ for all $b\in M$, i.e.,
$r_l(r_ka)\ne 0$, but this is impossible. Therefore
$r_k a\ne 0$ only for one $k\in \{ 1,\dots,n\}$. From $r_1+\dots+r_n=1$,
i.e., from $(r_1+\dots+r_n)a=a$, it follows that such~$k$ necessarily
exists and $r_k a=a$. Consequently, for some~$k$ we have
$r_k f_{i_k}(a_{i_k})=a$ and $r_kf_{i_k}(a_t)=0$ for $t\ne i_k$,
and for $l\ne k$ we have $r_l f_{i_l}(a_t)=0$ for all $t\in I^*$.
Thus, $f=f_{i_k}$ for some $k\in \{ 1,\dots,n\}$.

Thus we have shown that there exists such $\bar g^*$ that the formula
$\varphi(f,\bar g^*)$ defines in~$V$ some set consisting of
$\mu$ independent
projectors from~$V$ onto~$M$.
We shall obtain the required~$\bar g^*$
if we write a~formula stating
that
the space generated by the images of
those $i\in \Mor(M,V)$ that satisfy
$\exists f\losp (\varphi(f,\bar g^*) \logic\land f\circ i=1_M)$
is isomorphic to~$V$ and, if we exclude any pair $(f,i)$
from this space, the new space will not coincide with the initial one.

Recall that together with the simple module~$M$ we have fixed
a~progenerator~$P$ and an epimorphism $h\colon P\to M$,
and together with the module $V\cong M^{(\mu)}$ we have fixed
a~module~$V'$ which is an almost free module of rank~$\mu$ over~$P$
with an epimorphism $h'\colon V'\to V$
such that for every projection $i\colon V\to M$
there exists a~unique projection $i'\colon V'\to P$
such that $i\circ h'=h\circ i'$.

The set consisting of all projections $g\in \Mor(V,M)$
satisfying the formula $\varphi (\bar g^*, g)$ will be denoted by
$\Gen_{\bar g^*}(V,M)$. The set
consisting of all projections $g\in \Mor(V',P)$ satisfying the formula
$\exists f\in \Gen_{\bar g^*}(V,M)\losp (f\circ h'=h\circ g)$
will be denoted by $\Gen_{\bar g^*,h}(V',P)$.

\subsection[The Second Order Logic and the Structure
$\langle \protect\Cn, \protect\ring\rangle$]{The Second
Order Logic and the Structure
$\boldsymbol{\langle \Cn, \ring\rangle}$}\label{ss2.7}
Consider the structure $\langle \Cn, \ring\rangle$, consisting
of the class~$\Cn$ of all cardinal numbers and the ring $\ring$
with usual ring relations $+$~and~$\circ$.
The \emph{second-order logic} of the
structure ($L_2(\langle \Cn,\ring\rangle)$)
allows to use in formulas arbitrary predicate
symbols of the form
$$
P_{\lambda_1,\dots,\lambda_k}(c_1,\dots,c_k;v_1,\dots, v_n),
$$
where $\lambda_1,\dots,\lambda_k$ are fixed cardinal numbers,
$c_1,\dots,c_k$ are variables for elements
from $\lambda_1 ,\dots ,\lambda_k$,
respectively, and $v_1,\dots,v_n$ are variables for ring elements.

Therefore, in formulas of this language we can use
the following subformulas.
\begin{enumerate}
\renewcommand{\labelenumi}{\theenumi.}
\item
$\forall r\in \ring$.
\item
$\exists r\in \ring$.
\item
$\forall \varkappa \in \Cn$.
\item
$\exists \varkappa \in \Cn$.
\item
$\forall \alpha\in \varkappa$,
where $\varkappa$
is either a~free variable of the formula~$\varphi$
or is defined in the formula~$\varphi$ before~$\alpha$
(with the help of the subformula $\forall \varkappa \in \Cn$
or $\exists \varkappa \in \Cn$).
\item
$\exists \alpha \in \varkappa$,
where $\varkappa$
is either a~free variable of the formula~$\varphi$
or is defined in the formula~$\varphi$ before~$\alpha$
(with the help of the subformula $\forall \varkappa \in \Cn$
or $\exists \varkappa \in \Cn$).
\item
$r_1=r_2+r_3$, $r_1=r_2\times r_3$, $r_1=r_2$,
where each of the variables $r_1$, $r_2$, and~$r_3$
is either a~free variable of the formula~$\varphi$ or defined
in the formula~$\varphi$ before (with the help of the subformula
$\forall r_i\in \ring$ or $\exists r_i\in \ring$).
\item
$\varkappa_1=\varkappa_2$, where each of the variables
$\varkappa_1$,~$\varkappa_2$
is either a~free variable of the formula~$\varphi$
or defined in the formula~$\varphi$ before
(with the help of the subformula
$\forall \varkappa_i\in \Cn$ or
$\exists \varkappa_i\in \Cn$).
\item
$\alpha_1=\alpha_2$, where each of the variables $\alpha_1$,~$\alpha_2$
is either a~free variable of the formula~$\varphi$
or defined in the formula~$\varphi$
(with the help of the subformula $\forall \alpha_i\in \varkappa_i$
or $\exists \alpha_i\in \varkappa_i$).
\item\sloppy
$\forall P_{\varkappa_1,\dots\ttm,\varkappa_k}
(c_1,\dots\ttm,c_k;v_1,\dots\ttm,v_n)$,
$\exists P_{\varkappa_1,\dots\ttm,\varkappa_k}
(c_1,\dots\ttm,c_k;v_1,\dots\ttm,v_n)$,
where each of the variables $\varkappa_1,\dots\ttm,\varkappa_k$
is either a~free
variable of the formula~$\varphi$ or defined in the formula~$\varphi$
before (with the help of the subformula $\forall \varkappa_i\in \Cn$
or $\exists \varkappa_i\in \Cn$).
\item
$P_{\varkappa_1,\dots,\varkappa_k}(\alpha_1,\dots,\alpha_k;r_1,\dots,r_n)$,
where each of the variables
$\varkappa_1,\dots,\varkappa_k$, $\alpha_1,\dots,\alpha_k$,
$r_1,\dots,r_n$ and also the ``predicative variable''
$$
P_{\varkappa_1,\dots,\varkappa_k}(c_1,\dots,c_k;v_1,\dots,v_n)
$$
is either a~free variable of the formula~$\varphi$
or defined in the formula~$\varphi$ before
(with the help of the subformulas $\forall \varkappa_i\in \Cn$,
$\exists \varkappa_i\in \Cn$, $\forall \alpha_i\in \varkappa_i$,
$\exists \alpha_i\in \varkappa_i$, $\forall r_i\in \ring$,
$\exists r_i\in \ring$,
$\forall P_{\varkappa_1,\dots,\varkappa_k}(c_1,\dots,c_k;v_1,\dots,v_n)$,
or $\exists P_{\varkappa_1,\dots,\varkappa_k}(c_1,\dots,c_k;v_1,\dots,v_n)$),
$\varkappa_i$ is introduced in the formula before~$\alpha_i$ for every
$i=1,\dots,k$, and
$P_{\varkappa_1,\dots,\varkappa_k}(c_1,\dots,c_k;v_1,\dots,v_n)$
is introduced after all $\varkappa_1,\dots,\varkappa_k$.\sloppy
\end{enumerate}

\begin{theorem}\label{my2_7}\label{th5_2.7}
Let $R$~and~$S$ be rings.
Suppose that there exists a~sentence~$\psi$
of the language $L_2(\langle \Cn,\ring\rangle)$
which is true in the ring~$R$, false in any ring
similar to~$R$, and not equivalent to it in the language
$L_2(\langle \Cn,\ring \rangle )$.
If the categories $\textup{mod-}R$ and $\textup{mod-}S$ are
elementarily equivalent, then there exists a~ring~$S'$ which is similar
to~$S$ and such that the structures
$\langle \Cn ,R\rangle $ and $\langle \Cn,S'\rangle$ are equivalent in the
logic~$L_2$.
\end{theorem}

\begin{proof}
Suppose that some progenerator~$P$
in the category $\textup{mod-}T$, where $T$ is some ring, is fixed.
Then, according to the previous sections, we have formulas
defining a~simple module~$M$ which corresponds to the module~$P$,
modules~$M^{(\varkappa)}$ for all $\varkappa\in \Cn$,
modules~$M^{(n)}$ for all $n\in \omega$, modules~$M^{(\alpha)}$
for infinite $\alpha\in \Cn$, almost free modules~$V^\varkappa$ of rank
$\varkappa\in \Cn$, $\varkappa \in \omega$, $\varkappa\ge \omega$, and,
besides, for every module~$M^{(\varkappa)}$ (or~$V^\varkappa$)
its generating sets $\Gen_{g^*}(M^{(\varkappa)},M)$
(or $\Gen_{g^*}(V^\varkappa, R)$).
Further (see Sec.~2.3), for any $f,g\in \Mor(P,P)$ we suppose,
that their sum $f\oplus g\in \Mor(P,P)$ and product
$f\otimes g\in \Mor(P,P)$ are known.

Consider any arbitrary sentence~$\varphi$ in the language
$L_2(\langle \Cn,\ring \rangle)$.
As it was shown before, this sentence can contain
the following subformulas.
\begin{enumerate}
\renewcommand{\labelenumi}{\theenumi.}
\item
$\forall r\in \ring$.
\item
$\exists r\in \ring$.
\item
$\forall \varkappa \in \Cn$.
\item
$\exists \varkappa \in \Cn$.
\item
$\forall \alpha\in \varkappa$.
\item
$\exists \alpha \in \varkappa$.
\item
$r_1=r_2+r_3$.
\item
$r_1=r_2\cdot r_3$.
\item
$r_1=r_2$.
\item
$\varkappa_1=\varkappa_2$.
\item
$\alpha_1=\alpha_2$.
\item
$\forall P_{\varkappa_1,\dots,\varkappa_k}
(c_1,\dots,c_k; v_1,\dots,v_n)$.
\item
$\exists P_{\varkappa_1,\dots,\varkappa_k}
(c_1,\dots,c_k;v_1,\dots,v_n)$.
\item
$P_{\varkappa_1,\dots,\varkappa_k}
(\alpha_1,\dots,\alpha_k;r_1,\dots,r_n)$.
\end{enumerate}

We shall transform this sentence into a~sentence $\tilde \varphi_P$
(which depends on the fixed module~$P$) of the first order language of the
category theory by the following algorithm.

1. The subformula $\forall r\in \ring$ is transformed into the formula
$\forall f_r\in \Mor(P,P)$, i.e.,
every element of the ring $\ring$ corresponds
to an element of the ring $\Endom_T(P)$.

2. The subformula $\exists r\in \ring$ is transformed into the formula
$\exists f_r\in \Mor(P,P)$.

3. The subformula $\forall \varkappa \in \Cn$
is transformed into the subformula
$\forall X_\varkappa \in \Obj\losp
\forall \bar g_\varkappa^*\losp
(X_\varkappa =M^{(\varkappa)} \logic\land
\Gen_{\bar g_\varkappa^*}(X_\varkappa,M) \Rightarrow \dots)$,
i.e., every element $\varkappa\in \Cn$ corresponds to some module of the
form~$M^{(\varkappa)}$ for the simple module~$M$
(we have already mentioned that there exists a~natural
identification of the class $\Cn$ and the class of all
direct sums of the module~$M$), and immediately
the set $\Gen_{\bar g_{\varkappa}^*}(M^{(\varkappa)},M)$
of projectors from~$M^{(\varkappa)}$ onto~$M$ becomes fixed.

4. The subformula $\exists \varkappa \in \Cn$
is transformed into the subformula
$\exists X_\varkappa \in \Obj\losp
\exists \bar g_\varkappa^*\losp
(X_\varkappa =M^{(\varkappa)} \logic\land
\Gen_{\bar g_{\varkappa}^*}(X_\varkappa,M) \logic\land \ldots)$.

5. The subformula $\forall \alpha \in \varkappa$
is transformed into the subformula
$\forall f^\alpha_{X_\varkappa}\in
\Gen_{\bar g_\varkappa^*}(X_\varkappa,M)$,
i.e., elements of sets~$\varkappa$
correspond to those mappings of the set
$\Gen_{\bar g_\varkappa^*} (M^{(\varkappa)},M)$ which contain exactly
$\varkappa$ linearly independent projectors.

6. The subformula $\exists \alpha \in \varkappa$
is transformed into the subformula
$\exists f^\alpha_{X_\varkappa}\in \Gen_{\bar g_\varkappa^*}\losp
(X_\varkappa,M)$.

7. The subformula $r_1=r_2+r_3$ is transformed into the subformula
$f_{r_1}=f_{r_2}\oplus f_{r_3}$, i.e., the sum of elements of the ring
$\ring$ corresponds to the sum of elements of the ring $\Endom_T(P)$.

8. The formula $r_1=r_2\cdot r_3$ is transformed into the formula
$f_{r_1}=f_{r_2}\otimes f_{r_3}$, i.e.,
the product of elements of the ring $\ring$
corresponds to the product of elements of the ring $\Endom_T(P)$.

9. The subformula $r_1=r_2$ is transformed into the subformula
$f_{r_1}=f_{r_2}$, i.e., equal elements of the ring $\ring$ correspond
to equal elements of the ring $\Endom_T(P)$.

10. The subformula $\varkappa_1=\varkappa_2$ is transformed
into the subformula
$\exists g_{\varkappa_1,\varkappa_2}
\in \Mor(X_{\varkappa_1},X_{\varkappa_2})\losp
(\text{$g$ is}$ an isomorphism),
i.e., equal sets of the class~$\Cn$ correspond
to isomorphic modules of the form $M^{(I)}$~and~$M^{(J)}$,
i.e., such modules that $|I|=|J|=\varkappa$.

11. The subformula $\alpha_1=\alpha_2$ for
$\alpha_1,\alpha_2\in \varkappa$ is transformed into the subformula
$f^{\alpha_1}_{X_\varkappa}=f^{\alpha_2}_{X_\varkappa}$,
and the subformula $\alpha_1=\alpha_2$ for
$\alpha_1\in \varkappa_1$, $\alpha_2\in \varkappa_2$, and
$\varkappa_1 \ne \varkappa_2$ is transformed into the subformula
$f_{X_{\varkappa_1}}^{\alpha_1}=f^{\alpha_2}_{X_{\varkappa_2}}\circ g$,
i.e., equal elements of the set $\varkappa\in \Cn$ are mapped
to corresponding to each other projections in
$\Gen_{\bar g_{\varkappa_1}^*}(M^{(I)},M)$ and
$\Gen_{\bar g_{\varkappa_2}^*}(M^{(J)},M)$,
and the correspondence is fixed by the isomorphism between
$M^{(I)}$ and $M^{(J)}$.

Before the last three transformations we shall
introduce the following new formulas.

For every $f_t\in \Gen_{\bar g_\varkappa^*}(M^{(\varkappa)},M)$ by $f_t'$
we shall denote the corresponding mapping from
$\Gen_{\bar g_\varkappa^*}(V^\varkappa,P)$, by $\bar f_t$ we shall denote
a~mapping from $\Mor(M,M^{(\varkappa)})$ such that $f_t\circ \bar f_t=1_M$,
by $\bar f_t'$ we shall denote a~mapping from $\Mor(P,V^\varkappa)$ such that
$f_t'\circ \bar f_t'=1_P$. Given a~mapping
$f\in \Mor(V^\varkappa,V^\varkappa)$
we shall write
$$
f\in \Ring_{\bar g_\varkappa^*}(V^\varkappa)
$$
if
$$
\forall f_t',f_s'\in \Gen_{\bar g_\varkappa^*, h}(V^\varkappa,P)\losp
(f_t'\ne f_s' \Rightarrow f_t'\circ f\circ \bar f_s'=0).
$$
Given a~mapping $f\in \Mor(M^{(\varkappa_1)},M^{(\varkappa_2)})$
we shall write
$$
f\in \Sets_{\bar g_{\varkappa_1}^*,\bar g_{\varkappa_2}^*}
(M^{(\varkappa_1)},M^{(\varkappa_2)})
$$
if
$$
\forall f_t\in \Gen_{\bar g_{\varkappa_1}^*}(M^{(\varkappa_1)},M)\losp
\forall f_s\in \Gen_{\bar g_{\varkappa_2}^*}(M^{(\varkappa_2)},M)\losp
(f_s\circ f\circ \bar f_t=1_M \logic\lor f_s\circ f\circ \bar f_t=0).
$$

Therefore the elements from $\Ring_{\bar g_\varkappa^*}(V^\varkappa)$
are those endomorphisms of the module~$V^\varkappa$
that are diagonal in some initially fixed basis,
and so these endomorphisms can be considered as mappings from~$\varkappa$
into the ring $\Endom_T(P)$, mapping every $\alpha\in \varkappa$
to the element on the diagonal at position~$\alpha$.
Elements from $\Sets_{\bar g_{\varkappa_1}^*,\bar g_{\varkappa_2}^*}
(M^{(\varkappa_1)},M^{(\varkappa_2)})$ are those morphisms
from $M^{(\varkappa_1)}$ and $M^{(\varkappa_2)}$ that
in the given fixed basis have matrices
consisting only of zeros and units. These matrices
can be understood as correspondences~$F$
between the sets $\varkappa_1$~and~$\varkappa_2$,
where a~pair $\langle x,y\rangle$
belongs to the correspondence~$F$ if and only if
the intersection of the row with index~$x$
and the column with index~$y$ in this matrix is a~unit.

We use these remarks for the remaining transformations.

12. Let $\varkappa=\max\{ \varkappa_1,\dots,\varkappa_k,|\ring|\}$.
Then the subformula
$$
\forall P_{\varkappa_1,\dots,\varkappa_k} (c_1,\dots,c_k;v_1,\dots,v_n)
$$
is transformed into the subformula
\begin{multline*}
\forall f_P^{c_1}\in \Sets_{\bar g_{\varkappa}^*,\bar g_{\varkappa_1}^*}
(M^{(\varkappa)}, M^{(\varkappa_1)})\ldots
\forall f_P^{c_k} \in \Sets_{\bar g_{\varkappa}^*,\bar g_{\varkappa_k}^*}
(M^{(\varkappa)}, M^{(\varkappa_k)})\\
\forall f_P^{v_1}\in \Ring_{\bar g_\varkappa^*} (V^{\varkappa})\ldots
\forall f_P^{v_n}\in \Ring_{\bar g_\varkappa^*} (V^{\varkappa}),
\end{multline*}
i.e., every predicate symbol of the form
$P_{\varkappa_1,\dots,\varkappa_k}(c_1,\dots,c_k;v_1,\dots,v_n)$
corresponds to $k$~mappings for sets
$\varkappa_1,\dots,\varkappa_k$ and $n$~mappings for elements
of the ring which are connected to each other
with the help of the module~$M^{(\varkappa)}$.

13. The subformula
$$
\exists P_{\varkappa_1,\dots,\varkappa_k}(c_1,\dots,c_k;v_1\dots,v_n)
$$
is transformed into the subformula
\begin{multline*}
\exists f_P^{c_1}\in \Sets_{\bar g_\varkappa^*,\bar g_{\varkappa_1}^*}
(M^{(\varkappa)}, M^{(\varkappa_1)}) \ldots
\exists f_P^{c_k}\in \Sets_{\bar g_\varkappa^*,\bar g_{\varkappa_k}^*}
(M^{(\varkappa)}, M^{(\varkappa_k)})\\*
\exists f_P^{v_1} \in \Ring_{\bar g_\varkappa^*} (V^{\varkappa})\ldots
\exists f_P^{v_n}\in \Ring_{\bar g_\varkappa^*} (V^{\varkappa}).
\end{multline*}

14. The subformula
$$
P_{\varkappa_1,\dots,\varkappa_k} (\alpha_1,\dots,\alpha_k;r_1,\dots,r_n)
$$
is transformed into the subformula
\begin{multline*}
\exists f\in \Gen_{\bar g_\varkappa^*}(M^{(\varkappa)}, M)\\
(f_{X_{\varkappa_1}}^{\alpha_1}\circ f_P^{c_1}\circ \bar f=1
\logic\land \dots \logic\land
f_{X_{\varkappa_k}}^{\alpha_k}\circ f_P^{c_k}\circ \bar f=1
\logic\land
f'\circ f_P^{v_1}\circ \bar f'=f_{r_1}
\logic\land \dots \logic\land
f'\circ f_P^{v_n}\circ \bar f'=f_{r_n}).
\end{multline*}

Let now some sentence~$\varphi$ be true in the model
$\langle \Cn, \Endom_T P\rangle$.
Let all bound variables of the sentence~$\varphi$
be contained in the set
$x_1,\dots,x_q$ (where $x_1,\dots,x_q$ is either a~variable for
ring elements, or for elements of the class~$\Cn$,
or for elements of some~$\varkappa\in \Cn$,
or a~predicate variable). Since the sentence~$\varphi$
is true in the model $\langle \Cn,\Endom_T P\rangle$,
there exists some sequence $y_1,\dots,y_q$ of elements
of this model such that the sentence~$\varphi$ holds on it.
Transform the sequence $y_1,\dots, y_q$
of elements of the model $\langle \Cn,\Endom_T P\rangle$
into a~sequence $z_1,\dots,z_s$
of elements of the model $\textup{mod-}T$.

If $y_i\in \Endom_T(P)$, then we transform the element~$y_i$
to the element $z_i \prisv y_i=f_{y_i}\in \Mor(P,P)$.

If $y_i\in \Cn$ and $y_i=\varkappa$, then transform $y_i$
to the pair $z_i^{(1)}\prisv M^{(\varkappa)}\in \Obj$ and
$z_i^{(2)}\prisv \bar g_\varkappa^*$ such that
$\Gen_{\bar g_\varkappa^*}(M^{(\varkappa)},M)$.

If $y_i\in \varkappa$ and $y_i=\alpha$, where $\alpha$ is an ordinal
number, then transform $y_i$ to
$z_i \prisv f^\alpha\in \Gen_{\bar g_\varkappa^*}(M^{(\varkappa)},M)$,
i.e., to the projection from this set having the index~$\alpha$.

If $y_i=P_{\varkappa_1,\dots,\varkappa_k}(c_1,\dots,c_k;v_1,\dots,v_n)$,
i.e., it is a~relation $\bar P$ on the set
$$
\varkappa_1\times \dots \times \varkappa_k\times
\Endom_TP\times \dots \times \Endom_TP,
$$
then we shall set
$\varkappa \prisv \max\{ \varkappa_1,\dots,\varkappa_k, |\Endom_TP|\}$
and transform $y_i$ to a~sequence
$z_i^1,\dots,z_i^k;z_i^{k+1},\dots, z_i^{k+n}$ of morphisms from the sets
$$
\Sets_{\bar g_\varkappa^*, \bar g_{ \varkappa_1}^*}
(M^{(\varkappa)},M^{(\varkappa_1)}),\dots,
\Sets_{\bar g_\varkappa^*,\bar g_{\varkappa_k}^*}
(M^{(\varkappa)},M^{(\varkappa_k)}),
\Ring_{\bar g_\varkappa^*}(V^\varkappa),\dots,
\Ring_{\bar g_\varkappa^*}(V^\varkappa) ,
$$
respectively, such that
$\langle \alpha_1,\dots,\alpha_k;r_1,\dots,r_n\rangle\in \bar P$
if and only if there exists $\alpha\in \varkappa$ such that in every
matrix $z_i^l$, where $1\le l\le k$, the intersection of the column
with index~$\alpha$ and the row with index~$\alpha_l$
is a~unit, and in every matrix $z_i^l$,
where $k< l\le k+n$,
the element on the diagonal at position~$\alpha$
is~$r_l$.

Therefore, we have a~new sequence $z_1,\dots,z_s$. We show that
the sentence $\tilde \varphi_P$ holds
on this sequence in the model $\textup{mod-}T$.

We shall prove this by induction by the length of the formula.

1. If the formula has the form
$$
r_1=r_2+r_3,
$$
then its transformation has the form
$$
f_{r_1}=f_{r_2}\oplus f_{r_3},
$$
and $r_1=r_2+r_3$ in $\Endom_TP$ if and only if
$f_{r_1}=f_{r_2}\oplus f_{r_3}$ in $\Mor_T(P,P)$
because the rings $\Endom_TP$ and $\Mor_T(P,P)$ are isomorphic. Thus
$$
\langle \Cn,\Endom_T P\rangle_{L_2}\vDash r_1=r_2+r_3
$$
if and only if
$$
\textup{mod-}T\vDash f_{r_1}=f_{r_2}\oplus f_{r_3}.
$$

2. The proof in the case of formulas $r_1=r_2\cdot r_3$
and $r_1=r_2$ is similar to the previous one.

3. If the formula has the form
$$
\varkappa_1=\varkappa_2,
$$
then its transformation has the form
$$
\exists g\in \Mor(M^{(\varkappa_1)},M^{(\varkappa_2)})\losp
(\text{$g$ is an isomorphism}).
$$
If the cardinal numbers $\varkappa_1$~and~$\varkappa_2$
coincide, then the modules $M^{(\varkappa_1)}$ and $M^{(\varkappa_2)}$
are isomorphic, and if the modules $M^{(I)}$~and~$M^{(J)}$
are isomorphic, then $|I|=|J|$. Therefore
$$
\langle \Cn, \Endom_TP\rangle_{L_2}\vDash \varkappa_1=\varkappa_2
$$
if and only if
$$
\textup{mod-}T\vDash
\exists g_{\varkappa_1,\varkappa_2}\in
\Mor (M^{(\varkappa_1)},M^{(\varkappa_2)})\losp
(\text{$g$ is an isomorphism}).
$$

4. The proof of a~similar statement about the formula $\alpha_1=\alpha_2$
is the same.

5. If the formula has the form
$$
P_{\varkappa_1,\dots,\varkappa_k}(\alpha_1,\dots,\alpha_k;r_1,\dots,r_n)
$$
and its transformations has the form
$$
\tilde P_{\varkappa_1,\dots,\varkappa_k}
(\alpha_1,\dots,\alpha_k; r_1,\dots,r_n)_P
$$
and, further,
$$
\langle \Cn,\Endom_T P\rangle_{L_2}\vDash
P_{\varkappa_1,\dots,\varkappa_k}(\alpha_1,\dots,\alpha_k;r_1,\dots,r_n),
$$
then for the sequence
$$
\langle \alpha_1,\dots,\alpha_k,r_1,\dots,r_n\rangle
\in \varkappa_1\times \dots \times \varkappa_k\times
\Endom_T P\times \dots \times \Endom_T P
$$
we have
$$
\langle \alpha_1,\dots,\alpha_k;r_1,\dots,r_n\rangle\in \bar P,
$$
where $\bar P$ is the relation corresponding to the predicate
$P_{\varkappa_1,\dots,\varkappa_k}$, i.e.,
$$
\bar P\subset \varkappa_1\times \dots \times \varkappa_k\times
\Endom_TP\times \dots \times \Endom_TP.
$$
This relation is a~set of sequences that has cardinality at most
$$
|\varkappa_1\times \dots\times \varkappa_k\times |T|\times \dots\times |T||
\le |\varkappa\times \dots \times \varkappa|=\varkappa.
$$
Therefore all sequences from~$\bar P$ can be enumerated by elements
of~$\varkappa$. Let $\bar P(\alpha)$ be a~sequence
from~$\bar P$ with the number~$\alpha$ and let it have the form
$\langle \alpha_1,\dots,\alpha_k;r_1,\dots,r_n\rangle$.
Then the $\alpha$-th column of the matrix~$z_i^l$
for $l=1,\dots,k$ will contain~$1$ at position~$\alpha_l$
and $0$ at all other positions, and the $\alpha$-th column of the
matrix~$z_i^l$ for $l=k+1,\dots, k+n$ will contain
$r_{l-n}$ at the $\alpha$-th position
and $0$ at all other positions. Consequently
$$
\langle \Cn,\Endom_TP\rangle_{L_2}\vDash
P_{\varkappa_1,\dots,\varkappa_k} (\alpha_1,\dots,\alpha_k;r_1,\dots,r_n)
$$
if and only if
$$
\textup{mod-}T\vDash
\tilde P_{\varkappa_1,\dots,\varkappa_k}
(\alpha_1,\dots,\alpha_k;r_1,\dots,r_n)_P.
$$
All other parts of induction are proved similarly.

Now we can easily see that the sentence~$\varphi$
holds in the structure $\langle \Cn, \Endom_T(P)\rangle$ if and only if
the corresponding sentence $\tilde \varphi_P$ holds in $\textup{mod-}T$.

According to the condition of the theorem,
the formula
$$
\Select(P) \prisv P\in \Obj \logic\land \Proobr(P) \logic\land
\tilde \psi^P \logic\land \forall P'\in \Obj\losp
(\Proobr(P') \logic\land P'\not\cong P\Rightarrow \neg \tilde \psi^{P'})
$$
is true in $\textup{mod-}R$ only for $P\cong R$.

Let now categories $\textup{mod-}R$ and $\textup{mod-}S$
be elementarily equivalent and $\varphi$ be a~sentence
in the second-order language~$L_2$ of the structure
$\langle \Cn,\ring\rangle$ which is true in $\langle \Cn,R\rangle$.
Then the sentence
${\forall P\in \Obj}\losp\allowbreak
(\Select(P)\Rightarrow \tilde \varphi^P)$
is true in the category $\textup{mod-}R$, and, therefore, in the
category $\textup{mod-}S$.
Thus the sentence~$\varphi$ is true in $\langle \Cn,\Endom_S(P)\rangle$
for every module~$P$ satisfying
the formula $\Select(P)$
in the category $\textup{mod-}S$.
But for all modules~$P$ satisfying
the formula~$\varphi$ the rings of the form $\Endom_S P$ are
equivalent in the logic $\langle \Cn,\ring\rangle$.
Consequently if we set $S'\prisv \Endom_SP$ for some~$P$ satisfying
the formula $\Select(P)$, then we shall have that the sentence~$\varphi$
is true $\langle \Cn,S'\rangle$, and the ring~$S'$
does not depend on the sequence~$\varphi$.
Therefore the structures
$\langle \Cn,R\rangle$ and $\langle \Cn,S'\rangle$ are equivalent
in the logic~$L_2$.
\end{proof}

\subsection{The Inverse Theorem}\label{ss2.8}
Before proving the inverse theorem we introduce different notions
which we shall need later, and transfer them to the language
$L_2(\langle \Cn,\ring\rangle)$.

A~one-place relation $P_{\varkappa_1}(c)$ will be called
a~\emph{subset of the cardinal number}~$\varkappa_1$. The set
$\{ \alpha\in \varkappa_1\mid P_{\varkappa_1}(\alpha)\}$
will be denoted by~$P_{\varkappa_1}$. We shall use the
notation $\alpha\in P_{\varkappa_1}$ for it.

A~one-place relation $P(v)$ will be called
a~\emph{subset of the ring} $\ring$, and, similarly to the previous
notation, we shall use the notation $r\in P$.

Any two-place relation $F_{\varkappa_1,\varkappa_2}(c_1,c_2)$
(or $F_{\varkappa_1}(c_1,v_1)$, or $F(v_1,v_2)$) will be called
a~\emph{correspondence}
between cardinal numbers $\varkappa_1$ and $\varkappa_2$ (or
between a~cardinal number~$\varkappa_1$ and the ring, or in the ring).
We shall use the notation
$\langle \alpha_1,\alpha_2\rangle \in P_{\varkappa_1,\varkappa_2}$
($\langle \alpha_1,v_1\rangle \in P_{\varkappa_1}$
or $\langle v_1,v_2\rangle \in P$) for the formula
$P_{\varkappa_1,\varkappa_2}(\alpha_1,\alpha_2)$ (and so on).

A~correspondence $F_{\varkappa_1,\varkappa_2} (c_1,c_2)$
(or $F_{\varkappa_1}(c_1,v_1)$, or $F(v_1,v_2)$) for which the formula
$$
\forall \alpha\in \varkappa_1\losp \exists \beta\in \varkappa_2\losp
(\langle \alpha,\beta\rangle\in F_{\varkappa_1,\varkappa_2})
\logic\land \forall \alpha \in \varkappa_1\losp
\forall \beta_1,\beta_2 \in \varkappa_2\losp
(\langle \alpha,\beta\rangle \in F_{\varkappa_1,\varkappa_2}
\logic\land \langle \alpha,\beta_2\rangle \in F_{\varkappa_1,\varkappa_2}
\Rightarrow \beta_1=\beta_2)
$$
holds (similarly for other types of correspondences) is called
a~\emph{mapping from a~cardinal number~$\varkappa_1$
into a~cardinal number~$\varkappa_2$}
(respectively, from a~cardinal number~$\varkappa_1$ into the ring, or
from the ring into itself).
The fact that $F_{\varkappa_1,\varkappa_2}$ ($F_{\varkappa_1}$ or $F$)
is a~mapping will be denoted by $\Func(F_{\varkappa_1,\varkappa_2})$
($\Func(F_{\varkappa_1})$ or $\Func(F)$).

A~mapping $F_{\varkappa_1,\varkappa_2} (c_1,c_2)$
(or $F_{\varkappa_1}(c_1,v_1)$, or $F(v_1,v_2)$) for which the formula
$$
\forall \beta\in \varkappa_2\losp \exists \alpha\in \varkappa_2\losp
(\langle \alpha,\beta\rangle\in F_{\varkappa_1,\varkappa_2})
$$
holds (similarly for other types of mappings) is called
\emph{surjective}
(notation: $\Surj(F)$, or $\Surj(F_{\varkappa_1})$, or $\Surj(F)$).

A~mapping $F_{\varkappa_1,\varkappa_2}(c_1,c_2)$
(or $F_{\varkappa_1}(c_1,v_1)$, or $F(v_1,v_2)$) for which the formula
$$
\forall \alpha_1,\alpha_2\in \varkappa_1\losp
\forall \beta\in \varkappa_2\losp
(\langle \alpha_1,\beta\rangle \in F_{\varkappa_1,\varkappa_2}
\logic\land \langle \alpha_2,\beta\rangle \in F_{\varkappa_1,\varkappa_2}
\Rightarrow \alpha_1=\alpha_2)
$$
holds (similarly for other types of mappings) is called
\emph{injective} (notation: $\Inj(F_{\varkappa_1,\varkappa_2})$,
or $\Inj(F_{\varkappa_1})$, or $\Inj(F)$).

A~mapping which is simultaneously surjective and injective is called
\emph{bijective} (notation: $\Bij(F_{\varkappa_1,\varkappa_2})$, or
$\Bij(F_{\varkappa_1})$, or $\Bij(F)$).

For a~given mapping $F_{\varkappa_1,\varkappa_2}(c_1,c_2)$
(or $F_{\varkappa_1}(c_1,v_1)$, or $F(v_1,v_2)$) the \emph{inverse} mapping
is the mapping $F_{\varkappa_1,\varkappa_2}'(c_1,c_2)$
(or $F_{\varkappa_1}'(c_1,v_1)$, or $F'(v_1,v_2)$)
satisfying the formula
$$
\forall \alpha\in \varkappa_1\losp
\forall \beta\in \varkappa_2\losp
(\langle \alpha,\beta\rangle \in F_{\varkappa_1,\varkappa_2}
\Leftrightarrow
\langle \beta,\alpha\rangle\in F_{\varkappa_1,\varkappa_2}').
$$

The \emph{domain} of a~correspondence $F_{\varkappa_1,\varkappa_2}(c_1,c_2)$
(or $F_{\varkappa_1}(c_1,v_1)$, or $F(v_1,v_2)$) is the set
$A_{\varkappa_1}\subset \varkappa_1$ ($A\subset \ring$) satisfying
the formula
$$
\forall \alpha\in \varkappa_1\losp
(\alpha\in A_{\varkappa_1}\Leftrightarrow \exists \beta\in \varkappa_2\losp
\langle \alpha,\beta\rangle\in F_{\varkappa_1,\varkappa_2}).
$$
The domain is denoted by $\Dom(F_{\varkappa_1,\varkappa_2})$.

The \emph{image} of a~correspondence
$F_{\varkappa_1,\varkappa_2} (c_1,c_2)$
(or $F_{\varkappa_1}(c_1)$, or $F(v_1,v_2)$) is the set
$A_{\varkappa_2}\subset \varkappa_2$ ($A\subset \ring$)
satisfying the formula
$$
\forall \beta\in \varkappa_2\losp
(\beta\in A_{\varkappa_2} \Leftrightarrow
\exists \alpha \in \varkappa_1\losp
\langle \alpha,\beta\rangle\in F_{\varkappa_1,\varkappa_2})
$$
(notation: $\Rng(F_{\varkappa_1,\varkappa_2})$).

A~cardinal number $\mu\in \Cn$ will be called \emph{infinite}
(notation: $\mu\in \Inf$ or $\Inf(\mu)$) if it satisfies
the formula
$$
\exists F_{\mu,\mu}(c_1,c_2)\losp
(\Inj(F_{\mu,\mu}) \logic\land \Rng(F_{\mu,\mu})\ne \mu).
$$

A~cardinal number $\mu\in \Cn$ will be called
\emph{finite} (notation: $\mu\in \Fin$ or $\Fin(\mu)$) if
$\mu\notin \Inf$.

The \emph{cardinality} of a~set
$M_{\varkappa}\subset \varkappa$ ($M\subset \ring$) is
the cardinal number $\mu\in \Cn$ satisfying the formula
$$
\exists F_{\mu,\varkappa}(c_1,c_2)\losp
(\Inj(F_{\mu,\varkappa}) \logic\land
\Dom(F_{\mu,\varkappa})=\mu
\logic\land \Rng(F_{\mu,\varkappa})=M_\varkappa).
$$
The cardinality of a~set $M_\varkappa$ ($M$)
will be denoted by~$|M_\varkappa|$
($|M|$).

A~set~$M_\varkappa$ ($M$) will be called \emph{finite}
if its cardinality is a~finite cardinal number.

Consider some finite set $M_\varkappa$ ($M$).
A~correspondence $\bar M_{\varkappa,\varkappa} (c_1,c_2)$
($\bar M(v_1,v_2)$) will be called a~\emph{relation of
consecutive order} on this set if
\begin{align*}
& \forall \alpha_1,\alpha_2,\alpha_3\in M_\varkappa\losp
(
(\langle \alpha_1,\alpha_2\rangle \in \bar M_{\varkappa,\varkappa}
\logic\land
\langle \alpha_1,\alpha_3\rangle\in \bar M_{\varkappa,\varkappa}
\Rightarrow \alpha_2=\alpha_3)\\
& \quad {}\logic\land
(\langle \alpha_1,\alpha_3\rangle \in \bar M_{\varkappa,\varkappa}
\logic\land
\langle \alpha_2,\alpha_3\rangle\in \bar M_{\varkappa,\varkappa}
\Rightarrow \alpha_1=\alpha_2)
)\\
& \quad {}\logic\land
\exists \alpha_{\min}, \alpha_{\max}\in M_\varkappa\losp
\forall \alpha\in M_\varkappa\losp
(
(\alpha=\alpha_{\max} \logic\lor
\exists \alpha'\in M_\varkappa\losp
(\langle \alpha,\alpha'\rangle \in \bar M_{\varkappa,\varkappa}))\\
& \quad {}\logic\land
(\alpha=\alpha_{\min} \logic\lor
\exists \alpha'\in M_\varkappa\losp
(\langle \alpha',\alpha\rangle \in \bar M_{\varkappa,\varkappa}))
)
\logic\land
\forall \alpha\in M_\varkappa\losp
(\langle \alpha_{\max},\alpha\rangle \notin \bar M_{\varkappa,\varkappa}
\logic\land
\langle \alpha,\alpha_{\min} \rangle \notin \bar M_{\varkappa,\varkappa}).
\end{align*}

The property of a~predicate $\bar M_{\varkappa,\varkappa}(c_1,c_2)$
($\bar M(v_1,v_2)$) to be a~consecutive order on a~set
$M_\varkappa$ ($M$)
will be denoted by $\Next_{M_\varkappa}(\bar M_{\varkappa,\varkappa})$
($\Next_M(\bar M)$).

If $\bar M_{\varkappa,\varkappa}(c_1,c_2)$ ($\bar M(v_1,v_2)$)
is a~fixed consecutive order on a~set $M_\varkappa$ ($M$),
then for $\alpha_1,\alpha_2\in \varkappa$ ($r_1,r_2\in \ring$) such that
$\langle \alpha_1,\alpha_2\rangle \in \bar M_{\varkappa,\varkappa}$
($\langle r_1,r_2\rangle \in \bar M$) we shall write
$\alpha_2=\alpha_1\oplus_{\bar M} 1$ ($r_2=r_1\oplus_{\bar M} 1$).

Let $M\subset \ring$ be some subset of the ring $\ring$.
By $\smash[b]{\sum\limits_{r\in M} r}$
we shall denote the element~$\bar r$ of the ring
$\ring$ satisfying the formula
\begin{align*}
& \exists \bar M (v_1,v_2)\losp \exists S(v_1,v_2)\losp
(
\Next_M (\bar M) \logic\land \Bij(S) \logic\land
\langle r_{\min} (\bar M), r_{\min}(\bar M)\rangle \in S\\
& \quad \logic\land
\forall r_1,r_2,r_3,r_4\in M\losp
(
r_2=r_1\oplus_{\bar M} 1 \logic\land
\langle r_1,r_3\rangle \in S \logic\land
\langle r_2,r_4\rangle \in S \Rightarrow r_4=r_3+r_2
)\\
& \quad {}\logic\land\langle r_{\max}(\bar M),\bar r\rangle \in S
).
\end{align*}

It is clear that the formula $\sum\limits_{r\in M} r$ introduces
the usual addition in the ring $\ring$.

A~\emph{matrix} of size $\varkappa_1\times \varkappa_2$ is a~relation
$M_{\varkappa_1,\varkappa_2}(c_1,c_2,v_1)$ satisfying the formula
\begin{align*}
& \forall \alpha\in \varkappa_1\losp
\forall \beta\in \varkappa_2\losp
\exists r\in \ring\losp
(
\langle \alpha,\beta,r\rangle \in M_{\varkappa_1,\varkappa_2}
)\\
& \quad \logic\land
\forall \alpha\in \varkappa_1\losp
\forall \beta\in \varkappa_2\losp
\forall r_1,r_2\in \ring\losp
(
\langle \alpha,\beta,r_1\rangle\in M_{\varkappa_1,\varkappa_2}
\logic\land
\langle \alpha,\beta,r_2\rangle \in M_{\varkappa_1,\varkappa_2}
\Rightarrow r_1=r_2
)\\
& \quad \logic\land
\forall \beta\in \varkappa_2\losp
\forall M_{\varkappa_1} \subset \varkappa_1\losp
(
\forall \alpha\in \varkappa_1\losp
(\alpha\in M_{\varkappa_1}\Leftrightarrow
\exists r\in \ring\losp
(\langle \alpha,\beta,r\rangle \in M_{\varkappa_1,\varkappa_2}
\logic\land r\ne 0))
\Rightarrow
|M_{\varkappa_1}|\in \Fin
).
\end{align*}

Relations $M_{\varkappa_1,\varkappa_2} (c_1,c_2;v_1)$ which are matrices
will be denoted by $\Matrix(M_{\varkappa_1,\varkappa_2})$.

\begin{theorem}\label{my2_8}\label{th6_2.8}
If structures $\langle \Cn,R\rangle$ and $\langle \Cn, S\rangle$
are equivalent in the second-order logic~$L_2$, then the
categories $\textup{mod-}R$ and $\textup{mod-}S$ are elementarily
equivalent.
\end{theorem}

\begin{proof}
Consider an arbitrary sentence~$\varphi$ in the first order
language of category theory which is true in the
category $\textup{mod-}R$.

We shall transform it to a~sentence of the second-order language
of the structure $\langle \Cn,R\rangle$.

At the beginning we shall give an informal description of this
transformation.

Every object variable is transformed into a~pair where the first element
is a~cardinal number~$\varkappa$ (which corresponds to the rank
of a~free module over~$R$) and
the second element is a~matrix of size $\varkappa\times \varkappa$
with elements from the ring~$R$ such that
the matrix contains only a~finite number of nonzero elements in every column.
This matrix naturally corresponds
to a~submodule of the module~$R^{(\varkappa)}$
(the columns are the generating elements of this submodule).
We shall associate such a~pair with
a~factormodule of the free module~$R^{(\varkappa)}$ by this
submodule.

Every morphism variable is transformed into a~triplet
consisting of two objects encoded as described above
(we shall denote the corresponding cardinal numbers by
$\varkappa$~and~$\varkappa'$ and the corresponding submodules
by $A$~and~$A'$)
and of a~matrix of size $\varkappa\times \varkappa'$, defining
a~linear mapping from~$R^{(\varkappa)}$ into~$R^{(\varkappa')}$ such that
the image of the submodule~$A$ is a~submodule of the module~$A'$.

Every identity morphism is transformed into a~triplet where
the first and the second components coincide
and the third component
is the identity matrix.

The composition of two morphisms (two triplets) is transformed
into a~triplet where the first object is the first object
of the first triplet, the second object is the second
object of the second triplet,
and the third object is the composition of the matrices
from the first and the second triplet.

Now we shall go on to the formal translation.

We shall perform the following replacements in the sentence~$\varphi$.

1. A~subformula $\forall X\in \Obj$ will be replaced by the subformula
$$
\forall \varkappa_X\in \Cn\losp
\forall P_{\varkappa_X,\varkappa_X}^X(c_1,c_2,v)\losp
(\Matrix(P_{\varkappa_X,\varkappa_X}^X)\Rightarrow \dots).
$$

2. A~subformula $\exists X\in \Obj$ will be replaced by the subformula
$$
\exists \varkappa_X\in \Cn\losp
\exists P_{\varkappa_X,\varkappa_X}^X(c_1,c_2,v)\losp
(\Matrix(P_{\varkappa_X,\varkappa_X}^X) \logic\land \dots).
$$

Now we need to write a~condition for the matrix of a~morphism.
The condition will state that
this matrix moves the first object to the second
one, i.e., all columns of the matrix of the first object
will be transformed
by the action of this matrix
into linear combinations
of columns of the matrix of the second object.
To write this sentence we need to introduce
a~formula expressing the sum of an infinite set
of elements of a~ring if it is known
that only a~finite number of them are nonzero.

For convenience,
given a~matrix $M_{\varkappa_1,\varkappa_2}(c_1,c_2,v_1)$
and fixed $\alpha\in \varkappa_1$ and $\beta\in \varkappa_2$,
we shall denote by
$M_{\varkappa_1,\varkappa_2}(\langle \alpha,\beta\rangle)$
the unique $r\in \ring$
for which $\langle \alpha,\beta,r\rangle\in M_{\varkappa_1,\varkappa_2}$.

Suppose that we have some mapping $F_\varkappa (c,v)$,
whose image is a~subset of the ring $\ring$,
and there exist only a~finite number of $\alpha\in \varkappa$
such that
$\langle \alpha,r\rangle\in F_\varkappa$
for a nonzero $r\in \ring$.
Then by
$$
\sum_{\alpha\in \varkappa} F_\varkappa (\langle \alpha\rangle)
$$
we shall denote the element $r\in \ring$
satisfying the formula
$$
\forall M_\varkappa(c,v)\losp
(
\forall \alpha\in \varkappa\losp \forall r'\in \ring\losp
(
\langle \alpha,r'\rangle \in M_\varkappa \Leftrightarrow
r'\ne 0
\logic\land
\langle \alpha,r'\rangle \in F_\varkappa
)
\Rightarrow
r=\sum_{\alpha\in \Dom(M_\varkappa)} M_\varkappa (\langle \alpha\rangle)
).
$$

Now we are ready to give the translation~3.

3. A~subformula $\forall f\in \Mor$ will be replaced by
the subformula
\begin{align*}
& \forall \varkappa_f\losp
\forall P_{\varkappa_f,\varkappa_f}^f\in \widetilde{\Obj}\losp
\forall \varkappa_f'\losp
\forall {P_{\varkappa_f',\varkappa_f'}^f}'\in \widetilde{\Obj}\losp
\forall Q_{\varkappa_f,\varkappa_f'}^f (c_1,c_2,v)\losp
\biggl(
\Matrix(Q_{\varkappa_f,\varkappa_f'}^f)\\
& \quad {}\logic\land
\forall \beta\in \varkappa_f\losp
\exists S_{\varkappa_f'}(c,v)\losp
\biggl(
\Func(S_{\varkappa_f'})\langle \Dom\rangle
\logic\land
|\Dom(S_{\varkappa_f'})|\in \Fin\\
& \quad {}\logic\land
\forall \gamma\in \varkappa_f'\losp
\biggl(\!
\biggl(
\gamma\in \Dom(S_{\varkappa_f'})
\logic\land
\sum_{\alpha \in \varkappa_f}
Q_{\varkappa_f,\varkappa_f'}^f
(\langle \alpha,\gamma \rangle)\cdot
P_{\varkappa_f,\varkappa_f}^f(\langle \alpha,\beta\rangle)=
\sum_{\xi \in \varkappa_f'}
S(\gamma) \cdot {P_{\varkappa_f',\varkappa_f'}^f}'(\langle \xi,\gamma\rangle)
\biggr)\\
& \quad {}\logic\lor
\biggl(
\gamma\notin \Dom(S_{\varkappa_f'})
\logic\land
\sum_{\alpha\in \varkappa_f}
Q_{\varkappa_f,\varkappa_f'}(\langle \alpha,\gamma\rangle)\cdot
P_{\varkappa_f,\varkappa_f}^f(\langle \alpha,\beta\rangle)=0
\biggr)\!
\biggr)\!
\biggr)
\Rightarrow \ldots
\biggr).
\end{align*}

4. Similarly to the previous case, a~subformula $\exists f\in \Mor$
will be replaced by the subformula
\begin{align*}
& \exists \varkappa_f\losp
\exists P_{\varkappa_f,\varkappa_f}^f\in \widetilde{\Obj}\losp
\exists \varkappa_f'\losp
\exists {P_{\varkappa_f',\varkappa_f'}^f}'\in \widetilde{\Obj}\losp
\forall Q_{\varkappa_f,\varkappa_f'}^f (c_1,c_2,v)\losp
\biggl(
\Matrix(Q_{\varkappa_f,\varkappa_f'}^f)\\
& \quad {}\logic\land
\forall \beta\in \varkappa_f\losp
\exists S_{\varkappa_f'}(c,v)\losp
\biggl(
\Func(S_{\varkappa_f'})\langle \Dom\rangle
\logic\land
|\Dom(S_{\varkappa_f'})|\in \Fin\\
& \quad {}\logic\land
\forall \gamma\in {\varkappa_f'}\losp
\biggl(\!
\biggl(
\gamma\in \Dom(S_{\varkappa_f'})
\logic\land
\sum_{\alpha \in \varkappa_f}
Q_{\varkappa_f,\varkappa_f'}^f(\langle \alpha,\gamma \rangle)\cdot
P_{\varkappa_f,\varkappa_f}^f(\langle \alpha,\beta\rangle)=
\sum_{\xi \in \varkappa_f'} S(\gamma) \cdot
{P_{\varkappa_f',\varkappa_f'}^f}'(\langle \xi,\gamma\rangle)
\biggr)\\
& \quad {}\logic\lor
\biggl(
\gamma\notin \Dom(S_{\varkappa_f'})
\logic\land
\sum_{\alpha\in \varkappa_f}
Q_{\varkappa_f,\varkappa_f'}(\langle \alpha,\gamma\rangle)\cdot
P_{\varkappa_f,\varkappa_f}^f(\langle \alpha,\beta\rangle)=0
\biggr)\!
\biggr)\!
\biggr)
\logic\land \ldots
\biggr).
\end{align*}

5. A~subformula $X=Y$ for $X,Y\in \Obj$ will be replaced by
the subformula
$$
\varkappa_X=\varkappa_Y \logic\land
\forall \alpha,\beta\in \varkappa_X\losp
\forall r\in \ring\losp
(P_{\varkappa_X,\varkappa_X}^X(\alpha,\beta,r)\Leftrightarrow
P_{\varkappa_X,\varkappa_X}^Y(\alpha,\beta,r)),
$$
and the subformula $f=g$ for $f,g\in \Mor$ will be replaced by the formula
\begin{multline*}
\varkappa_f=\varkappa_g \logic\land
\varkappa_f'=\varkappa_g' \logic\land
\forall \alpha_1,\alpha_2 \in \varkappa_f\losp
\forall \beta_1,\beta_2\in \varkappa_f'\losp
\forall r\in \ring\losp
(
(P_{\varkappa_f,\varkappa_f}^f(\alpha_1,\alpha_2,r)\Leftrightarrow
P_{\varkappa_f,\varkappa_f}^g(\alpha_1,\alpha_2,r))\\
{}\logic\land
({P_{\varkappa_f',\varkappa_f'}^f}'(\beta_1,\beta_2,r)\Leftrightarrow
{P_{\varkappa_f',\varkappa_f'}^f}'(\beta_1,\beta_2,r))
\logic\land
(Q_{\varkappa_f,\varkappa_f'}^f(\alpha_1,\beta_1,r)\Leftrightarrow
Q_{\varkappa_f,\varkappa_f'}^g(\alpha_1,\beta_1,r))
).
\end{multline*}

6. A~subformula $f\in \Mor(X,Y)$ for given $f\in \Mor$ and $X,Y\in \Obj$
will be replaced by the formula
\begin{multline*}
\varkappa_f=\varkappa_X \logic\land
\varkappa_f'=\varkappa_Y \logic\land
\forall \alpha_1,\alpha_2\in \varkappa_X\losp
\forall \beta_1,\beta_2\in \varkappa_Y\losp
\forall r\in \ring\\
(P_{\varkappa_X,\varkappa_X}^f(\alpha_1,\alpha_2,r)\Leftrightarrow
P_{\varkappa_X,\varkappa_X}^X(\alpha_1,\alpha_2,r))
\logic\land
({P_{\varkappa_Y,\varkappa_Y}^f}' (\beta_1,\beta_2,r)\Leftrightarrow
P_{\varkappa_Y,\varkappa_Y}^Y(\beta_1,\beta_2,r)).
\end{multline*}

7. A~subformula $f=1_X$ for given $f\in \Mor$ and $X\in \Obj$
will be replaced by the subformula
\begin{multline*}
\varkappa_f=\varkappa_X \logic\land
\varkappa_f'= \varkappa_X \logic\land
\forall \alpha,\beta\in \varkappa_X\losp
\forall r\in \ring\losp
(P_{\varkappa_X,\varkappa_X}^X(\alpha,\beta,r)\Leftrightarrow
P_{\varkappa_X,\varkappa_X}^f(\alpha,\beta,r)\Leftrightarrow
{P_{\varkappa_X,\varkappa_X}^f}'(\alpha,\beta,r))\\
{}\logic\land
\forall \gamma \in \varkappa_X\losp
(Q_{\varkappa_X,\varkappa_X}^f(\gamma,\gamma,1))
\logic\land
\forall \gamma,\eta\in \varkappa_X\losp
(\gamma\ne \eta\Rightarrow Q_{\varkappa_X,\varkappa_X}^f(\gamma,\eta,0)).
\end{multline*}

8. A~subformula $f=g\circ h$ for given $f,g,h\in \Mor$ will be replaced by
the formula
\begin{align*}
& \varkappa_f=\varkappa_h \logic\land
\varkappa_f'=\varkappa_g' \logic\land
\varkappa_h'= \varkappa_g\\
& \quad \logic\land
\forall \alpha_1,\alpha_2\in \varkappa_f\losp
\forall \beta_1,\beta_2\in \varkappa_f'\losp
\forall \gamma_1,\gamma_2\in \varkappa_g\losp
\forall r\in \ring\losp
(
(P_{\varkappa_f,\varkappa_f}^f(\alpha_1,\alpha_2,r)\Leftrightarrow
P_{\varkappa_f,\varkappa_f}^h(\alpha_1,\alpha_2,r))\\
& \quad {}\logic\land
({P_{\varkappa_f',\varkappa_f'}^f}'(\beta_1,\beta_2,r)\Leftrightarrow
{P_{\varkappa_f',\varkappa_f'}^g}'(\beta_1,\beta_2,r))
\logic\land
(P_{\varkappa_g,\varkappa_g}^g(\gamma_1,\gamma_2,r)\Leftrightarrow
{P_{\varkappa_g,\varkappa_g}^h}'(\gamma_1,\gamma_2,r))
)\\
& \quad {}\logic\land
\forall \xi \in \varkappa_f\losp
\forall \eta\in \varkappa_f'\losp
\biggl(
Q_{\varkappa_f,\varkappa_f'}^f(\langle \xi,\eta\rangle)=
\sum_{\alpha\in \varkappa_g}
Q_{\varkappa_g,\varkappa_g'}^g(\langle \alpha,\eta\rangle)\cdot
Q_{\varkappa_h,\varkappa_h'}^h(\langle \xi,\alpha\rangle)
\biggr).
\end{align*}

Thus every sentence~$\varphi$ in the first order logic
of the category theory can be translated to a~sentence~$\tilde \varphi$
of the second-order logic~$L_2$
of the structure $\langle \Cn,\ring\rangle$,
and the algorithm of this translation does not depend on the basic ring.
The sentence~$\varphi$ holds in the category $\textup{mod-}R$
if and only if the sentence~$\tilde \varphi$
holds in the structure $\langle \Cn,R\rangle$.

Consider some sentence~$\varphi$ (or some formula~$\varphi$)
in the first order language of the category theory.

Let all bound (free and bound) variables of the sentence
(formula)~$\varphi$ be contained in the set $x_1,\dots,x_q$
(every $x_l$ is either a~variable for
elements of the class $\Obj$ or for elements of the class $\Mor$).
Consider some sequence of elements of the model $\textup{mod-}R$
$y_1,\dots,y_q$ such that if $x_l$ is a~variable for objects,
then $y_l\in \Obj$ and if $x_l$ is a~variable for
morphisms, then $y_l\in \Mor$.

We shall translate the sequence $y_1,\dots,y_q$ into a~sequence
$z_1,\dots,z_s$ of elements of the model $\langle \Cn,R\rangle_{L_2}$
as follows.

If $y_l\in \Obj$, then $y_l$ is some module over a~ring~$R$.
As we know, in this case there exist $\varkappa_l\in \Cn$
and a~submodule~$M_l$ of the module~$R^{(\varkappa_l)}$ such that
$$
y_l\cong R^{(\varkappa_l)}/M_l.
$$
Then we transform the element~$y_l$
into a~pair $\langle z_l^1,z_l^2\rangle$,
where $z_l^1=\varkappa_l$, $z_l^2$ is a~matrix of size
$\varkappa_l\times \varkappa_l$ over the ring~$R$,
and every column of~$z_l^2$ is a~vector from the generating set of
vectors of the module~$M_l$.
Naturally, in this case every column of the matrix~$M_l$
contains only a~finite number of nonzero elements.

If $y_l\in \Mor$, then $y_l$ is a~morphism from the module~$M_1$
into the module~$M_2$. Let
$$
M_1\cong R^{(\varkappa_1)}/N_1,\quad
M_2\cong R^{(\varkappa_2)}/N_2.
$$
Then for $m\in M_1$
$$
m=r_1e_{\alpha_1}+\dots+r_ke_{\alpha_k}+N_1,
$$
where $r_1,\dots,r_k\in R$ and $e_{\alpha_1},\dots,e_{\alpha_k}$
are elements of the basis of the module~$R^{(\varkappa_1)}$.
Let $y_l(m)=n\in M_2$, i.e., $n=s_1e_{\beta_1}+\dots+s_n e_{\beta_n}$,
where $s_1,\dots,s_n\in R$ and $e_{\beta_1},\dots,e_{\beta_n}$
are elements of the basis of the module~$R^{(\varkappa_2)}$.

We see that such a~morphism is completely defined by a~matrix
of size $\varkappa_1\times \varkappa_2$ such that
${y_l(N_1)\subset N_2}$.
Therefore we shall translate the morphism~$y_l$ to the elements
$z_l^1$, $z_l^2$, $z_l^3$, $z_l^4$, and $z_l^5$
where $z_l^1$~and~$z_l^2$ are the translations of the object
from which we are making this morphism, $z_l^3$~and~$z_l^4$ are
the translations of the object into which we are making our morphism,
and $z_l^5$ is the matrix of size $\varkappa_1\times \varkappa_2$
defined by the following formula: for every $\alpha\in \varkappa_1$
the $\alpha$-th column of the matrix~$z_l^5$ contains~$r_i$ in the row
with number $\beta_i\in \varkappa_2$
if $y_l(e_\alpha)=\sum r_ie_{\beta_i}$
(the column contains $0$ in all other rows).

Thus we obtain some new sequence $z_1,\dots,z_s$.
As it was done in the previous theorem, it is easy to show
by induction that
the sentence~$\tilde \varphi$
is true
on this sequence in the model
$\langle \Cn, R\rangle_{L_2}$
if and only if the sentence~$\varphi$
is true in the model $\textup{mod-}R$ on the sequence
$y_1,\dots,y_q$.
Thus, similarly to the previous subsection, we deduce that if
$\langle \Cn,R\rangle\equiv_{L_2}\langle \Cn,S\rangle$, then
$\textup{mod-}R\equiv \textup{mod-}S$.
\end{proof}

\subsection{An Analogue of the Morita Theorem and Its Corollaries}
\label{ss2.9}
The following theorem directly follows from Theorems
\ref{my2_7}~and~\ref{my2_8}.

\begin{theorem}\label{my2_9osn}\label{th7_2.9}
Let $R$~and~$S$ be rings.
Suppose that there exists
a~sentence~$\psi$ of the language $L_2(\langle \Cn,\ring\rangle)$
which is true in the ring~$R$ and is false in any ring similar
to~$R$ and not equivalent to it in the language
$L_2(\langle \Cn,\ring\rangle )$.
Then the categories $\textup{mod-}R$ and $\textup{mod-}S$ are
elementarily equivalent if and only if
there exists a~ring~$S'$ similar to the ring~$S$ and such that
the structures $\langle \Cn ,R\rangle $ and $\langle \Cn,S'\rangle$
are equivalent in the logic~$L_2$.
\end{theorem}

The most evident corollaries from Theorem~\ref{my2_9osn} are
the following two statements.

\begin{corollary}\label{co1_2.9}
For any skewfields $F_1$~and~$F_2$ the categories $\textup{mod-}F_1$ and
$\textup{mod-}F_2$ are elementarily equivalent if and only if
the structures $\langle \Cn,F_1\rangle$ and $\langle \Cn,F_2\rangle$
are equivalent in the second-order logic~$L_2$.
\end{corollary}

\begin{corollary}\label{co2_2.9}
For any commutative rings $R_1$~and~$R_2$
the categories $\textup{mod-}R_1$ and
$\textup{mod-}R_2$ are elementarily equivalent if and only if
the structures $\langle \Cn,R_1\rangle$ and
$\langle \Cn,R_2\rangle$ are equivalent in the second-order
logic~$L_2$.
\end{corollary}

\begin{proof}
In a~category $\textup{mod-}R$, where $R$ is a~commutative ring,
the formula $\Proobr(X)$ defines all progenerators~$X$, and
the formula
$$
\mathrm{Comm}(X)\prisv \Proobr(X) \logic\land
\forall f,g\in \Mor(X,X)\losp (f\circ g=g\circ f)
$$
defines all objects which are isomorphic to the ring~$R$
(see Theorem~\ref{kommut}).
\end{proof}

Also the corollaries from Theorem~\ref{my2_9osn}
for \emph{local rings} and \emph{integral domains} are
not difficult.

A~\emph{local ring} is a~ring in which the set of all
noninvertible elements is a~left ideal (see \cite[Lemma~1.2, p.~15]{9}).

\begin{proposition}\label{localring}
If $R$ is a~local ring, then every finitely generated projective
$R$-module is free.
\end{proposition}

\begin{proof}
We show that if a~ring~$R$ is local, then the set~$M$
of all noninvertible elements is also a~right ideal.
Indeed, suppose that some product $m\lambda$, where
$m\in M$ and $\lambda\in R$, is invertible. Then there exists $r\in R$
such that $m\cdot r=1$. It is clear that $r$ can not belong
to the left ideal~$M$. But $r$ can not be invertible either,
since in the opposite case the formula
$$
m=m(vv^{-1})=(mv)v^{-1}=v^{-1}
$$
shows that $m$ is also invertible.

This contradiction proves that $M$ is a~two-sided ideal.
It is clear that the factor ring $R/M$ is a~skewfield.

Note that a~square matrix over~$R$ is invertible if and only if
its reduction modulo the ideal~$M$ is invertible. To prove this
let us multiply this matrix from the left side by a~matrix
that represents an invertible matrix modulo~$M$,
then diagonalize this product with the help of elementary
transformations
of rows. Therefore the matrix has a~left inverse matrix;
similarly we can construct the right inverse matrix.

Suppose that a~module~$P$ is finitely generated and projective over~$R$.
Then we can find a~module~$Q$ such that $P\oplus Q\cong R^{(n)}$.
Choose bases in $P/MP$ and $Q/MQ$ (as in spaces over the skewfield $R/M$).
We shall lift up
every element of these bases
to~$P$ or to~$Q$, respectively.

This obtained set of elements is a~basis of the module $P\oplus Q$.
It is clear that therefore the module~$P$ is free.
\end{proof}

\begin{corollary}\label{co3_2.9}
For arbitrary local rings $R_1$~and~$R_2$ the categories
$\textup{mod-}R_1$ and $\textup{mod-}R_2$
are elementarily equivalent if and only if
the structures $\langle \Cn,R_1\rangle$ and $\langle \Cn,R_2\rangle$
are equivalent in the second-order logic~$L_2$.
\end{corollary}

\begin{proof}
In the category $\textup{mod-}R$, where $R$ is a~local ring, the formula
\begingroup
\setlength{\multlinegap}{0pt}
\begin{multline*}
\mathrm{Local}(X) \prisv \Proobr(X) \logic\land
\forall f,g,h\in \Mor(X,X)\losp
((\forall f'\in \Mor(X,X)\losp \neg(f\circ f'=f'\circ f=1_X))\\
\logic\land
(\forall g'\in \Mor(X,X)\losp \neg(g\circ g'=g'\circ g=1_X)) \logic\land
h=f\oplus g\Rightarrow
(\forall h'\in \Mor(X,X)\losp \neg(h\circ h'=h'\circ h=1_X)))
\end{multline*}
\endgroup%
holds only for modules which are isomorphic to the module~$R_R$.

Indeed, from Proposition~\ref{localring} it follows that the formula
$\Proobr(X)$ holds only for $X\circ R^{(n)}$. Let
$e_1,\dots,e_n$ be a~basis of the ring $R^{(n)}$, where $n\ge 1$.
Then consider $f,g,h\in \Mor(X,X)$ such that
$f(e_1)=e_1$, $f(e_i)=0$ for $i \ne 1$, $g(e_1)=0$, $g(e_i)=e_i$
for $i\ne 1$, and $h(e_i)=e_i$ for every $i=1,\dots, n$.

Then for morphisms $f$, $g$, and $h$
\begin{multline*}
(\forall f'\in \Mor(X,X)\losp \neg (f\circ f'=f'\circ f=1_X))
\logic\land
(\forall g'\in \Mor(X,X)\losp \neg (g\circ g'=g'\circ g=1_X))\\
\logic\land
(h=f\oplus g) \logic\land
\exists h'\in \Mor(X,X)\losp (h\circ h'=h'\circ h=1_X),
\end{multline*}
where $h'=h$.

Therefore in the module~$X$ the formula $\mathrm{Local}(X)$ does not hold.
\end{proof}

A~ring~$R$ is called an \emph{integral domain} if
it does not contain any zero divisors and each of its
ideals is \emph{principal} (is generated by an element).

\begin{proposition}[{see \cite[Chap.~XV, Sec.~2]{10}}]\label{princideal}
Let $P$ be a~progenerator over an integral domain. Then
the module~$P$ is free.
\end{proposition}

\begin{proof}
Since $P$ is a~progenerator, it is a~submodule of the module~$R^{(n)}$.
Let the module $R^{(n)}$ have a~basis $e_1,\dots ,e_n$, and let $P_r$
be the intersection of the module~$P$ with the module
$\langle e_1,\dots,e_r\rangle$.
Then $P_1=P\cap \langle e_1\rangle$ is a~submodule in $\langle e_1\rangle$
and hence has the form $\langle r_1e_1\rangle$ for some $r_1\in R$.
Thus the module $P_1$ is either nonzero or free of rank~$1$.
Suppose by induction that the module~$P_r$ is free of rank~${\le}\,r$.
Let $M$ be the set of all elements $m\in R$ such that
there exists $x\in P$ which can be written in the form
$$
x=b_1e_1+\dots+b_re_r+me_{r+1},
$$
where $b_i\in R$.

It is clear that $M$ is an ideal in~$R$ and, therefore,
is a~principal ideal, generated by some
$r_{r+1}\in R$. If $r_{r+1}=0$, then $P_{r+1}=P_r$ and
the induction step is proved.
If $r_{r+1}\ne 0$, then let $w\in P_{r+1}$ be such that
its $e_{r+1}$-th coefficient is equal to $r_{r+1}$. If $x\in P_{r+1}$,
then its $e_{r+1}$-th coefficient can be divided by $r_{r+1}$ and,
therefore, there exists such $c\in R$ that
$x-cw\in P_r$. Consequently,
$$
P=P_r+\langle w\rangle.
$$
On the other hand,
$P_r\cap \langle w\rangle=0$, and therefore
this sum is direct.
\end{proof}

\begin{corollary}\label{co4_2.9}
For arbitrary integral domains $R_1$~and~$R_2$ the categories
$\textup{mod-}R_1$ and $\textup{mod-}R_2$ are elementarily
equivalent if and only if the structures
$\langle \Cn,R_1\rangle$ and $\langle \Cn,R_2\rangle$
are equivalent in the logic~$L_2$.
\end{corollary}

\begin{proof}
In a~category $\textup{mod-}R$, where $R$ is an integral domain,
the formula
$$
\mathrm{Principal}(X)\prisv
\Proobr(X) \logic\land
\forall f\in \Mor(X,X)\losp \forall g\in \Mor(X,X)\losp
(f\circ g\ne 0 \logic\land g\circ f\ne 0)
$$
holds only for modules which are isomorphic to the module~$R_R$.
This follows easily from Proposition~\ref{princideal}.\looseness1
\end{proof}

A module $M$ over a~ring $R$ is called \emph{Artinian}
if the following equivalent conditions are fulfilled:
\begin{enumerate}
\item
every nonempty set of submodules of the module~$M$,
ordered by inclusion, contains a~minimal element;
\item
every decreasing sequence of submodules of the module~$M$
is stationary.
\end{enumerate}

A~ring~$R$ is called \emph{Artinian} if the module~$R_R$
is Artinian.

A~module~$M$ is called \emph{decomposable} if there exist such modules
$M_1$~and~$M_2$ that $M=M_1\oplus M_2$. In the opposite case
a~module~$M$ is called \emph{indecomposable}.

In \cite[p.~139]{11} the following theorem is proved.

\begin{theorem}\label{artin}
Let $M$ be a~finitely generated module over an Artinian ring~$R$.
\begin{enumerate}
\renewcommand{\theenumi}{\alph{enumi}}
\item
The module~$M$ can be represented as a~direct sum
of a~finite family $(M_i)_{1\le i\le m}$
of indecomposable nonzero submodules.
\item
If the module~$M$ is a~direct sum of another
family $(M_j')_{1\le j\le n}$ of indecomposable nonzero submodules,
then $m=n$ and there exist a~substitution~$\pi$ of
the set $\{ 1,\dots,n\}$ and an automorphism~$\alpha$
of the set~$M$ such that
$$
\alpha (M_j')=M_{\pi(j)},\quad 1\le j\le n.
$$
\end{enumerate}
\end{theorem}

Now introduce the following sentences
of the second-order language of the structure $\langle \Cn,\ring\rangle$.

1. For a~subset~$M$ of the ring the formula
$$
\Mod(M) \prisv \forall r\in \ring\losp \forall m\in M\losp
\exists n\in M\losp (rm=n) \logic\land
\forall l,m\in M\losp \exists n\in M\losp (n=l+m)
$$
means that the set~$M$ is a~module over the ring $\ring$.

2. For sets $M$~and~$N$ the formula
\begin{align*}
& (M\cong N) \prisv \Mod(N) \logic\land \Mod(M)
\logic\land
\exists F(v_1,v_2)\losp
(\Dom (F)=M \logic\land \Rng(F)=N \logic\land \Bij(F)\\
& \quad \logic\land
\forall r_1,r_2\in \ring\losp \forall m_1,m_2\in M\losp
\forall n_1,n_2\in N\\
& \quad (\langle m_1,n_1\rangle \in F \logic\land
\langle m_2,n_2\rangle \in F \Rightarrow
\langle r_1m_1+r_2m_2,r_1n_1+r_2n_2\rangle\in F))
\end{align*}
means that the sets $M$~and~$N$ are $\ring$-modules and
that they are isomorphic.

3. For sets $L,M,N\subset \ring$ the formula
\begin{multline*}
(N=M\oplus L) \prisv \Mod(M) \logic\land \Mod(L) \logic\land
\Mod(N)\\
\logic\land
\forall n\in N\losp \exists m\in M\losp \exists l\in L\losp
(n=m+l) \logic\land
\forall m\in M\losp \forall l\in L\losp (m=l\Rightarrow m=0)
\end{multline*}
means that the module~$N$ is a~direct sum of the modules $M$~and~$L$.

4. For a~set $M\subset \ring$ the formula
$$
\Undir(M) \prisv \Mod(M) \logic\land
\forall L(c) ,N(c)\losp \neg (M=L\oplus N)
$$
means that the module~$M$ is indecomposable.

5. For a~set $M\subset \ring$ the formula
\begin{multline*}
\mathrm{Dir}_N (M) \prisv \Mod(M) \logic\land
\exists M_1(c)\ldots\exists M_N(c)\losp
(\Mod(M_1) \logic\land \ldots \logic\land \Mod(M_N))\\
\logic\land \bigwedge_{i\ne j} \neg (M_i\cong M_j)
\logic\land M=M_1\oplus \ldots\oplus M_N
\logic\land
(\Undir (M_1)\logic\land\ldots \logic\land \Undir(M_N))
\end{multline*}
means that the module~$M$ is a~direct sum of indecomposable
modules $M_1, \dots, M_N$ which are not isomorphic to each other.

Suppose that we have some Artinian ring~$R$. Then the module~$R_R$
is Artinian, and therefore it is a~direct sum
of $n$~indecomposable modules. Let it be modules
$$
M_1^{1},\dots,M_1^{i_1},M_2^1,\dots,M_2^{i_2},\dots,M_k^{1},\dots,M_k^{i_k},
$$
and for $k\ne l$
$$
M_k^i\not\cong M_l^j,
$$
but for every~$k$
$$
M_k^i\cong M_k^j.
$$

Consider the module
$$
M \prisv M_1^1\oplus \dots \oplus M_k^1.
$$
Since the module~$M$ is a~direct summand of the module~$R_R$,
it is projective and finitely generated. Since the module~$R_R$
is a~direct summand of the module $M^{(\max(i_1,\dots,i_k))}$,
we see that $M$ is a~generator. Therefore the module~$M$ is a~progenerator
and the ring $\Endom_RM$ is similar to the ring~$R$.

Thus for some $N\in \omega$ the formula
$$
\psi(P) \prisv \Proobr(P)\land \Undir_N(P)
$$
defines a~unique, up to an isomorphism, progenerator
$$
M \prisv M_1^1\oplus \dots\oplus M_k^1.
$$

Consequently, we have proved the following corollary.

\begin{corollary}\label{co5_2.9}
For any Artinian rings $R_1$~and~$R_2$
the categories $\textup{mod-}R_1$ and
$\textup{mod-}R_2$ are elementarily equivalent if and only if there exist
rings $S_1$~and~$S_2$ such that the ring~$S_1$ is similar to
the ring~$R_1$, the ring~$S_2$ is similar to the ring~$R_2$,
and the structures $\langle \Cn, S_1\rangle$
and $\langle \Cn,S_2\rangle$ are equivalent in the
second-order logic~$L_2$.
\end{corollary}

\section{Elementary Equivalence of Endomorphism Rings of Modules
of Infinite Ranks}\label{s3}
\subsection[Endomorphism Rings of Modules and Categories
$C_{M(V)}$]{Endomorphism Rings of Modules and Categories
$\boldsymbol{C_{M(V)}}$}\label{ss3.1}
Suppose that we have some associative ring~$R$ with~$1$, an infinite
cardinal number~$\varkappa$, and a~free module $V=V_\varkappa^R$ of
rank $\varkappa$ over~$R$.

In this section, we assume that every ideal
of the ring~$R$ is generated by at most $\varkappa$ elements.
This is always so if $\varkappa\ge |R|$, or if $R$ is an integral domain,
or if the ring~$R$ is semisimple.

In the ring $\Endom_R(V)$ we want to interpret the category~$C_V$,
consisting of the modules~$V$,
all quotient modules of the module~$V$,
and all homomorphisms between them, i.e., to give an algorithm,
transforming every formula $\varphi$ of the first order language
of the category theory to a~formula~$\tilde \varphi$
of the first order language of the ring theory in such a~way that
the formula~$\varphi$ holds in $C_{M(V)}$ if and only if
$\tilde \varphi$ holds in $\Endom_R(V)$.

At the beginning we shall give an informal description of this translation.

1. To every object~$X$ of the category $C_{M(V)}$ we associate
an element~$\tilde X$ of the ring $\Endom_RV$ in the following way:
if $X\in C_{M(V)}$, then $X=V/X'$ for some~$X'$
which is a~submodule of the module~$V$.
Every submodule of the module~$V$ can be defined by the generating
vectors, and the cardinality of the set of generating vectors
is not greater than~$\varkappa$. These vectors can be written
as columns of a~matrix of size
$\varkappa\times \varkappa$ (if this cardinality is less than~$\varkappa$,
then we can extend this matrix by zero columns),
i.e., as an element of the ring $\Endom_RV$.
Conversely, if $\tilde X \in \Endom_RV$, then we can consider the module
generated by the columns of the matrix~$\tilde X$,
and then the factormodule $X\prisv V/\tilde X$.

2. To every morphism~$f$ of the category $C_{M(V)}$ we associate
a~triplet $\langle X_f,Y_f,\tilde f\rangle$
of elements of the ring $\Endom_RV$ such that if $f\in \Mor(X,Y)$,
then $X_f=\tilde X$, $Y_f=\tilde Y$, and $\tilde f$
is a~matrix, establishing a~homomorphism $\tilde f\in \Mor(V,V)$ such that
$$
\tilde f\circ p_Y=p_X\circ f,
$$
where $p_X$ and $p_Y$ are standard epimorphisms from the module~$V$
onto the modules $X$~and~$Y$, respectively.

This condition shows that the matrix~$\tilde f$ has to translate
vectors of the module~$X'$ into vectors of the module~$Y'$, i.e.,
the matrix $\tilde f\tilde X$ has to generate
a~submodule of the module generated by the matrix~$\tilde Y$.
This means that there exists $A\in \Endom_RV$ such that
$$
\tilde f\tilde X=\tilde YA.
$$

Two endomorphisms of the module~$V$ define the same morphism
from the module~$X$ into the module~$Y$ if their difference
defines a~zero morphism from the module~$X$ into the module~$Y$,
i.e., the image of this morphism belongs to the module~$Y'$.

Therefore triplets $\langle X_f, Y_f, \tilde f_1\rangle$
and $\langle X_f,Y_f,\tilde f_2\rangle$ are considered as equal
if
$$
\exists A\losp (f_1-f_2=Y_fA).
$$

Now we shall give a~formal description.

1. A~subformula $\forall X\in \Obj$ is translated to
the subformula $\forall \tilde X$
(similarly for a~subformula $\exists X\in \Obj$).

2. A~subformula $\forall f\in \Mor$ is translated to the subformula
$$
\forall X_f\losp \forall Y_f\losp \forall \tilde f\losp
(\exists A\losp (\tilde f\circ X_f=Y_f\circ A) \Rightarrow \ldots)
$$
(similarly for a~subformula $\exists f\in \Mor$).

3. A~subformula $f\in \Mor(X,Y)$ is translated to the subformula
$X_f=\tilde X \logic\land Y_f=\tilde Y$.

4. A~subformula $h=f\circ g$ is translated to the subformula
$\tilde h=\tilde f\circ \tilde g$.

5. A~subformula $f=1_X$ is translated to the subformula
$X_f=Y_f=\tilde X \logic\land \tilde f=1$.

The algorithm is constructed. Similarly to the previous sections,
we can show that the sentence~$\varphi$ holds in the category
$C_{M(V)}$ if and only if the sentence~$\tilde \varphi$
holds in the ring $\Endom_RV$.

Now note that we shall consider not simply the structure
$C_{M(V)}$ with the language of the category theory,
but the structure $C_{M(V)}$ with the selected module~$V$,
i.e., in formulas we can use the subformula $X=V$ for $X\in \Obj$.
This subformula will be translated to the subformula
$\tilde X=0$.

Therefore, if rings $\Endom_RV$ and $\Endom_SW$ are elementarily equivalent,
then the categories $C_{M(V)}$ and $C_{M(W)}$ are also
elementarily equivalent.

We now prove the inverse implication.

To do this we need
to interpret the ring $\Endom_RV$
inside the category $C_{M(V)}$ with the selected
object~$V$.

Indeed, in the category $C_{M(V)}$ we shall fix some
$V^2\in \Obj$ such that
$V^2\cong V\oplus V$ (for example, $V\cong V\oplus V$) and morphisms
$i_1,i_2\in \Mor(V,V^2)$ and $p_1,p_2\in \Mor(V^2,V)$ such that
$$
p_1\circ i_1=p_2\circ i_2=1_V \logic\land
p_1\circ i_2=p_2\circ i_1 =0 \logic\land
\forall i\in \Mor(V,V^2)\losp
(i\ne 0_{V,V^2}\Rightarrow p_1 \circ i\ne 0_{V} \logic\lor
p_2\circ i\ne 0_V).
$$

It is clear that in this case the morphisms $i_1$~and~$i_2$
are embeddings of the module~$V$ into the module $V\oplus V$,
their images do not intersect, and their sum is $V\oplus V$.

Now translate the subformulas $\forall f$ and $\exists f$ to
the subformulas $\forall f\in \Mor(V,V)$ and $\exists f\in \Mor(V,V)$;
and the subformulas $h=f\cdot g$ and $h=f+g$ to the subformulas
$h=f\circ g$ and $h=f\oplus g$ (see Sec.~2.5).

Therefore we have that $C_{M(V_1)}\equiv C_{M(V_2)}$ implies
$\Endom_{R_1}(V_1)\equiv \Endom_{R_2}(V_2)$.

Consequently,
the question of elementary equivalence of endomorphism rings
$\Endom_{R_1}(V_1)$ and $\Endom_{R_2}(V_2)$
is equivalent to the question of elementary equivalence of the categories
$C_{M(V_1)}$ and $C_{M(V_2)}$ with selected objects $V_1$~and~$V_2$,
respectively.
{\sloppy

}

\subsection[Elementary Equivalence in Categories
$C_{M(V)}$]{Elementary Equivalence in Categories
$\boldsymbol{C_{M(V)}}$}\label{ss3.2}
Note that our new situation is very close
to the situation of Sec.~2.
We have the category $C_{M(V)}$, which is a~subcategory in
$\textup{mod-}R$ and is closed under taking
quotient modules and direct products of cardinality
at most $\varkappa$.
This category resembles the category
$\textup{mod-}R$, but it is small and bounded by the given
cardinal number~$\varkappa$.
Moreover, in this category the module~$V$ is selected.

We generalize all possible results from Sec.~2 to this case.

The formula $\Simp(M)$ also defines in the category $C_{M(V)}$
simple modules because this category
is closed under taking factormodules. The formula $\tSum^\omega(X,M)$
also defines the module $X\simeq M^{(\omega)}$
because the cardinal number~$\varkappa$ by the condition
is at most~$\omega$.
It is clear that the formula $\tSum^{\Fin}(X,M)$ holds for finite
direct sums of the module~$M$,
and the formula $\tSum(X,M)$ holds for all direct sums of the modules~$M$
which belong to the category $C_{M(V)}$. Similarly
we can generalize for the case of the category~$C_V$ all formulas
from Sec.~2.2, and even the formula $\Proobr(P)$
which defines in this category all progenerators.

After selecting some progenerator~$P$ completely similarly to Sec.~2.3
we can construct an analogue of the ring $\Endom_RP$ because in Sec.~2.3
we used only closedness of the category $\textup{mod-}R$ under finite
direct sums.

Since all results of Sec.~2.4 also can be easily generalized
to our case, we have the following theorem.

\begin{theorem}\label{my3_2fin}
Categories $C_{M(V_1^R)}$ and $C_{M(V_2^S)}$, where $R$ is a~finite ring,
are elementarily equivalent if and only if $R\cong \Endom_S P$
for some progenerator~$P$ of the category $C_{M(V_2^S)}$.
\end{theorem}

It is also clear that following Secs.\ 2.5~and~2.6 we can find
a~formula $\varphi(f)$ which holds for some independent set
of mappings $f\colon V\to P$ of cardinality~$\varkappa$ such that
for every~$f$ there exists
$g\colon P\to V$ such that $f\circ g=1_P$ and $g\circ f$ is a~projector
from~$V$ into~$V$.

Indeed, for these objects we get similar results.
For this purpose we consider, together with the
full language $L_2(\langle \Cn, \ring\rangle)$,
its part which can be described as follows.

As we said before (see Sec.~1), the \emph{theory}
of a~given model~$\mathcal U$ in a~language~$\mathcal L$
is the set of all sentences of the language~$\mathcal L$
which are true in the model~$\mathcal U$.
It is clear that two models $\mathcal U$~and~$\mathcal V$
in the same language~$\mathcal L$
are equivalent in the language~$\mathcal L$ if and only if
their theories in this language coincide.

The theory of the structure $\langle \Cn,R\rangle$ in the language~$L_2$
is denoted by $\Th_2(\langle \Cn,R\rangle)$.

We can also consider the structure $\langle \varkappa,R\rangle$,
consisting of a~set of cardinality~$\varkappa$ and the ring~$R$
with ring operations $+$~and~$\circ$.

By $\Th_2^\varkappa(\langle \varkappa,R\rangle)$
we shall denote the part of the theory
$\Th_2(\langle \varkappa,R\rangle)$
bounded by the cardinal number~$\varkappa$,
i.e., the sentences $\varphi\in \Th_2(\langle \varkappa,R\rangle)$
such that the quantifiers $\forall $~and~$\exists $ appear only with
the predicate symbols
$$
P(c_1,\dots,c_k;v_1,\dots,v_n),
$$
where the set
$$
\{ \langle \alpha_1,\dots,\alpha_k,r_1,\dots,r_n\rangle \mid
\alpha_1,\dots,\alpha_k\in \varkappa \logic\land
r_1,\dots,r_n\in R \logic\land
P (\alpha_1,\dots,\alpha_k;r_1,\dots,r_n)\}
$$
is of cardinality at most~$\varkappa$.

Then we can write the following analogue of
Theorem~\ref{my2_9osn} from Sec.~2.

\begin{theorem}\label{my3_7}
Let $V_1$~and~$V_2$ be free modules of infinite ranks
$\varkappa_1$~and~$\varkappa_2$ over rings $R_1$~and~$R_2$,
respectively. Suppose that
there exists a~sentence
$\psi\in \Th_2^{\varkappa_1}(\langle \varkappa_1,R_1\rangle)$
such that $\psi\notin \Th_2^{\varkappa_1}(\langle \varkappa_1,R'\rangle)$
for every ring~$R'$
such that $R'$ is similar to~$R_1$ and
$\Th_2^{\varkappa_1} (\langle \varkappa_1,R_1\rangle )\ne
\Th_2^{\varkappa_1} (\langle \varkappa_1, R'\rangle )$.
Then if the categories $C_{V_1}$ and $C_{V_1}$ are elementarily
equivalent, then
there exists a~ring~$S$ similar to the ring~$R_2$ and such that
the theories $\Th_2^{\varkappa_1}(\langle \varkappa_1 ,R_1\rangle) $
and $\Th_2^{\varkappa_2}(\langle \varkappa_2,S\rangle)$ coincide.
\end{theorem}

\begin{proof}
The proof of this theorem resembles the proof of
Theorem~\ref{my2_7} from Sec.~2, but we shall give it in
detail to show differences.

At the beginning we assume that we fix some progenerator~$P$
in the category $C_{M(V)}$, where
$V=V_T^{\varkappa}$, $\varkappa$~is an infinite cardinal number, and
$T$~is a~ring (it is clear that all progenerators of the
category $\textup{mod-}T$ are contained in the category $C_{M(V)}$).
Then we have formulas defining a~simple module~$M$
corresponding to the module~$P$,
modules~$M^{(\alpha)}$ for all $\alpha\in \Cn\cap \varkappa+1$,
modules~$M^{(n)}$ for all $\alpha\in \omega$,
modules~$M^{(\alpha)}$ for infinite $\alpha\in \Cn\cap \varkappa+1$,
almost free modules~$V^\alpha$ of ranks
$\alpha \in \Cn\cap \varkappa+1$, $\alpha\in \omega$,
$\alpha\in \Cn\cap \varkappa+1\setminus \omega$,
and also the selected free module~$V$, which is
almost free over the module~$P$.

For the module $M^{(\varkappa)}$ ($V$) we shall define
(see Sec.~2.5) its generating set of projectors
$\Gen_{\bar g^*}(M^{(\varkappa)},M)$ (or $\Gen_{\bar g^*}(V,P)$).

Further (see Sec.~2.3), for every $f,g\in \Mor(P,P)$ we assume that
their sum $f\oplus g\in \Mor(P,P)$
and their product $f\otimes g\in \Mor(P,P)$ are known.

Consider an arbitrary sentence~$\varphi$
of the language $L_2(\langle \varkappa,\ring\rangle)$.
This sentence can contain the following subformulas.
\begin{enumerate}
\renewcommand{\labelenumi}{\theenumi.}
\item
$\forall (\exists ) r\in \ring$.
\item
$\forall (\exists )\alpha\in \varkappa$.
\item
$r_1=r_2+r_3$.
\item
$r_1=r_2\cdot r_2$.
\item
$r_1=r_2$.
\item
$\alpha_1=\alpha_2$.
\item
$\forall (\exists) P(c_1,\dots,c_k;v_1,\dots,v_n)$.
\item
$P(\alpha_1,\dots,\alpha_k;r_1,\dots,r_n)$.
\end{enumerate}

Translate this sentence to the sentence~$\tilde \varphi_P$
(depending of the fixed module~$P$) of the first order language
of the category theory by the following algorithm.
\begin{enumerate}
\renewcommand{\labelenumi}{\theenumi.}
\item
A~subformula $\forall (\exists) r\in \ring$ is translated
to the subformula $\forall (\exists )f_r\in \Mor(P,P)$,
i.e., every element of the ring $\ring$
corresponds to an element of the ring $\Endom_TP$.
\item
A~subformula $\forall (\exists )\alpha\in \varkappa$
is translated to the subformula
$\forall (\exists) F^\alpha \in \Gen_{\bar g^*}(M^{(\varkappa)},M)$.
\item
A~subformula $r_1=r_2+r_3$ is translated to the subformula
$f_{r_1}=f_{r_2}\oplus f_{r_3}$.
\item
A~subformula $r_1=r_2\cdot r_3$ is translated to the subformula
$f_{r_1} =f_{r_2}\otimes f_{r_3}$.
\item
A~subformula $r_1=r_2$ is translated to the subformula
$f_{r_1}=f_{r_2}$.
\item
A~subformula $\alpha_1=\alpha_2$ is translated to the subformula
$f^{\alpha_1}=f_{\alpha_2}$.
\item
A~subformula $\forall (\exists) P(c_1,\dots,c_k;v_1,\dots,v_n)$
is translated to the subformula
\begin{multline*}
\forall (\exists) f_P^{c_1}\in \Sets(M^{(\varkappa)},M^{(\varkappa)})
\ldots
\forall (\exists) f_P^{c_k}\in \Sets(M^{(\varkappa)},M^{(\varkappa)})\\
\forall (\exists) f_P^{v_1}\in \Ring (V)\dots
\forall (\exists) f_P^{v_n}\in \Ring(V).
\end{multline*}
\item
A~subformula
$P(\alpha_1,\dots,\alpha_k;r_1,\dots,r_n)$
is translated to the subformula
\begin{multline*}
\exists f\in \Gen(M^{(\varkappa)},M)\losp
(f^{\alpha_1}\circ f_P^{c_1}\circ \bar f=1
\logic\land \ldots \logic\land
f^{\alpha_k}\circ f_P^{c_k}\circ \bar f=1\\
\logic\land
f'\circ f_P^{v_1}\circ \bar f'= f_{r_1}
\logic\land \ldots \logic\land
f'\circ f_P^{v_n}\circ \bar f'=f_{r_n}).
\end{multline*}
\end{enumerate}

As it was done in Theorem~\ref{my2_7} of Sec.~2, we can show that
the sentence~$\varphi$ holds in the theory
$\langle \varkappa,\Endom_TP\rangle$
if and only if the sentence~$\tilde \varphi_P$
holds in the model $C_{M(V_T^\varkappa)}$, whence,
similarly to Theorem~\ref{my2_7} from Sec.~2, we prove the theorem.
\end{proof}

\begin{theorem}\label{my3_8}
If $\varkappa_1$~and~$\varkappa_2$ are infinite cardinal numbers,
$V_1$~and~$V_2$ are free modules of ranks $\varkappa_1$~and~$\varkappa_2$
over the rings $R_1$~and~$R_2$, respectively, and
the theories $\Th_2^{\varkappa_1}(\langle \varkappa_1,R_1\rangle)$
and $\Th_2^{\varkappa_2}(\langle \varkappa_2, R_2\rangle)$
coincide, then the categories $C_{M(V_1)}$
and $C_{M(V_2)}$ are elementarily equivalent.
\end{theorem}

\begin{proof}
The proof of this theorem differs from the proof of
Theorem~\ref{my2_8} from Sec.~2 only in the moment that
the module~$V$ has to be the selected object of the category
$C_{M(V)}$. But since by the theorem condition we consider only
free modules
(only at this point it is important that the modules a~free,
but not almost free), we have that the selected object of the category
will be the zero matrix.
\end{proof}

A~direct corollary from Theorems \ref{my3_7}~and~\ref{my3_8}
is Theorem~\ref{3_9osn}.

\begin{theorem}\label{3_9osn}
Let $V_1$~and~$V_2$ be free modules of infinite ranks
$\varkappa_1$~and~$\varkappa_2$ over rings $R_1$~and~$R_2$,
respectively. Suppose that
there exists a~sentence
$\psi\in \Th_2^{\varkappa_1}(\langle \varkappa_1,R_1\rangle)$
such that
$\psi\notin \Th_2^{\varkappa_1}(\langle \varkappa_1,R'\rangle)$
for every ring~$R'$
such that $R_1$ is similar to~$R'$ and
$\Th_2^{\varkappa_1}(\langle \varkappa_1,R_1\rangle)\ne
\Th_2^{\varkappa_1} (\langle \varkappa_1, R'\rangle)$.
Then the categories $C_{M(V_1)}$ and $C_{M(V_1)}$
are elementarily equivalent if and only if
there exists a~ring~$S$ similar to the ring~$R_2$ and such that
the theories $\Th_2^{\varkappa_1}(\langle \varkappa_1 ,R_1\rangle) $
and $\Th_2^{\varkappa_2}(\langle \varkappa_2,S\rangle)$
coincide.
\end{theorem}

\subsection{The Main Theorem}\label{ss3.3}
The previous results imply the following theorem.

\begin{theorem}\label{3osn}\label{th5_3.3}
Let $V_1$~and~$V_2$ be free modules of infinite ranks
$\varkappa_1$~and~$\varkappa_2$ over rings $R_1$~and~$R_2$,
respectively.
Suppose that there exists a~sentence
$\psi\in \Th_2^{\varkappa_1}(\langle \varkappa_1,R_1\rangle)$
such that
$\psi\notin \Th_2^{\varkappa_1}(\langle \varkappa_1,R'\rangle)$
for every ring~$R'$
such that $R_1$ is similar to~$R'$ and
$\Th_2^{\varkappa_1}(\langle \varkappa_1,R_1\rangle)\ne
\Th_2^{\varkappa_1} (\varkappa_1, R'\rangle)$.
Then the categories $C_{M(V_1)}$ and $C_{M(V_1)}$
are elementarily equivalent if and only if
there exists a~ring~$S$ similar to the ring~$R_2$ and such that
the theories $\Th_2^{\varkappa_1}(\langle \varkappa_1 ,R_1\rangle) $
and $\Th_2^{\varkappa_2}(\langle \varkappa_2,S\rangle)$ coincide.
\end{theorem}

\begin{corollary}\label{co1_3.3}
Let $V_1$~and~$V_2$ be two spaces of infinite ranks
$\varkappa_1$~and~$\varkappa_2$ over arbitrary skewfields
\textup{(}integral domains\textup{)} $F_1$~and~$F_2$.
Then the rings
$\Endom_{F_1} V_1$ and $\Endom_{F_2}V_2$
are elementarily equivalent if and only if the theories
$\Th_2^{\varkappa_1}(\langle \varkappa_1,F_1\rangle)$ and
$\Th_2^{\varkappa_2}(\langle \varkappa_2,F_2\rangle)$ coincide.
\end{corollary}

\begin{corollary}\label{co2_3.3}
Suppose that $\varkappa_1$~and~$\varkappa_2$ are infinite
cardinal numbers, $R_1$~and~$R_2$ are commutative \textup{(}local\textup{)}
rings, and every maximal ideal of the ring~$R_1$
is generated by at most~$\varkappa_1$ elements of the ring.
Then for free modules $V_1$~and~$V_2$ of ranks
$\varkappa_1$~and~$\varkappa_2$ over the rings $R_1$~and~$R_2$, respectively,
the rings $\Endom_{R_1}V_1$ and $\Endom_{R_2}V_2$ are elementarily
equivalent if and only if the theories
$\Th_2^{\varkappa_1}(\langle \varkappa_1,R_1\rangle)$ and
$\Th_2^{\varkappa_2}(\langle \varkappa_2,R_2\rangle)$ coincide.
\end{corollary}

\begin{corollary}\label{co3_3.3}
Suppose that $\varkappa_1$~and~$\varkappa_2$ are infinite cardinal numbers,
$R_1$~and~$R_2$ are Artinian rings, and every maximal
ideal of the ring~$R_1$ is generated by at most~$\varkappa_1$
elements of the ring. Then for free modules $V_1$~and~$V_2$ of ranks
$\varkappa_1$~and~$\varkappa_2$ over the rings $R_1$~and~$R_2$,
respectively, the rings $\Endom_{R_1}V_1$ and $\Endom_{R_2}V_2$
are elementarily equivalent if and only if
there exist rings $S_1$~and~$S_2$ similar to the rings $R_1$~and~$R_2$,
respectively, such that the theories
$\Th_2^{\varkappa_1}(\langle \varkappa_1,S_1\rangle)$ and
$\Th_2^{\varkappa_2}(\langle \varkappa_2,S_2\rangle)$ coincide.
\end{corollary}

\begin{corollary}\label{co4_3.3}
For free modules $V_1$~and~$V_2$ of infinite ranks
$\varkappa_1$~and~$\varkappa_2$ over semisimple rings
$R_1$~and~$R_2$, respectively, the rings
$\Endom_{R_1}(V_1)$ and $\Endom_{R_2}(V_2)$
are elementarily equivalent if and only if
there exist rings $S_1$~and~$S_2$ similar to the rings
$R_1$~and~$R_2$, respectively,
such that the theories
$\Th_2^{\varkappa_1}(\langle \varkappa_1,S_1\rangle)$ and
$\Th_2^{\varkappa_2}(\langle \varkappa_2,S_2\rangle)$ coincide.
\end{corollary}

\section[The Projective Space of the Module~$V$]{The Projective
Space of the Module~$\boldsymbol{V}$}\label{s4}
\subsection{The Language of the Projective Space and
Basic Notions, Definable in This Language}\label{ss4.1}
Suppose that we have some free module~$V$ of infinite
rank~$\varkappa$ over a~ring~$R$. The \emph{projective space}
$P(V)$ of the module~$V$ is an algebraic structure consisting
of all submodules of the module~$V$ with the relation $\subset$
(we write $M\subset N$
if the module~$M$ is a~submodule of the module~$N$).

In this section, we assume that every submodule of
the module~$V$ can be generated by at most $\varkappa$
elements of the module~$V$ (this is true if $\varkappa\ge |R|$,
or if the ring~$R$ is semisimple,
or if the ring~$R$ is an integral domain).

Let $M_1,M_2,M_3\in P(V)$.
We shall write that $ M_1=V$ if $\forall M\losp (M\subset M_1)$.
We shall also write that $M_1=\varnothing$
if $\forall M\losp (M_1\subset M)$.
The formula $M_1=M_2\cap M_3$ will denote the formula
$$
M_1\subset M_2 \logic\land
M_1\subset M_3 \logic\land
\forall M_4\losp
(M_4\subset M_2 \logic\land M_4\subset M_3 \Rightarrow M_4\subset M_1),
$$
the formula $M_1=M_2+M_3$ will denote the formula
$$
M_2\subset M_1 \logic\land
M_3\subset M_1 \logic\land
\forall M_4\losp
(M_2\subset M_4 \logic\land M_3\subset M_4\Rightarrow M_1\subset M_4),
$$
and the formula $M_1=M_2\oplus M_3$ will denote the formula
$$
M_1=M_2+M_3 \logic\land M_2\cap M_3=\varnothing.
$$
It is clear that if $M_1=M_2\cap M_3$, then the module~$M_1$ is
the intersection of the modules $M_2$~and~$M_3$,
if $M_1=M_2+M_3$, then it is the sum of the
modules $M_2$~and~$M_3$, and if $M_1=M_2\oplus M_3$, then it is
the direct sum of the modules $M_2$~and~$M_3$.

Consider now for given modules $P_1$~and~$P_2$ the formula
\begin{multline*}
P_1\cap P_2=\varnothing\\
\logic\land
\exists P\losp
(P\subset P_1\oplus P_2 \logic\land P\ne\varnothing \logic\land
P\cap P_1=\varnothing
\logic\land
P\cap P_2=\varnothing \logic\land P\oplus P_1=P_1\oplus P_2
\logic\land P\oplus P_2=P_1\oplus P_2).
\end{multline*}
Let the modules $P_1$~and~$P_2$ not intersect and let
there be a~module~$P$ satisfying all conditions from the formula.
Since $P\subset P_1\oplus P_2$, it follows that every $ x\in P$ has the form
$x=y+z$, where $y\in P_1$, $z\in P_2$, and the elements
$z$~and~$y$ are uniquely defined by the vector~$x$.
Consider the correspondence $F\subset P_1\times P_2$
which is defined by the formula
$$
\forall y\in P_1\losp \forall z\in P_2\losp
\langle y,z\rangle \in F \Leftrightarrow \exists x\in P\losp (x=y+z).
$$
We show that $F$ is an isomorphism between the modules $P_1$~and~$P_2$.

1. If $y_1,y_2\in P_1$, $z\in P_2$,
$\langle y_1,z\rangle \in F$, and
$\langle y_2,z\rangle \in F$, then
$\exists x_1,x_2\in P\losp (x_1=y_1+z \logic\land x_2=y_2+z)$, i.e.,
$x\prisv x_1-x_2=y_1-y_2\in P$. Since in this case
$y_1-y_2\in P$, it follows that
$y_1-y_2\in P\cap P_1\Rightarrow y_1-y_2=0\Rightarrow y_1=y_2$.

2. Similarly, from $y\in P_1$, $z_1,z_2\in P_2$,
$\langle y,z_1\rangle \in F$, and
$\langle y,z_2\rangle \in F$ it follows that $z_1=z_2$.

3. Consider an arbitrary vector $y\in P_1$. Since
$y\in P_1\oplus P_2$, it follows that $y\in P\oplus P_2$, i.e.,
$\exists x\in P\losp \exists z\in P_2\allowbreak\losp (y=x+z)$, i.e.,
$x=y-z$, whence $\langle y_1-z\rangle\in F$, i.e.,
$\Dom(F)=P_1$.

4. Similarly, we can prove that $\Rng(F)=P_2$.

5. We have proved that $F$ is a~bijection between
the modules $P_1$~and~$P_2$.
Now we only need to show that $F$ is a~homomorphism, i.e., that
$\langle y_1,z_1\rangle,\langle y_2,z_2\rangle \in F$ implies
$$
\langle \alpha_1y_1+\alpha_2y_2,\alpha_1z_1+\alpha_2z_2\rangle \in F.
$$
Indeed, $\langle y_1,z_1\rangle, \langle y_2,z_2\rangle \in F$ implies
\begin{multline*}
y_1+z_1,y_2+z_2\in P\Rightarrow
\alpha_1(y_1+z_1)+\alpha_2(y_2+z_2)\in P\\
\Rightarrow
(\alpha_1y_1+\alpha_2y_2)+(\alpha_1z_1+\alpha_2z_2)\in P\Rightarrow
\langle \alpha_1 y_1+\alpha_2 y_2,\alpha_1 z_1+\alpha_2 z_2\rangle \in F.
\end{multline*}

Therefore modules $P_1$~and~$P_2$ that satisfy our formula
do not intersect and are isomorphic. Conversely, if
two modules $P_1$~and~$P_2$ do not intersect and are isomorphic,
then they satisfy our formula. Hence we shall
denote it by $P_1\cong_d P_2$.

Suppose that modules $P_1$~and~$P_2$ ``are not too big'', i.e.,
there exist modules $P_1'$~and~$P_2'$ such that
$P_1\cap P_1'=P_2\cap P_2'=\varnothing$,
the module~$P_1'$ contains a~submodule which
is isomorphic to~$P_1$,
and the module~$P_2'$ contains a~submodule which is isomorphic to~$P_2$.
Then the formula
$$
\exists P\losp \exists P'\losp
(P\cong_d P_1 \logic\land P'\cong_d P_2 \logic\land P\cong_d P')
$$
holds if and only if the modules $P_1$~and~$P_2$
are isomorphic.

We know that a~module~$P$ is projective if and only if
it is isomorphic to a~direct summand of a~free module.
Therefore the formula
$$
\mathrm{Proj}(P) \prisv \exists Q\losp (V=P\oplus Q)
$$
defines in the space $P(V)$ all projective modules.

Consider some projective module~$P$. Its submodule~$M$
will be called a~\emph{maximal} submodule of the module~$P$
($M=\max (P)$) if the formula
$$
\forall P'\losp
(M\subset P' \logic\land P'\subset P\Rightarrow P'=M \logic\lor P'=P)
$$
holds.
For every finitely generated module~$P$ there exists
a~maximal submodule~$M$.

Let some projective module~$P$ and its maximal submodule~$M$
be fixed.

The formula $X\subset_\circ Y$ will denote that the module~$X$
is a~direct summand of the module~$Y$.

Consider a~pair of modules $\langle X,Y\rangle$
satisfying the following formula:
\begin{align*}
& \overline{\tSum}_{P,M}(X,Y)\prisv
Y\subset X \logic\land
\exists Q\losp \exists Q'\losp
(
Q\oplus P=X \logic\land Q\cong X \logic\land Q'\oplus M=Y
\logic\land Q'\cong Y\\
& \quad \logic\land
\forall N \subset_\circ X\losp
(N\cong P \Rightarrow N\cap Y \cong M \logic\land (N\cap Y) \subset_\circ Y)
)\\
& \quad
\logic\land \forall Z\losp
(Z\subset X \logic\land \forall N\losp
(N\subset_\circ Z\Rightarrow N\not\cong P) \Rightarrow Z\subset Y).
\end{align*}
Let us see which modules $X$~and~$Y$ satisfy the formula
$\overline{\tSum}_{P,M}$.

From the formula $\exists Q\losp (Q\oplus P=X \logic\land Q\cong X)$
we see that the module~$P$ is a~direct summand of the module~$X$
and the complement~$Q$ is isomorphic to~$X$. Therefore, there exist
some infinite cardinal number~$\alpha$ and modules $X_1$~and~$X_2$
such that $X_1\oplus X_2=X$, $X_1\cong P^{(\alpha)}$,
and the module~$P$ is not isomorphic to any direct summand of
the module~$X_2$. The part of the formula
$$
\forall Z\losp (Z\subset X \logic\land
(\forall N\losp
(N\subset_\circ Z\Rightarrow N\not\cong P)\Rightarrow Z\subset Y)
$$
shows that if $Z$ is some submodule of the module~$X$ such that
the module~$P$ is not isomorphic to its direct summand, then
$Z$ is also a~submodule in~$Y$. If we set $X_2\prisv Z$, then
$X_2\subset Y$. Take an arbitrary
$y\in Y$. Since $y\in X$, it follows that $y=x_1+x_2$, where
$x_1\in X_1$, $x_2\in X_2$.
Since $X_2\subset Y$, it follows that
$x_1\in Y$, i.e., $Y=(X_1\cap Y)\oplus X_2$.

Now other conditions imply $X_1\cap Y\cong M^{(\alpha)}$.
Therefore, if $X$~and~$Y$ satisfy the formula $\overline{\tSum}_{P,M}(X,Y)$,
then there exist a~module~$Q$ and an infinite
cardinal number~$\alpha$ such that
$X\cong Q\oplus P^{(\alpha)}$ and $Y\cong Q\oplus M^{(\alpha)}$.
The inverse implication is clear if the module~$X$ is ``not too big''.

Now consider the formula
\begin{multline*}
\tSum_{P,M}^\omega(X,Y)\prisv
\forall Z\losp \forall T\losp
(\overline{\tSum}_{P,M}(Z,T)\\
\Rightarrow
\exists X_1\losp \exists X_2\losp \exists Y'\losp
X_1\oplus X_2=Z \logic\land X_1\cap T=Y' \logic\land
X_1\cong X \logic\land Y'\cong Y)
\logic\land \overline{\tSum}_{P,M}(X,Y).
\end{multline*}
The subformula $\overline{\tSum}_{P,M}(X,Y)$ implies that
$X\cong Q\oplus P^{(\alpha)}$ and
$Y\cong Q\oplus M^{(\alpha)}$ for some cardinal number~$\alpha$.
The first part of the formula implies that $X$ is a~direct
summand in every submodule of the form $Q'\oplus P^{(\beta)}$
($\beta$~is an infinite cardinal number) and, therefore,
in the module~$P^{(\omega)}$. Hence
$\alpha=\omega$, the module~$Q$ is projective and countably
generated.

Now consider the formula
$$
\tSum_{P,M}^{\Fin}(X,Y)\prisv
\neg \tSum_{P,M}^\omega(X,Y) \logic\land
\exists X',Y'\losp (\tSum^\omega (X',Y') \logic\land
\exists X''\losp (X'=X\oplus X' \logic\land Y=X\cap Y')).
$$
Every module~$X$ satisfying the formula $\tSum_{P,M}^{\Fin}(X,Y)$
is a~direct summand in the module $Q\oplus P^{(\omega)}$, i.e.,
has the form $Q'\oplus P^{(n)}$ (possibly, $n=0$, but $n\in \omega$), and
$Q'$ is a~direct summand of the module~$Q$. Let modules
$X_1$, $X_2$, $Y_1$, and $Y_2$ be
such that $\tSum_{P,M}^{\Fin}(X_1,Y_1)$ and
$\tSum_{P,M}^{\Fin}(X_2,Y_2)$. If
$$
\exists X_1',Y_1'\losp (\tSum_{P,M}^{\Fin}(X_1',Y_1') \logic\land
X_1'\cong X_1 \logic\land Y_1'\cong Y_1)
$$ and
\begingroup
\setlength{\multlinegap}{0pt}
\begin{multline*}
\forall P'\losp \forall M'\losp
(P'\cong P \logic\land M'\cong M \logic\land M'=\max(P') \logic\land
\exists P''\losp (P'\oplus P''=X_1' \logic\land P'\cap Y_1'=M')
\Rightarrow P'\subset_\circ X_1'\cap X_2')\\
\logic\land
\forall P'\losp \forall M'\losp
(P'\cong P \logic\land M'\cong M \logic\land M'=\max(P') \logic\land
\exists P''\losp (P'\oplus P''=X_2' \logic\land P'\cap Y_2'=M')
\Rightarrow P'\subset_\circ X_1'\cap X_2'),
\end{multline*}
\endgroup%
then we shall call the pairs $(X_1,Y_1)$ and $(X_2,Y_2)$
\emph{equivalent} (notation: $(X_1,Y_1)\sim (X_2,Y_2)$). It is clear that
if $(X_1,Y_1)\sim (X_2,Y_2)$, $X_1\cong Q_1\oplus P^{(n_1)}$, and
$X_2\cong Q_2\oplus P^{(n_2)}$, then $n_1=n_2$. We shall denote
the equivalence classes of such pairs by $\Cl_{P,M}^n$.

For two classes $\Cl_{P,M}^m$ and $\Cl_{P,M}^n$ we shall write
$\Cl_{P,M}^m<\Cl_{P,M}^n$ if
$$
\forall (X_1,Y_1)\in \Cl_{P,M}^m\losp \exists (X_2,Y_2)\in \Cl_{P,M}^n\losp
\exists X_3\losp (X_1\cong X_3 \logic\land X_1\subset_\circ X_2).
$$
It is clear that the condition $\Cl_{P,M}^m< \Cl_{P,M}^n$
is equivalent to the condition $m<n$.

Similarly to modules of the form $Q\oplus P^{(n)}$,
with the help of the same formula,
we can introduce the equivalence classes $\Cl_{P,M}^{(\alpha)}$ also for
infinite cardinal numbers~$\alpha$ and we can introduce the relation~$<$
between them.

A~module~$P$ will be called a~\emph{generator} if
$$
\exists \Cl_{P,M}^\alpha\losp
\forall V_1\losp \forall V_2\losp \forall X\losp \forall Y\losp
(V_1\oplus V_2=V \logic\land (X,Y)\in \Cl_{P,M}^\alpha \Rightarrow
V_1\subset_\circ X \logic\lor V_2\subset_\circ X).
$$
This formula will be denoted by $\Gener(P)$.

The formula
$$
\Pret(P)\prisv \mathrm{Proj}(P) \logic\land \Gener(P) \logic\land
\exists M\subset P\losp (M=\max(P))
$$
holds for all projective generators that have maximal submodules,
and it necessarily holds for all progenerators.

The formula
$$
\mathrm{FDSum}_{P,M}(X)\prisv
\exists \Cl_{P,M}^{(n)}\losp \exists Y (X,Y)\in \Cl_{P,M}^{(n)}
\logic\land
\forall X',Y'\losp (X',Y')\in \Cl_{P,M}^{(n)} \Rightarrow
(X,Y)\subset_\circ (X',Y')
$$
defines for a~given~$n$ a~module $Q\oplus P^{(n)}$ with a~submodule
$Q\oplus M^{(n)}$ such that for every pair
$(Q'\oplus P^{(n)},\allowbreak Q'\oplus M^{(n)})$
the module $Q\oplus P^{(n)}$ is a~direct summand in
$Q'\oplus P^{(n)}$ and $Q'\oplus M^{(n)}\subset Q\oplus P^{(n)}$.
Consider
the pair $(P^{(n)},M^{(n)})$
as the modules $Q'\oplus P^{(n)}$
and $Q'\oplus M^{(n)}$.
Then
$P^{(n)}\cong P^{(n)}\oplus Q$ and
$M^{(n)}\cap P^{(n)}\oplus Q=M^{(n)}\oplus Q$.

This formula defines all modules of the form $P^{(n)}$, where $n\in \omega$,
and some other \emph{finitely generated} modules.

Every projective finitely generated module is a~direct summand
of the module~$R^{(n)}$ for some $n\in \omega$. Therefore, if
$P$ is a~finitely generated projective module, then for every generator~$S$
$$
P\oplus Q\cong S^{(m)}
$$
for some $m\in \omega$ and some module~$Q$. But if a~module~$P$
is not finitely generated, but is a~projective generator,
then it can not be embedded into~$R^{(n)}$ for any $n\in \omega$.

Therefore, the formula
$$
\Proobr(P)\prisv \Pret(P) \logic\land \forall S\losp
(\Pret(S)\Rightarrow \exists M\losp\exists X\losp
(\mathrm{FDSum}_{S,M}(X) \logic\land P\subset_\circ X))
$$
holds for progenerators, and only for them.

Note that, selecting some fixed progenerator~$P$ with the help
of the formula $\Proobr()$, we have also the set
of all almost free modules of ranks ${\le}\,\varkappa$
over the ring~$R$.

\subsection[The Ring $\protect\Endom_RP$]{The Ring
$\boldsymbol{\Endom_RP}$}\label{ss4.2}
In this section, we assume that we have some
fixed progenerator~$P$.

Let $P_1$, $P_2$, and $P_3$ be three mutually disjoint modules,
and let each of them be isomorphic to the module~$P$.

A~module $U_{1,2}$ is defined by the formula
$$
U_{1,2}\subset P_1\oplus P_2 \logic\land
P_1\subset U_{1,2}\oplus V_2 \logic\land
V_2\subset U_{1,2}\oplus V_1.
$$
As we know in this case the module~$U_{1,2}$ consists of sums
$e+f(e)$, where $e\in P_1$ and $f\colon P_1\to P_2$ is an isomorphism
between the modules $P_1$~and~$P_2$.
Evidently one can suppose that the isomorphism~$f$
coincides with the isomorphism which identifies
the modules $P_1$~and~$P_2$, i.e., the module~$U_{1,2}$
consists of vectors $f_1(e)+f_2(e)$,
where $f_1\colon P\to P_1$ and $f_2\colon P\to P_2$
are isomorphisms identifying the modules $P$, $P_1$, and~$P_2$.

Similarly, let us introduce a~module~$U_{2,3}$, consisting of
vectors of the form $f_2(e)+f_3(e)$.

A~module $U_{1,2,3}$ will be introduced by the formula
$$
U_{1,2,3} \prisv (P_1\oplus U_{2,3})\cap (P_3\oplus U_{1,2}).
$$
If $v\in U_{1,2,3}$, then $v\in P_1\oplus U_{2,3}$, i.e.,
$$
v=f_1(e)+f_2(e')+f_3(e'),
$$
and $v\in P_3\oplus U_{1,2}$ implies
$$
v=f_1(g)+f_2(g)+f_3(g').
$$
Therefore,
$$
f_1(e)+f_2(e')+f_3(e')=f_1(g)+f_2(g)+f_3(g'),
$$
and so
$$
v=f_1(e)+f_2(e)+f_3(e).
$$

A~module $U_{1,3}$ is introduced by the formula
$(P_1\oplus P_3)\cap (U_{1,2,3}\oplus P_2)$.

Thus we have the modules generated by the elements $f_1(e)=e_1$,
$f_2(e)=e_2$, $f_3(e)=e_3$, $f_1(e)+f_2(e)=e_1+e_2$,
$f_2(e)+f_3(e)=e_1+e_3$,
$f_1(e)+f_2(e)+f_3(e)=e_1+e_2+e_3$ for $e\in P$.

Introduce now a~set~$V_q^3$ of modules with the help of the formula
$$
V_q^3\subset U_{1,2}\oplus P_3 \logic\land U_{1,2}\subset V_q^3\oplus P_3.
$$
Since $V_q^3\subset U_{1,2}\oplus P_3$, it follows that $v\in V_q$ implies
$$
v=f_1(e)+f_2(e)+f_3(e').
$$
From $U_{1,2}\subset V_q\oplus P_3$ it follows that for every
$e\in P$ there exists $v=f_1(e)+f_2(e)+f_3(e')$.

For every $e\in P$ there exists a~unique $e'\in P$
such that $f_1(e)+f_2(e) +f_3(e')\in V_q^3$.

It is clear that the correspondence which maps an element~$e$
into the element~$e'$
is a~homomorphism of the module~$P$ into itself.
We shall denote it by~$q$.
For every module~$V_q^3$ by $W_q^{1,3}$ we denote the module
defined by the formula
$$
W_q^{1,3}\subset (P_1\oplus P_3)\cap (V_q^3\oplus P_2)
\logic\land P_1\subset W_q^{1,3}\oplus P_3.
$$
If $w\in W_q^{1,3}$, then $w\in P_1\oplus P_3$ implies
$w=f_1(e)+f_3(e')$, and $w\in V_q^3\oplus P_2$ implies
$$
w=f_1(e'')+f_2(e'')+f_3(qe'')+f_2(e''').
$$
Therefore $w=f_1(e)+f_3(qe)$.

Similarly we can introduce a~module $W_q^{2,3}$, consisting of
vectors
$$
w=f_2(e)+f_3(qe).
$$
For given $V_q^3$~and~$V_r^3$ consider the module~$V$ defined
by the formula
$$
V\prisv (U_{1,2}\oplus P_3)\cap (W_q^{1,3}\oplus W_r^{2,3}).
$$
If $v\in V$, then, on the one hand,
$$
v=f_1(e)+f_2(e)+f_3(e'),
$$
and on the other hand,
$$
v=f_1(e'')+f_3(qe'')+f_2(e''')+f_3(re''').
$$
We see that $E''=e'''=e$, i.e.,
$$
v=f_1(e)+f_2(e)+f_3(qe)+f_3(re)=f_1(e)+f_2(e)+f_3((q+r)e),
$$
whence
$$
V=V^3_{q+r}.
$$

Hence, on the set $\{ V_q\mid q\in \Endom_RP\}$ of modules we have
the operation of addition $\langle V_q,V_r\rangle \mapsto V_{q+r}$.
It is clear that in this case we also have the operation of
taking an opposite element
$$
V_q\mapsto V_{-q}.
$$

By $X_q^{2,3}$ we denote the module
$$
(W_q^{2,3}\oplus P_2)\cap U_{2,3}.
$$
It consists of vectors of the form
$$
f_2(qe)+f_3(qe),\quad e\in P.
$$

Now consider the module $W$ defined by the formula
$$
W\subset P_2\oplus P_3 \logic\land
P_3\subset P_2\oplus W \logic\land
X_q^{2,3}=
(((W\oplus P_3)\cap P_2)\oplus ((W_q^{2,3}\oplus P_2)\cap P_3))\cap U_{2,3}.
$$
It is easy to see that such a~module consists of vectors of the form
$$
f_3(e)+f_2(qe).
$$
We shall denote it by $W_q^{3,2}$.

The module
$$
(W_q^{3,2}\oplus P_1)\cap (U_{1,3}\oplus P_2)
$$
will be denoted by~$V_q^2$. It consists of vectors of the form
$$
f_1(e)+f_2(qe)+f_3(e).
$$
The module $V_q^2\oplus P_3\cap P_1\oplus P_2$ is denoted by $W_q^{1,2}$
and consists of vectors of the form $f_1(e)+f_2(qe)$.

If we have a~module $W_q^{1,2}$, then the formula
$$
(W_{q'}^{1,2}\oplus W_q^{1,3})\cap U_{2,3}=X_q^{2,3}
$$
gives $q'=q$, i.e., having a~module $W_q^{1,2}$, we automatically
have the module $W_q^{1,3}$, and, therefore, the module~$V_q^3$.

Now, writing the formula
$$
W_s^{1,2}=(W_q^{3,2}\oplus W_{-r}^{1,3})\cap (P_1\oplus P_2),
$$
we shall have for $w\in W_s^{1,2}$
$$
w=f_1(e)+f_3(-re)+f_3(e')+f_2(qe)=f_1(e'')+f_2(e'').
$$
Thus we have
$$
f_3(-re)+f_3(e')=0,
$$
i.e., $e'=re$, and so
$$
w=f_1(e)+f_2(qre),
$$
i.e., $s=qr$.

Therefore, given two modules $V_r^3$~and~$V_q^3$
we can construct the module~$V_{qr}^3$,
i.e., on the set $\{ V_q^3\mid q\in \Endom_RP\}$ we have introduced
the operation of addition and multiplication in such a~way that
it becomes isomorphic to the ring $\Endom_RP$.

\subsection[Construction of the Ring $\protect\Endom_RV$]{Construction
of the Ring $\boldsymbol{\Endom_RV}$}\label{ss4.3}
For a~given progenerator~$P$ select in the module~$V$ two
disjoint submodules $V_1$~and~$V_2$ and one equivalence class
$\Cl_{P,M}^\alpha$ which is maximal among all other $\Cl_{P,M}^\beta$.
It is clear that in this case $\alpha=\varkappa$. Let, further,
$V_1\oplus V_2\oplus P=V$.

Let $V_1=Q_1\oplus \sum\limits_{i\in \varkappa} P_i$
and
$V_2=Q_2\oplus \sum\limits_{i\in \varkappa} P_i'$, where for every
$i\in \varkappa$
$$
P_i\cong P_i'\cong P.
$$
Fix isomorphisms
\begin{align*}
f_i &\colon P\to P_i,\\
f_i'&\colon P\to P_i'.
\end{align*}

Let a~formula $\Endom(X)$ state about a~module~$X$ the following.

1. $\forall T\losp
(T\subset_\circ V_1 \logic\land
T\cong P\Rightarrow \exists T'\losp
(T'\subset_\circ V_2 \logic\land T'\cong P \logic\land
\exists V_q^3(P,T,T')\subset_\circ X))$, i.e.,
for every direct summand~$P_i$ of the module~$V_1$ there exists a~direct
summand~$P'$ (a~linear combination of some~$P_i'$)
of the module~$V_2$ such that
for some $q\in \Endom_RP$ the module
$$
\{ e +f_i(e)+f'(qe)\mid e\in P\}
$$
is a~direct summand of the module~$P$.

2. $X\cap V_2=0$, and it implies that for every direct summand~$P_i$
of the module~$V_1$ there exists only one direct summand~$P'$
of the module~$V_2$ such that the module
$$
\{ e+f_i(e)+f'(qe)\mid e\in P\}
$$
is a~direct summand of the module~$X$ for some $q\in \Endom_RP$.

3. $X\cap P=0$.
Such a~module presents an endomorphism of the module~$P^{(\varkappa)}$
over the ring $\Endom_RP$ in the following form.

For every vector $v\in P^{(\varkappa)}$ there exists~$P'$
(a~direct summand of the module $P^{(\varkappa)}$)
which is isomorphic to~$P$ and such that $v\in P'$.
By condition~1, in the module~$V_2$ there exists a~direct summand~$P''$,
and also there exists an endomorphism $q\in \Endom_RP$ such that
$V_q^3(P,P',P'')\subset X$.
Then the module~$V_q^3$ contains a~unique element
$$
(f')^{-1}(v)+v+ f''(q(f')^{-1}(v)) .
$$
We assume that $X(v)\prisv f''(q(f')^{-1}(v))$.
We show that the obtained mapping is well defined and linear.

Indeed, the simplicity of decomposition
follows from condition~2. Check the linearity.

If $v_1,v_2\in P_i$ for some $i\in \varkappa$, then for every $q_1,q_2
\in R$ the condition $X(q_1v_1+q_2v_2)=q_1X(v_1)+q_2 X(v_2$ follows from
the linearity of the corresponding endomorphism~$q$:
\begin{align*}
& V_q^3(P,P_i,P')\subset X
\Rightarrow \begin{cases}
f_i^{-1}(v_1)+v_1+f'(qf_i^{-1}(v_1))\in X,\\
f_i^{-1}(v_2)+v_2+f'(qf_i^{-1}(v_2))\in X
\end{cases}\\
& \quad \Rightarrow \begin{cases}
q_1f_i^{-1}(v_1)+q_1v_1+q_1f'(qf_i^{-1}(v_1))\in X,\\
q_2f_i^{-1}(v_2)+q_2v_2+q_2f'(qf_i^{-1}(v_2))\in X
\end{cases}\\
& \quad \Rightarrow
q_1f_i^{-1}(v_1)+q_2 f_i^{-1}(v_2)+ q_1v_1+q_2v_2+q_1f'(qf_i^{-1}(v_1))+
q_2f'(qf_i^{-1}(v_2))\in X\\
& \quad \Rightarrow
f_i^{-1}(q_1v_1+q_2v_2)+(q_1v_1+q_2v_2)+ f'(q_1 q(f_i^{-1}(v_1))+
q_1q(f_i^{-1}(v_2)))\in X,
\end{align*}
i.e.,
$$
X(q_1v_1+q_2v_2)=q_1X(v_1)+q_2 X(v_2).
$$

Two modules $X_1$ and $X_2$ satisfying the formula $\Endom(X)$
will be called \emph{equivalent} if
$$
\forall T\subset_\circ V_1\losp \forall S\subset_\circ V_2\losp
(V_q^3(P,T,S)\subset X_1\Leftrightarrow V_q^3(P,T,S)\subset X_2).
$$
We see that in every equivalence class there exists a~module of the form
$$
\sum_{i\in \varkappa} V_{q_i}^3(P,P_i,T_i),
$$
where $T_i$ is a~unique module for~$P_i$ such that
$$
V_{q_i}^3(P,P_i,T_i)\subset X.
$$

Consider some module~$X_0$ satisfying the formula $\Endom(X)$ and
such that $X_0\subset P\oplus V_1$. It is clear that the
endomorphism corresponding to the module~$X_0$
is the zero endomorphism of the module~$V_1$. We shall
now consider only modules~$X$ satisfying the formula
$$
\Endom^{X_0}(X) \prisv \Endom(X) \logic\land X\subset X_0\oplus V_2.
$$
Now define the sum of two modules $X_1$~and~$X_2$ satisfying
the formula $\Endom^{X_0}(X)$.
\begin{multline*}
(X=X_1+X_2)\prisv
\forall T\subset_\circ V_1\losp \forall V_q^3(P,T,S_1)\subset_\circ X_1\losp
\forall V_r^3(P,T,S_r)\subset_\circ X_2\\
(X_0\oplus V_2)\cap
(W_q^{1,3}(V_q^3(P,T,S_q))\oplus W_r^{2,3}(V_r(P,T,S_r)))
\subset_\circ X.
\end{multline*}
It is easy to see (compare to Sec.~4.2) that the module~$X$ satisfying
the formula $X=X_1+X_2$ is the sum of the endomorphisms $X_1$~and~$X_2$.

Now introduce some module~$X_e$ satisfying the formula $\Endom_{X_0}(X)$
and such that
\begin{multline*}
X_e\cap X_0=0 \logic\land
\forall S\subset_\circ V_2\losp
(S\cong P\\
\Rightarrow \exists T\subset_\circ V_1\losp
(\exists V_q(P,T,S)\subset X_e \logic\land
V_q^3(P,T,S)\text{ is an isomorphism between $T$ and $S$})).
\end{multline*}
It is clear that such a~module~$X_e$ establishes an isomorphism between
the modules $V_1$~and~$V_2$. Therefore,
$X_e$ will be the unit of $\Endom_RV$.

Now consider three modules $X_1$, $X_2$, and $X$
satisfying the formula $\Endom^{X_0}(X)$.
We need to define the formula $X=X_1\circ X_2$. We describe this formula
by words to understand its essence.

Let $V_q^3(P,T,S_q)\subset_\circ X_1$ for some $T\subset_\circ V_1$ and
$S_q\subset_\circ V_2$. As we have already said, for every
sum $v+f_T(v)+f_{S_q}(qv)$ we suppose that $X_1$ maps the vector
$f_T(v)\in T\subset V_1$ to the vector $f_{S_q}(qv)\in S_q\subset V_2$.

For a~given $S_q\subset_\circ V_2$ there exists a~unique~$T_q$ such that
$V_e(P,T_q,S_q)\subset_\circ X_e$. For an arbitrary vector $v\in P$, if
$v+f_{T_q}(v)+f_{S_q}(v)\in V_e(P,T_q,S_q)$,
then the vectors $f_{T_q}(v)$ and
$f_{S_q}(v)$ coincide if we identify $V_1$~and~$V_2$, i.e.,
$X_1$ maps $f_T(v)$ to $f_{T_q}(qv)$.

Then, for a~given $T_q \subset_\circ V_1$ there exists a~unique
$S_{rq}\subset_\circ V_2$ such that
$$
V_r(P,T_q,S_{rq})\subset_\circ X_2.
$$
If
$$
v+f_{T_q}(v)+f_{S_{rq}}(rv)\in V_r(P,T_q,S_{rq}),
$$
then the mapping~$X_2$ maps the vector $f_{T_q}(v)\in T_q\subset V_1$
to the vector $f_{S_{rq}}(rv)$, i.e., the composition $X_2X_1$
maps the vector $f_{T}(v)$ to the vector $f_{S_{rq}}(zqv)$, i.e.,
the mapping~$X$ is the composition of the mappings $X_1$~and~$X_2$
if and only if for every $T\subset_\circ V_1$
$$
V_3(P,T,S_{rq})\subset_\circ X,
$$
and with it $V_s(P,T,S_{rq})$ consists of vectors of the form
$$
v+f_T(v)+f_{S_{rq}}(rqv).
$$
We can easily make certain that the formula
$$
(V_q(P,T,S_q)\oplus V_e(P,T_q,S_q))\cap X_0
=(V_r(P,T_q,S_{rq})\oplus V_s(P,T,S_q))\cap X_0
$$
holds, and therefore we have a~formula which is equivalent
to the formula
$$
X=X_2\circ X_1.
$$
Thus, in the lattice of submodules of the module~$V$ we have interpreted
a~ring which is isomorphic to the ring $\Endom_{\Endom_RP}V$.
Consequently, as before, we have
that if two lattices of submodules $P(R_1,V_1)$ and $P(R_2,V_2)$
are elementarily equivalent, then
for some progenerators $P_1$~and~$P_2$
the rings
$\Endom_{\Endom_{R_1}P_1}(V_1)$ and
$\Endom_{\Endom_{R_2}P_2}(V_2)$
are also elementarily equivalent, and, therefore,
the rings $\Endom_{R_1}V_1$ and $\Endom_{R_2}V_1$
are elementarily equivalent.
Now we see that we have proved the following theorem.

\begin{theorem}\label{poject1}\label{th1_4.3}
For free modules $V_1$~and~$V_2$ of infinite ranks
over arbitrary rings $R_1$~and~$R_2$, respectively,
elementary equivalence of the lattices of submodules $P(V_1)$ and
$P(V_2)$ implies elementary equivalence of the endomorphism rings
$\Endom_{R_1}(V_1)$ and $\Endom_{R_2}(V_2)$.
\end{theorem}

\subsection{The Inverse Theorem}\label{ss4.4}
Now we need to prove the inverse theorem.

\begin{theorem}\label{project2}\label{th2_4.4}
Suppose that $V_1$~and~$V_2$ are free modules of infinite ranks
$\varkappa_1$~and~$\varkappa_2$ over rings $R_1$~and~$R_2$, respectively,
and every submodule of the module $V_1$~\textup{(}$V_2$\textup{)}
has at most~$\varkappa_1$~\textup{(}$\varkappa_2$\textup{)}
generating elements
\textup{(}for example,
this is true if $\varkappa_1\ge |R_1|$
and $\varkappa_2\ge |R_2|$, or if $R_1$ and~$R_2$ are semisimple rings
or integral domains\textup{)}.
Then
$\Endom_{R_1}(V_1)\equiv \Endom_{R_2}(V_2)$
implies $P(V_1)\equiv P(V_2)$.
\end{theorem}

\begin{proof}
Suppose that we have an associative ring~$R$ with a~unit, an infinite
cardinal number~$\varkappa$, and a~free module
$V=V_\varkappa^R$ of rank~$\varkappa$ over~$R$.
Further, let every ideal of the ring~$R$ be generated by
at most~$\varkappa$ elements of the ring.

We want to interpret in the ring $\Endom_RV$ the space $P(V)$,
consisting of all submodules of the module~$V$, with the relation~$\subset$.
As before, by the word ``interpret'' we understand
existence of some algorithm mapping every formula~$\varphi$
of the first order language of the theory of projective spaces
to a~formula~$\tilde \varphi$ of the first order language
of the ring theory in such a~way that
the formula~$\varphi$ holds in $P(V)$ if and only if
$\tilde \varphi$ holds in $\Endom_R(V)$.

At the beginning we shall give an informal description of the translation.

1. We know that every object of the space $P(V)$ is a~submodule
of the module~$V$, but it is generated by at most~$\varkappa$
vectors of the module~$V$.
Each of these vectors is a~linear combination of some finite number
of elements of a~basis of the module~$V$, i.e., every such vector
can be written as a~column of a~matrix
which has only a~finite set of nonzero elements. If we write in this matrix
all generating vectors, we shall get a~matrix of size
$\varkappa\times \varkappa$, i.e., an element of
$\Endom_RV$. In the case where a~submodule is generated less
than~$\varkappa$ vectors, we extend the matrix by zero columns.
Two such matrices $X_1$~and~$X_2$
describe the same submodule of the module~$V$ if
$$
\exists A\losp \exists B\losp (X_1=X_2A \logic\land X_2=X_1B).
$$
In this case, the elements $X_1$~and~$X_2$ will be considered equivalent.

Therefore, every submodule of the module~$V$ maps to
the corresponding equivalence class of elements of the ring~$\Endom_RV$.

2. It is clear that the module~$Y_1$ generated by a~matrix~$X_1$
is a~submodule of the module~$Y_2$ generated by a~matrix~$X_2$
if and only if
$$
\exists A\losp (X_1=X_2A).
$$
This formula will be denoted by $X_1\subset X_2$.

From all these statements we obtain the statement of the theorem.
\end{proof}

\section{Elementary Equivalence of Automorphism Groups
of Modules of Infinite Ranks}\label{s5}
\subsection[An Isomorphism of Groups $\protect\Aut_R(V)$]{An
Isomorphism of Groups $\boldsymbol{\Aut_R(V)}$}\label{ss5.1}
In this section, we are based on the paper~\cite{1} of I.~Z.~Go\-lub\-chik
and A.~V.~Mikhalev.

Consider some ring $R$ and a~free module $V(=V_\varkappa^R)$
of infinite rank~$\varkappa$ over this ring.

Let $I_\varkappa$ be a~set of cardinality~$\varkappa$.

As above, by $\Endom_R(V)$ we shall denote the endomorphism ring
of the module~$V$,
and by $\Aut_R(V)$ we shall denote the automorphism group of the module~$V$.

Let, further, $E_R(V)$ be the group generated by the automorphisms
$E_{\gamma \beta}$ of the form
$$
v_\gamma\mapsto n_\gamma+rv_\beta ,\quad
\gamma,\beta\in I_\varkappa,\ \ \gamma\ne \beta,\ \ r\in R,
$$
and
$$
v_\alpha \mapsto v_\alpha,\quad \alpha\in I_\varkappa,\ \ \alpha \ne \gamma,
$$
where $\{ v_\alpha\}$ is a~basis of the module~$V$;
$D_R(V)$ is the diagonal group (the automorphisms
of the form $v_\gamma\mapsto r_\gamma v_\gamma$
$\forall \gamma\in I_\varkappa$);
$DE_R(V)$ is the group generated by $E_R(V)$ and $D_R(V)$.

A~subset $\{ e_{ij}\}_{i,j\in I_\varkappa}$ of the ring $\Endom_R(V)$
is called a~\emph{system of matrix units} if
\begin{enumerate}
\item
$e_{ij}\circ e_{st}=\delta_{js} e_{it}$ ($\delta_{js}$
is the Kronecker delta);
\item
for every $a\in \Endom_R(V)$ and every $k\in I$ there exist
$i_1,\dots,i_n\in I$ such that
$(e_{i_1i_1}+\dots+e_{i_ni_n}) a = a(e_{i_1i_1}+\dots+e_{i_ni_n}= a$.
\end{enumerate}

Let $I$ be an ideal of the ring~$R$;
$E_R(V,I)$ be the subgroup of the group $\Aut_R(V)$
generated by the automorphisms $1+e_{ij}\circ \lambda$, where
$\lambda\in I$, $i\ne j\in I_\varkappa$,
$\Aut_R(V,I)$ be the kernel of the canonical homomorphism
$\varphi_I\colon \Aut_R(V)\to \Aut_{R/I}(V)$, $C_R(V,I)$ be the inverse
image of the center in the homomorphism~$\varphi_I$. Let, further,
$[A,B]\equiv A^{-1}\circ B^{-1}\circ A\circ B$.

\begin{lemma}\label{mih-l-1}
Let $R$ be an associative ring with $1/2$, $N$~and~$M$ be normal
subgroups of the group $\Aut_R(V)$ such that
$N\cap M=\{ 1\}$ and $NM=\Aut_R(V)$.
Then there exist ideals $I$ and $J$ of the ring~$R$ such that
$$
R=I\oplus J,\quad
E_R(V,I)\subseteq N\subseteq C_R(V,I),\quad
E_R(V,J)\subseteq M\subseteq C_R(V,J).
$$
\end{lemma}

\begin{proof}
By the condition,
\begin{equation}\label{e1}
(1-2e_{ii})=a_i\circ b_i,\quad a_i\in N,\ \ b_i\in M,
\end{equation}
for all $i\in I_\varkappa$. Since $N\cap M=\{1\}$ and
$[1-2 e_{11}, 1-2 e_{ii}]=1$,
it follows that
$[a_1,1-2e_{ii}]=1$. Since
$1/2\in R$, the element~$a_1$ is diagonal. This means that
$e_{ii}\circ a_1\circ e_{jj}=0$ for all $i\ne j\in I_\varkappa$.
The same holds also for~$b_1$. Let for all $i\in I_\varkappa$
\begin{equation}\label{e2}
e_{ii}\circ a_1\circ e_{ii}=\lambda_i,\quad
e_{ii}\circ b_1\circ e_{ii}=\mu_i.
\end{equation}
From \eqref{e2} it follows that
$$
a_1\circ (1-e_{12})\circ a_1^{-1}\circ (1+e_{12})=
1+(1-\lambda_1\lambda_2^{-1})\circ e_{12}\in N.
$$
Since $[1+\lambda e_{12}, 1+re_{2k}]=1+ \lambda r e_{1k}$
for all $\lambda,r\in R$ and $k\in I_\varkappa$,
it follows that
if the group~$N$ is normal, then $E_R(V,I)\subseteq N$,
where $I=R(\lambda_1-\lambda_2) R$. Similarly, $E_R(V,J)\subseteq M$,
where $J=R(\mu_1-\mu_2)R$.
From \eqref{e1} and \eqref{e2} it follows that
$$
\lambda_1\mu_1=-1,\quad \lambda_2\mu_2=1,\quad
\mu_1=-\lambda_1^{-1},\quad \mu_2=\lambda_2^{-1}.
$$
By the definition of ideals $I$~and~$J$,
$$
1-\lambda_1\lambda_2^{-1}\in I,\quad
1-\mu_1\mu_2^{-1}=1+\lambda_1^{-1}\lambda_2\in J,
$$
and
$$
\lambda_1(1+\lambda_1^{-1}\lambda_2)\lambda_2^{-1}=1+\lambda_1\lambda_2^{-1}
\in J.
$$
Consequently,
$1=1/2(1-\lambda_1\lambda_2^{-1}+1+\lambda_1\lambda_2^{-1})\in I+J$
and $R=I+J$. Further, $E_R(V,I\cap J)\subseteq N\cap M=\{ 1\}$,
and, therefore, $I\cap J=\{0\}$. Thus, $I\oplus J=R$.

If $a\in N$, then $a=a_1\circ a_2$,
where $a_1\in \Aut_R(V,I)$ and $a_2\in \Aut_R(V,J)$.
Further, we have $[a,E_R(V,J)]\subseteq N\cap M=\{ 1\}$.
Thus $a_2$ is a~central idempotent of $\Aut_R(V,J)$
and $N\subseteq C_R(V,I)$. Similarly, $M\subseteq C_R(V,J)$.
\end{proof}

The following lemma is basic in the proof.

\begin{lemma}\label{mih-l-2}
Let $R$ and $S$ be associative rings with $1/2$,
$I_1=I_\varkappa$ and $I_2=I_{\varkappa'}$ be infinite sets
of cardinalities $\varkappa$~and~$\varkappa'$, respectively,
$V=V_{I_1}^R$ and $V'=V_{I_2}^S$ be free modules
over the rings $R$~and~$S$ and the sets $I_1$~and~$I_2$, respectively,
$\{ e_{ij}\}_{i,j\in I_\varkappa}$
be a~system of matrix units of the ring $\Endom_R(V)$, and
$\varphi\colon \Aut_R(V)\to \Aut_S(V')$ be a~group isomorphism.
Then there exist
a~central idempotent $q\in \Endom_S (V')$ and systems of matrix units
$\{ f_{ij}\}_{i,j\in I_2}$ and $\{ h_{ij}\}_{i,j\in I_2}$ of the rings
$q\circ \Endom_s(V')$ and $(1-q)\circ \Endom_S(V')$, respectively,
such that
$$
\varphi (1-2 e_{ii})=(q-2f_{ii})-(1-q-2 h_{ii}),\quad i\in I_1.
$$
\end{lemma}

\begin{proof}
Consider $b_i\equiv \varphi(1-2 e_{ii})$. We know that $b_i^2=1$.
Therefore, for $f_i\equiv 1/2(1-b_i)\in \Endom_S(V')$ we have
$f_i^2=f_i$. Define such $f_i$ for all $i\in I_1$.
We shall get
\begin{equation}\label{e3}
\varphi(1-2e_{ii})=1-2 f_i.
\end{equation}
Since $1-2e_{11}$ and $1-2e_{22}$ commute,
$b_1$~and~$b_2$ also commute,
and, thus, $f_1$~and~$f_2$ commute. Thus, $(1-2 f_1 f_2)^2=1$, i.e.,
$1-2f_1 f_2\in \Aut_S(V')$.
Set
\begin{equation}\label{e4}
1-2e=\varphi^{-1}(1-2f_1f_2).
\end{equation}
Then $e\in \Endom_R(V)$, $e^2=e$, and from~\eqref{e3} it follows that
if $[a, 1-2 e_{ii}]=1$ for $i=1,2$, then
\begin{equation}\label{e5}
[a,1-2e]=1;
\end{equation}
if $b(1-2e_{11})b^{-1}=1-2 e_{22}$ and $b(1-2 e_{22})b^{-1}=1-2 e_{11}$,
then
\begin{equation}\label{e6}
[b,1-2e]=1.
\end{equation}
Applying \eqref{e5} and \eqref{e6}, we get
\begin{equation}\label{e7}
(1-2e)=\varepsilon_1(e_{11}+e_{22}+\varepsilon_2 (1-e_{11}-e_{22})),
\end{equation}
where $\varepsilon_1,\varepsilon_2\in R$,
$\varepsilon_1^2=\varepsilon_2^2=1$, and
the elements $\varepsilon_1$,~$\varepsilon_2$ are permutable with
all invertible elements of the ring~$R$. Then
\begin{equation}\label{e8}
\varepsilon_1=1-2 e_1,\quad \varepsilon_2=1-2e_2,
\end{equation}
and $e_1$,~$e_2$ are central idempotents of the ring~$R$.

Set
\begin{equation}\label{e9}
N\equiv \varphi(\Aut_R (V,e_2R)),\quad
M\equiv \varphi(\Aut_R(V, (1-e_2)R)).
\end{equation}
By Lemma~\ref{mih-l-1},
\begin{equation}\label{e10}
E_S(V',I)\subseteq N\subseteq C_S(V',I),\quad
E_S(V',J)\subseteq M\subseteq C_S(V',J),
\end{equation}
then $\Endom(I^{(\varkappa')})=(1-q)\Endom(V')$,
$\Endom(J^{(\varkappa')}=q\Endom(V')$,
and $q$ is some central idempotent of the ring $\Endom(V')$.
From \eqref{e7}~and~\eqref{e8} it follows that
\begin{align*}
& e_{11}+e_{22}+(1-2e_2)(1-e_{11}-e_{22})\in \Aut_R(V, e_2R),\\
& -e_{11}-e_{22}+(1-2e_2)(1-e_{11}-e_{22})\\
& \quad = -(e_{11}+e_{22}+(1-2(1-e_2))(1-e_{11}-e_{22}))
\in (-1)\Aut_R(V, (1-e_2)R)
\end{align*}
and, therefore,
\begin{equation}\label{e11}
1-2e\in C_R(V,e_2R),\quad
(1-2e_{11})(1-2e_{22})(1-2e)\in C_R(V,(1-e_2)R).
\end{equation}
From \eqref{e3}, \eqref{e4}, \eqref{e9}, \eqref{e10}, and \eqref{e11}
it follows that $1-2f_1f_2=a+b$, where $a\in \Endom(I^{(I_2)})$ and
$b\in \Endom(J^{(I_2)}$. Consequently, $b$ is a~central element
of the ring $\Endom(J^{(I_2)})$ and
$a_1\equiv a(1-2f_1)(1-2f_2)$
is a~central element of the ring $\Endom(I^{(I_2)})$.
Further, $(1-2f_1f_2)^2=1$, and, therefore, $b^2=q$, $a_1^2=a^2=1-q$,
$a_1 =1-q-2q_2$, and $b=q-2q_1$, where $q$, $q_1$, and $q_2$
are central idempotents of the rings $\Endom_S(V')$, $q\Endom(V')$,
and $(1-q)\Endom(V')$, respectively.
Thus,
\begin{equation}\label{e12}
(1-2f_1f_2)=(q-2q_1)+(1-q-2q_2)(1-2f_1)(1-2f_2).
\end{equation}
We shall show that $q_1=0$ and $q_2=1-q$. Indeed, multiplying
the equality~\eqref{e12} by~$q_1$, we get
$q_1({1-2f_1f_2})=-q_1$, i.e., $q_1f_1f_2=q_1$.

Multiplying the last equality by~$f_1$, we see that
$q_1 f_1f_2=q_1f_1$ and $q_1f_1=q_1$. Similarly, $q_1f_2=q_1$.

Hence, $q_1(1-2f_1)(1-2f_2)=q_1$ and, according to~\eqref{e3},
\begin{equation}\label{e13}
q_1\varphi(1-2e_{11}-2e_{22})=q_1.
\end{equation}
Since
$$
\left[
\begin{pmatrix}
1& -1/2 r\\
0& 1
\end{pmatrix},
\begin{pmatrix}
1& 0\\ 0& -1
\end{pmatrix}
\right]=
\begin{pmatrix}
1& r\\0&1
\end{pmatrix},
$$
we have that a~normal divisor of the group $\Aut_R(V)$ containing
the matrix $1-2e_{11}-2e_{22}$ contains also the subgroup $E_R(V)$.

Therefore, from~\eqref{e13} it follows that
\begin{equation}\label{e14}
\varphi(E_R(V))\subset \Aut_S(V',(1-q_1)S).
\end{equation}
By condition~\eqref{e12}, $q_1$~is a~central idempotent
of the ring $\Endom_S(V')$. By Lemma~\ref{mih-l-1},
$$
E_R(V,I_1)\subseteq \varphi^{-1}(\Aut_S(V',q_1S))\subseteq C_R(V,I_1).
$$
On the other hand, \eqref{e14} implies that
$$
\varphi^{-1}(\Aut_S(V',q_1S))\cap E_R(V)=\{ 1\}.
$$
Consequently, $I_1=\{ 0\}$, and the group $\varphi^{-1}(\Aut_S(V',q_1S))$
belongs to the center of the group $\Aut_R(V)$, i.e.,
\begin{equation}\label{e15}
q_1=0.
\end{equation}
Multiplying the equality~\eqref{e12} by $q_3\equiv 1-q-q_2$, we get
$$
(1-2f_1f_2)q_3=(1-2f_1)(1-2f_2)q_3
$$
and
\begin{equation}\label{e16}
2f_1f_2q_3=2f_1q_3+2f_2q_3-4f_1f_2q_3 .
\end{equation}
Multiplying the equality~\eqref{e16} by $1/2(1-f_1)$, we shall see
that $(1-f_1)f_2q_3=0$ and $f_2q_3=f_1f_2q_3$. Similarly,
$(1-f_2)f_2q_3=0$ and $f_1q_3=f_1f_2q_3=f_2q_3$. Hence
$$
2f_1f_2q_3=2f_1q_3+2f_2q_3-4f_1f_2q_3=4f_1q_3-4f_1q_3=0.
$$
Thus,
$$
f_1q_3=f_2q_3=f_1f_2q_3=0,\quad
q_3(1-2f_1)(1-2f_2)=q_3,
$$
and
\begin{equation}\label{e17}
q_3\varphi(1-2e_{11}-2e_{22})=q_3.
\end{equation}
Similarly as from the equality~\eqref{e13} we obtained
$q_1=0$, from the equality~\eqref{e17} we shall now find
\begin{equation}\label{e18}
0=q_3=1-q-q_2.
\end{equation}
From \eqref{e12}, \eqref{e15}, and \eqref{e18} it follows that
\begin{equation}\label{e19}
1-2f_1f_2=q-(1-q)(1-2f_1)(1-2f_2),\quad
f_1f_2q=0,\quad (1-f_1)(1-f_2)(1-q)=0.
\end{equation}
Since the group $\Aut_R(V)$ acts transitively on the set
$$
\{ 1-2e_{ii},1-2e_{jj}\}_{i\ne j;\ i,j\in I_1},
$$
from the conditions \eqref{e3} and \eqref{e19} we obtain
\begin{equation}\label{e20}
f_if_jq=0,\quad (1-f_i)(1-f_j)(1-q)=0
\end{equation}
for all $i,j\in I_1$, where $q$ is a~central idempotent
of the ring $\Endom_S(V')$ from condition~\eqref{e10}.

According to condition~\eqref{e20},
$\{ f_i q=1/2 (1-\varphi(1-2e_{ii}))q\}$ is an orthogonal
system of conjugate idempotents of the ring $q\Endom_S(V')$, and,
therefore, there exist elements $f_{ij}\in q \Endom_S(V')$ such that
$f_{ii}=q f_i$ and $f_{ij}f_{ks} =\delta_{jk} f_{is}$.

Now we show that if $a\in \Endom_S(V')$ and $m\in I$, then
there exist $i_1,\dots,i_n\in I_2$ such that
$$
f_{i_1i_1}+\dots+f_{i_ni_n} (ae_{mm})=
(ae_{mm})f_{i_1i_1}+\dots +f_{i_ni_n}=(ae_{mm})q.
$$
Fix some $a\in \Endom_S(V')$ and $m\in I$.
It is clear that in this case there exists a~set $i_1,\dots,i_n\in I_1$
such that $a$ commutes with the element
$\varphi(-1_{i_1,\dots,i_n})\equiv
\varphi\Bigl(\prod\limits_{1\le k\le n} (1-2e_{i_ki_k})\Bigr)$
and
$\varphi(-1_{i_1,\dots,i_n})\circ ae_{mm}=
ae_{mm}\varphi(-1_{i_1,\dots,i_n}) =-ae_{mm}$.
Then
$q\circ (-1_{i_1,\dots,i_n})=
\prod\limits_{1\le k\le i}(q-2f_{i_ki_k})=
q-2f_{i_1i_1}-\dots-2f_{i_ni_n}$, i.e.,
$(q-2f_{i_1i_1}-\dots -2f_{i_ni_n})ae_{mm}=
ae_{mm}(q-2f_{i_1i_1}-\dots -2f_{i_ni_n})=-ae_{mm}q$.
Therefore,
$2ae_{mm}q=(2f_{i_1i_1}+\dots+2f_{i_ni_n})ae_{mm}=
ae_{mm}(2f_{i_1i_1}+\dots+2f_{i_ni_n})$, i.e.,
$ae_{mm}(f_{i_1i_1}+\dots+f_{i_ni_n})=
(f_{i_1i_1}+\dots+f_{i_ni_n})ae_{mm}=ae_{mm}q$,
as required.

Thus we have shown that $\{ f_{ij}\}_{i,j\in I_2}$ is a~system
of matrix units of the ring $q\Endom_S(V')$. In a~similar way,
there exists
a~system of matrix units $\{ h_{ij}\}_{i,j\in I_2}$ of the ring
$(1-q)\Endom_S(V')$ such that $h_{ii}=(1-f_i)(1-q)$. Consequently,
$$
\varphi(1-2e_{ii})=1-2f_{ii}=(1-2f_{ii})q-(1-2(1-f_{ii}))(1-q)=
(q-2f_{ii})-(1-q-2h_{ii}).\qed
$$
\renewcommand{\qed}{}
\end{proof}

\begin{theorem}\label{mih-t-1}
Let $R$ and $S$ be associative rings with $1/2$,
$V=V_{\varkappa}^R$ and $V'=V_{\varkappa'}^S$ be free modules over
$R$~and~$S$ of infinite ranks $\varkappa$~and~$\varkappa'$ respectively, and
$\varphi\colon \Aut_R(V)\to \Aut_S(V')$ be a~group isomorphism.
Then there exist central idempotents $e$~and~$f$
of the rings $\Endom_R(V)$ and $\Endom_S(V')$, respectively,
a~ring isomorphism
$\theta_1\colon e\Endom_R(V)\to f \Endom_S(V')$, a~ring antiisomorphism
$\theta_2\colon (1-e)\Endom_R(V)\to (1-f) \Endom_S(V')$, and
a~group homomorphism
$\chi\colon DE_R(V)\to C(\Aut_S(V'))$ such that
$\varphi(A)=\chi(A)(\theta_1(A)+\theta_2(A^{-1}))$ for all $A\in E_R(V)$.
\end{theorem}

\begin{proof}
By Lemma~\ref{mih-l-2},
\begin{equation}\label{e21}
\varphi(1-2e_{ii})=(q-2f_{ii})-(1-q-2h_{ii}),
\end{equation}
where $q$ is a~central idempotent of the ring $\Endom_S(V')$,
$e_{ij}$, $f_{ij}$, and $h_{ij}$
are matrix units of the rings
$\Endom_R(V)$, $q\Endom_S(V')$, and $(1-q)\Endom_S(V')$, respectively.

Set
$$
f\equiv f_{11}+f_{22}+h_{11}+h_{22}.
$$

1. Let $\{ e_{ij}'\}_{i,j\in I_1}$ be some system of matrix units
of the ring $\Endom_R(V)$ and $\forall i\ne 1,2\losp (e_{ii}'=e_{ii})$.
Then
\begin{equation}\label{e22}
\varphi(1-2e_{ii}')=q-(-q)+x,\ \ \text{where}\ \ x\in f \Endom_S(V') f.
\end{equation}
By the condition, $[1-2e_{kk}',1-2e_{ii}']=1$ for $k=1,2$, $i\ne 1,2$.
By \eqref{e21} and~\eqref{e22},
\begin{equation}\label{e23}
\varphi(1-2e_{kk}')=1-2e_k+c_k,
\end{equation}
where $k=1,2$, $e_k\in f\Endom_S(V')f$, $c_k\in (1-f)\Endom_S(V')(1-f)$.
Note that
$$
(1-2e_{11}')(1-2e_{22}')=(1-2e_{11})(1-2e_{22}).
$$
According to the equalities \eqref{e21}, \eqref{e22}, and \eqref{e23},
$$
(f-2e_1)(f-2e_2)=-f
$$
and
\begin{equation}\label{e24}
e_1+e_2=f, \quad e_1e_2=0.
\end{equation}
By Lemma~\ref{mih-l-2}, there exists a~central idempotent~$q'$
of the ring $\Endom_S(V')$ such that
$$
(q'-2f_{ii}')-(1-q'-2h_{ii}')=\varphi(1-2e_{ii}').
$$
Consequently, for $k=1,2$ we have
\begin{align}
q'(1-\varphi(1-2e_{kk}'))(1-\varphi (1-2e_{33}'))&=0,\label{e25}\\
(1-q')(1+\varphi(1-2e_{kk}'))(1+\varphi(1-2e_{33}'))&=0.\label{e26}
\end{align}
Multiplying \eqref{e25} from the left side by $1-f$ and from the right
side by~$q$ and using the conditions
\eqref{e21} and~\eqref{e23}, we have that $q'c_k\cdot 2f_{33}=0$,
\begin{equation}\label{e27}
q'c_kf_{33}=0.
\end{equation}
Multiplying \eqref{e26} from the left side by~$f$
and from the right side by~$fg$ and using \eqref{e21},~\eqref{e23},
we have
$$
(1-q')2(f-e_k)2 fq=0.
$$
According to the equalities~\eqref{e24}, $f=e_1+e_2$.
Thus, $(1-q')e_k q=0$ and $(1-q')fq=0$.
Since $f=f_{11}+f_{22}+h_{11}+h_{22}$,
it follows that
$\Endom_S(V') f \Endom_S(V')=\Endom_S(V')$ and,
by the equalities $(1-q')fq=0$,
$$
0=(1-q')q \Endom_S(V') f \Endom_S(V')=(1-q')q\Endom_S(V').
$$
Therefore,
\begin{equation}\label{e28}
(1-q')q=0.
\end{equation}
From \eqref{e27} and~\eqref{e28} it follows that
$$
c_kf_{33}=c_kqf_{33}=q(q'c_kf_{33})+(1-q')qc_kf_{33}=0+0=0.
$$
Similarly, $c_kf_{ii}-0$ for all $i\in I_2$. By~\eqref{e23},
$c_k\in (1-f)\Endom_S(V')(1-f)$ and $c_k q=c_k(1-f)q$, i.e.,
\begin{equation}\label{e29}
c_kq=0.
\end{equation}
Multiplying the equality~\eqref{e25} from the left side by~$f$
and from the right side by $(1-q)f$, we have
$$
q'\cdot 2e_{kk}\cdot 2f(1-q)=0.
$$

Therefore,
\begin{equation}\label{e30}
q'(1-q)=0.
\end{equation}
From \eqref{e29} and \eqref{e30} it follows that $q=q'$, and from
\eqref{e21}, \eqref{e23}, and \eqref{e26} it follows that
$(2-2e_k+c_k)2h_{33}=0$.

Since $e_kh_{33}=e_kf(1-f)h_{33}=0$, it follows that $2h_{33}+c_kh_{33}=0$.
Similarly, $2h_{ii}+c_kh_{ii}=0$ for all $i\in I_2$. Thus
$c_k(1-q)=c_k(1-f)(1-q)$ and for any $i,j\in I_2$ for $i\ne 1,2$
we have $c_kq\cdot h_{ij}=c_k h_{ij}=-2 h_{ij}$, and for any
$j\in I_2$, $i=1,2$ we have  $c_kq\cdot h_{ij}=0$.
Thus, it is shown that
\begin{equation}\label{e31}
c_k(1-q)=-2(1-q)+2(1-q)f.
\end{equation}
From \eqref{e23}, \eqref{e29}, and \eqref{e31} it follows that
$1+c_k-q+(1-q)\in f \Endom_S(V') f$
and $\varphi(1-2e_{kk}')-q+(1-q)\in f\Endom_S(V')f$
for $k=1,2$.

2. We show that in \eqref{e21} matrix units can be chosen in such a~way
that
\begin{equation}\label{e32}
\varphi(1-e_{ii}-e_{jj}+e_{ij}+e_{ji})=
(q-f_{ii}-f_{jj}+f_{ij}+f_{ji})-(1-q-h_{ii}-h_{jj}+h_{ij}+h_{ji})
\end{equation}
for all $i,j\in I_1$, $i\ne j$.

Indeed, set
$$
e_{11}'=1/2(e_{11}+e_{22}-e_{12}-e_{21}),\quad
e_{22}'=1/2(e_{11}+e_{22}+e_{12}+e_{21}),\quad
e_{ii}'=e_{ii}\ \forall i\ne 1,2.
$$

The system $\{ e_{ii}'\}$ can be added to the system
of matrix units $\{ e_{ij}'\}_{i,j\in I_1}$
of the ring $\Endom_R(V)$. According to the argument in item~1,
$$
\varphi(1-e_{11}-e_{22}+e_{12}+e_{21})=\varphi(1-2e_{11}')=q-(1-q)+x,
$$
where $x\in f \Endom_S(V') f$ and $f$ is taken from~\eqref{e22}.
Consequently,
\begin{equation}\label{e33}
\varphi(1-e_{11}-e_{22}+e_{12}+e_{21})=
\begin{pmatrix}
a_{11}& a_{12}& 0& \vdots\\
a_{21}& a_{22}& 0& \vdots\\
\hdotsfor{2}& 1 & \vdots\\
\hdotsfor{3}& \ddots
\end{pmatrix}-
\begin{pmatrix}
b_{11}& b_{12}& 0& \vdots\\
b_{21}& b_{22}& 0& \vdots\\
\hdotsfor{2}& 1& \vdots\\
\hdotsfor{3}& \ddots
\end{pmatrix},
\end{equation}
where $a_{ij}\in f_{11}\Endom_S(V') f_{11}$ and
$b_{ij}\in h_{11} \Endom_S(V') h_{11}$.
Since
$$
(1-e_{11}-e_{22}+e_{12}+e_{21})(1-2e_{11})=
(1-2e_{22})(1-e_{11}-e_{22}+e_{12}+e_{21}),
$$
we have that \eqref{e21}~and~\eqref{e33} imply
\begin{align*}
\begin{pmatrix}
a_{11}& a_{12}\\
a_{21}& a_{22}
\end{pmatrix}
\begin{pmatrix}
-1& 0\\ 0& 1
\end{pmatrix} &=
\begin{pmatrix}
1& 0\\
0& -1
\end{pmatrix}
\begin{pmatrix}
a_{11}& a_{12}\\
a_{21}& a_{22}
\end{pmatrix},\\
\begin{pmatrix}
b_{11}& b_{12}\\
b_{21}& b_{22}
\end{pmatrix}
\begin{pmatrix}
-1& 0\\
0& 1
\end{pmatrix} &=
\begin{pmatrix}
1& 0\\
0& -1
\end{pmatrix}
\begin{pmatrix}
b_{11}& b_{12}\\
b_{21}& b_{22}
\end{pmatrix},
\end{align*}
and
\begin{equation}\label{e34}
a_{11}=a_{22}=0,\quad b_{11}=b_{22}=0.
\end{equation}

Then, $(1-e_{11}-e_{22}+e_{12}+e_{21})^2=1$. By
\eqref{e33}~and~\eqref{e34},
$$
a_{21}=a_{12}^{-1},\quad b_{21}=b_{12}^{-1}.
$$

Similarly,
\begin{multline*}
\varphi(1-e_{ii}+r_{i+1,i+1}+e_{i,i+1}+e_{i+1,i})\\
=(q-f_{ii}-f_{i+1,i+1}+a_i f_{ii+1}+a^{-1} f_{i+1,i})=
(1-q-h_{i+1,i+1}-h_{ii}+b_i h_{i,i+1}+b^{-1}h_{i+1,i})
\end{multline*}
for all $i\in I_1$.

Set, by transfinite induction,
$c_1\equiv 1$, $c_{i+1}\equiv c_i\cdot a_i^{-1}$, and $c_i\equiv 1$
for a~limit ordinal number~$i$.
Similarly, set
$d_1\equiv 1$, $d_{i+1}\equiv d_i\cdot b_i^{-1}$,
and $d_i\equiv 1$ for a~limit ordinal number~$i$. Let, further,
$C\equiv \mathrm{diag}( c_1,\dots, c_n,\dots)+
\mathrm{diag}(d_1,\dots,d_n,\dots)$, $h_{ij}'\equiv C h_{ij}C{^-1}$.
Then $h_{ii}'=h_{ii}$, $f_{ii}'=f_{ii}$, $f_{i,i+1}'=a_if_{i,i+1}$,
$f_{i+1,i}=a_i^{-1}f_{i+1,i}$, $h_{i,i+1}'=b_1 h_{i,i+1}$, and
$h_{i+1,i}'=b_i^{-1} h_{i+1,i}$.

Thus,
\begin{multline*}
\varphi(1-e_{ii}-e_{i+1,i+1}+e_{i,i+1}+e_{i+1,i})\\
=(q-f_{ii}'-f_{i+1,i+1}'+f_{i,i+1}'+f_{i+1,i}')-
(1-q-h_{ii}'-h_{i+1,i}'+h_{i,i+1}'+h_{i+1,i}').
\end{multline*}
Finally, the assertion~2 is proved.

3. Set $g_{ij}=f_{ij}+h_{ij}$, where $f_{ij}$,~$h_{ij}$ are
matrix units for which conditions \eqref{e21} and~\eqref{e32}
hold. Then $\{ g_{ij}\}_{i,j\in \varkappa}$ is a~system of matrix units
of the ring $\Endom_S(V')$.
An arbitrary element $C\in \Endom_S(V')$ will be written
in the form
$$
C\equiv \begin{pmatrix}
c_{11}& \dots& c_{1n}& \vdots\\
\dots& \ddots& \dots& \vdots\\
c_{n1}&\dots&c_{nn}& \vdots\\
\hdotsfor{3}& \ddots
\end{pmatrix},\ \
\text{where}\ \ c_{ij}\in g_{ij} \Endom_S(V') g_{ij}.
$$

4. We show that for every element $r\in R$,
\begin{equation}\label{e35}
\varphi(1+re_{12})=
\begin{pmatrix}
a_r& b_r& 0& \vdots\\
c_r& d_r& 0& \vdots\\
\hdotsfor{2}& 1& \vdots\\
\hdotsfor{3}&\ddots
\end{pmatrix},
\end{equation}
where $a_r,b_r,c_r,d-r\in g_{11}\Endom_S(V')g_{11}$, and that
\begin{equation}\label{e36}
\varphi(1-e_{ii}-e_{jj}+e_{ij}-e_{ji})=1-g_{ii}-g_{jj}+g_{ij}-g_{ji}
\end{equation}
for $i\ne j$.

Indeed,
$$
1-e_{ii}-e_{jj}+e_{ij}-e_{ji}=(1-e_{ii}-e_{jj}+e_{ij}+e_{ji})(1-2e_{ii}).
$$
By \eqref{e21}~and~\eqref{e32},
\begin{multline*}
\varphi(1-e_{ii}-e_{jj}+e_{ij}-e_{ji})\\
=(e-f_{ii}-f_{jj}+f_{ij}-f_{ji})+(1-e-h_{ii}-h_{jj}+h_{ij}-h_{ji})=
1-g_{ii}-g_{jj}+g_{ij}-g_{ji}.
\end{multline*}
Set
$$
e_{ij}''=(1+1/2re_{12})e_{ij}(1+1/2re_{12})^{-1}.
$$
Then, according to the assertion~1,
$$
\varphi(1-2e_{11}'')=q-(1-q)+x
$$
and
$$
x_1\in f \Endom_S(V') f,\ \ \text{where}\ \ f=h_{11}+h_{22}+f_{11}+f_{22}.
$$
Then,
$$
1-2e_{11}''=1-2e_{11}+r e_{12}=(1+re_{12})(1-2e_{11})
$$
and, by~\eqref{e21},
$$
\varphi(1+re_{12})=\varphi((1-2e_{11}'')(1-2e_{11}))=1+x_2,
$$
where $x_2\in f \Endom_S(V') f$.

But from $f \Endom_S(V') f=(g_{11}+g_{22})\Endom_S(V')(g_{11}+g_{22})$
it follows that
$$
\varphi(1+re_{12})=
\begin{pmatrix}
a_r& b_r& 0&\vdots\\
c_r& d_r& 0&\vdots\\
\hdotsfor{2}& 1& \vdots\\
\hdotsfor{3}&\ddots
\end{pmatrix}.
$$

5. Using the equalities \eqref{e35}~and~\eqref{e36}, and the equality
$$
\begin{pmatrix}
1& 0& 0\\
0& 0& -1\\
0& 1& 0
\end{pmatrix}
\begin{pmatrix}
1& r& 0\\
0& 1& 0\\
0& 0& 1
\end{pmatrix}
\begin{pmatrix}
1& 0& 0\\
0& 0& 1\\
0& -1& 0
\end{pmatrix}=
\begin{pmatrix}
1& 0& r\\
0& 1& 0\\
0& 0& 1
\end{pmatrix},
$$
we shall get
\begin{equation}\label{e37}
\varphi(1+re_{13})=
\begin{pmatrix}
a_r& 0& b_r& 0 & \vdots\\
0& 1& 0& 0& \vdots\\
c_r & 0& d_r& 0& \vdots\\
0& 0& 0& 1& \vdots\\
\hdotsfor{4}& \ddots
\end{pmatrix},
\end{equation}
where $a_r,b_r,c_r,d_r$ are taken from~\eqref{e35}.

From \eqref{e35}~and~\eqref{e37} we have that for all $r,s\in R$
$$
\begin{pmatrix}
a_r& b_r& 0\\
c_r& d_r& 0\\
0&0&1
\end{pmatrix}
\begin{pmatrix}
a_s& 0& b_s\\
0& 1& 0\\
c_s& 0& d_s
\end{pmatrix}=
\begin{pmatrix}
a_s& 0& b_s\\
0& 1& 0\\
c_s& 0& d_s
\end{pmatrix}
\begin{pmatrix}
a_r &b_r& 0\\
c_r& d_r& 0\\
0& 0& 1
\end{pmatrix}
$$
and
\begin{equation}\label{e38}
b_r=a_sb_r,\quad
c_ra_s=c_r,\quad
c_rb_s=0.
\end{equation}

Similarly, using the equalities
\begin{equation}
\label{e39}
\varphi(1+r e_{23})=
\begin{pmatrix}
1& 0& 0& 0& \vdots\\
0& a_r& b_r& 0& \vdots\\
0& c_r& d_r& 0& \vdots\\
0& 0& 0& 1& \vdots\\
\hdotsfor{4}& \ddots
\end{pmatrix}
\end{equation}
and $[1+se_{23}, 1+re_{13}]$, we have that for all
$r,s\in R$
\begin{equation}\label{e40}
b_r=b_r d_s,\quad
d_s c_r=c_r,\quad
b_sc_r=0.
\end{equation}
From the equalities
$$
(1+re_{ij})^{-1}=(1-re_{ij})=(1-2e_{ii})(1+r_{ij})(1-2e_{ii})
$$
and \eqref{e21},~\eqref{e35} it follows that for all $r\in R$
\begin{equation}\label{e41}
\varphi(1+re_{12})^{-1}=
\begin{pmatrix}
a_r& -b_r& 0&\vdots\\
-c_r& d_r& 0&\vdots\\
0& 0& 1& \vdots\\
\hdotsfor{3}& \ddots
\end{pmatrix},
\end{equation}
and, according to \eqref{e41} and~\eqref{e40},
\begin{equation}\label{e42}
a_r^2=d_r^2=1.
\end{equation}
From the equalities
$$
\begin{pmatrix}
0& 1\\
-1& 0
\end{pmatrix}
\begin{pmatrix}
1& 1\\
0& 1
\end{pmatrix}
\begin{pmatrix}
0& -1\\
1& 0
\end{pmatrix}=
\begin{pmatrix}
1& 0\\
-1& 1
\end{pmatrix}
$$
and \eqref{e35},~\eqref{e36}, we have that
\begin{equation}\label{e43}
\varphi(1-e_{21})=
\begin{pmatrix}
d_1& -c_1& \vdots\\
-b_1& a_1& \vdots\\
\hdotsfor{2}& \ddots
\end{pmatrix}.
\end{equation}

Then,
$$
\begin{pmatrix}
1& 1\\
0& 1
\end{pmatrix}
\begin{pmatrix}
1& 0\\
-1& 1
\end{pmatrix}=
\begin{pmatrix}
0& 1\\
-1& 0
\end{pmatrix}
\begin{pmatrix}
1& -1\\
0& 1
\end{pmatrix}.
$$
From \eqref{e35}, \eqref{e43}, \eqref{e36}, and \eqref{e40} it follows that
\begin{align}
a_1d_1-b_1^2&=-c_1,\label{e44}\\
-c_1^2+d_1a_1&=b_1.\label{e45}
\end{align}
Multiply the equality~\eqref{e45} from the right side by~$b_1$
and, using~\eqref{e38}, we shall get $b_1^2=d_1b_1$. Multiplying~\eqref{e45}
from the left side by~$b_1$, we shall get
$b_1^2=b_1a_1$. Therefore, we have shown that
$b_1a_1=d_1b_1=b_1^2$, $d_1b_1d_1b_1=d_1b_1^2=d_1^2b_1=b_1$,
$d_1b_1d_1b_1=d_1b_1^2a_1=d_1^2b_1a_1=b_1a_1b_1^2$, and
$b_1=b_1^2$.

From \eqref{e44} it follows that
\begin{equation}\label{e46}
a_1c_1=c_1d_1=c_1=-c_1^2.
\end{equation}
From \eqref{e45}~and~\eqref{e46} we have
$$
d_1a_1=b_1+c_1^2=b_1c_1.
$$
From \eqref{e38}~and~\eqref{e40} it follows that
$$
b_1c_1=c_1b_1=0.
$$

Therefore,
$$
(d_1a_1)^2=b_1^2+c_1^2=b_1-c_1=d_1a_1.
$$
According to \eqref{e42}, the element $d_1a_2$ is invertible.
Consequently,
\begin{equation}\label{e47}
1=d_1a_1=b_1-c_1.
\end{equation}
By \eqref{e38}, \eqref{e40}, and \eqref{e47}, $b_sc_r=c_rb_s=0$ and
\begin{equation}\label{e48}
b_r\in b_1 f_{11}\Endom_S(V')f_{11}b_1,\quad
c_r\in (1-b_1)f_{11}\Endom_S(V')f_{11}(1-b_1)
\end{equation}
for all $r,s\in R$. Then, according to~\eqref{e38},
$$
(a_s-1)b_1=c_1(a_s-1)=0.
$$
By \eqref{e47},
$$
a_s-1=-b_1(a_s-1)c_1.
$$
By \eqref{e38}, \eqref{e40}, and \eqref{e42},
$$
b_1c_1=c_1b_1=0,\quad
1=a_s^2=(1-b_1a_sc_1)^2=1-2b_1a_sc_1,
$$
and $a_s=1$. Similarly, $d_s=1$. Thus,
\begin{equation}\label{e49}
a_r=d_r=1
\end{equation}
for all $r\in R$. Set $e_1=b_1\cdot 1$, then $e_1$ is an idempotent
of the ring $\Endom_S(V')$.
By \eqref{e35}, \eqref{e48}, and \eqref{e49},
\begin{gather*}
e_1\varphi(1+re_{12})=\varphi(1+re_{12})e_1=e_1+b_rg_{12},\\
[1-2e_1,\varphi(1+re_{12})]=1.
\end{gather*}
Similarly,
$$
[1-2e_1,\varphi(1-e_{ii}-e_{jj}+e_{ij}-e_{ji})]=1.
$$
Consequently, the matrix $\varphi^{-1}(1-2e_1)$ belongs
to the centralizer of the group $E_R(V)$ and is a~central
matrix. Therefore, the matrix $1-2e_1$ belongs to the center
of the ring $\Endom_S(V')$, $e_1$ is a~central idempotent
of the ring $\Endom_S(V')$,
\begin{equation}\label{e50}
e_1\Endom_S(V')\oplus (1-e_1) \Endom_S(V')=\Endom_S(V').
\end{equation}
Set
$\theta_3(r)\equiv b_r$ and $\theta_4(r)\equiv -c_r$. From the equalities
\begin{align*}
[1+re_{12},1-se_{23}]&=1+(rs)e_{13},\\
[1+c_r g_{21},1-c_s g_{32}]&=1-(c_sc_r)g_{31}
\end{align*}
and \eqref{e35}, \eqref{e48}, \eqref{e49}, \eqref{e37}, \eqref{e39},
and \eqref{e41} it follows that
$$
\begin{pmatrix}
1& 0& b_{rs}\\
0& 1& 0\\
c_{rs}& 0& 1
\end{pmatrix}=
\begin{pmatrix}
1& 0& b_r b_s\\
0& 1& 0\\
-c_sc_r& 0& 1
\end{pmatrix}.
$$
Hence, $\theta_3\colon R\to b_1(f_{11}\Endom_S(V')f_{11})$ is a~homomorphism
of rings and $\theta_4\colon R\to (1-b_1)(f_{11}\Endom_S(V')f_{11})$ is
an antihomomorphism of rings. Further, by \eqref{e35}, \eqref{e36}, and
\eqref{e49},
\begin{equation}\label{e51}
\varphi(1+re_{ij})=1+\theta_3(r)g_{ij}-\theta_4(r)g_{ji}.
\end{equation}

Set, for every $a_{ij}e_{ij}\in \Endom_R(V)$,
\begin{align*}
\theta_1(a)&= \theta_3(a_{ij})g_{ij},\\
\theta_2(a)&=\theta_4(a_{ij})g_{ji},
\end{align*}
and for other elements of the ring $\Endom_R(V)$ we continue these
homomorphisms in the natural way. Then
$\theta_1\colon \Endom_R(V)\to e_1 \Endom_S(V')$ is
a~homomorphism of rings,
$\theta_2\colon \Endom_R(V)\to (1-e_1)\Endom_S(V')$ is
an antihomomorphism of rings, and, by \eqref{e47} and \eqref{e51},
\begin{equation}\label{e52}
\varphi(A)=\theta_1(A)+\theta_2(A^{-1})
\end{equation}
for all $A\in E_R(V)$. Let $I$,~$J$ be ideals of the ring~$S$
such that $\Endom_I(V')=e_1 \Endom_S(V')$ and
$\Endom_J(V')=(1-e_1)\Endom_S(V')$.
By~\eqref{e50}, $I\oplus J=S$. Set
$N_1\equiv \varphi^{-1}(\Aut_S(V',I))$ and
$M_1\equiv \varphi^{-1}(\Aut_S(V',J))$. By Lemma~\ref{mih-l-1},
$$
E_R(V,eR)\subseteq N_1,\quad
E_R(V,(1-e)R)\subseteq M_1,
$$
where $e$ is some central idempotent of the ring $\Endom_R(V)$.
Let $B\in E_R(V,eR)$. Then $\varphi(B)-1\in \Endom_I(V')=e_1\Endom_S(V')$.

By \eqref{e52},
$$
\varphi(B)-1=\theta_1 (B-1)+\theta_2(B^{-1}-1)
$$
and
$$
\theta_1(B-1)\in e_2 \Endom_S(V'),\quad
\theta_2 (B^{-1}-1)\in (1-e_2)\Endom_S(V').
$$
Consequently, $\theta_2(B^{-1}-1)=0$ and
$\Endom_{eR}(V)\subseteq \Ker \theta_2$. Similarly,
$\Endom_{(1-e)R}(V)\subseteq \Ker\theta_1$. Since $\varphi$ is
a~group isomorphism, we have by~\eqref{e52}
$$
\Ker \theta_1\cap \Ker\theta_2=\{ 0\}.
$$

Therefore,
$$
\Endom_{eR}(V)=\Ker\theta_2,\quad
\Endom_{(1-e)R}(V)=\Ker \theta_1,
$$
and
$$
\Ker \theta_1\oplus \Ker\theta=\Endom_R(V).
$$
A~similar argument for the mapping $\varphi^{-1}$ leads us to
$$
\Image \theta_1\oplus \Image \theta_2=\Endom_S(V').
$$
Set
$$
\varphi_1(B)=\varphi^{-1}(\theta_1(B)+\theta_2(B^{-1}))
$$
for all $B\in \Aut_R(V)$. Then $\varphi_1$ is an automorphism
of the group $\Aut_R(V)$, and, by~\eqref{e52},
$$
\varphi_1(A)=A\ \ \text{for all}\ \ A\in E_R(V).
$$
The theorem is proved.
\end{proof}

Suppose that rings $R$~and~$S$ with $1/2$ do not contain any
central idempotents which are not equal to $0$~or~$1$.
Then we have the following theorem.

\begin{theorem}\label{mih-t-2}
The groups $\Aut_R(V)$ and $\Aut_S(V')$ are isomorphic if and only if
$\Endom_R(V)\cong \Endom_S(V')$.
\end{theorem}

\begin{proof}
By Theorem~\ref{mih-t-1},
on the group $DE_R(V)$
every isomorphism~$\varphi$ of the groups $\Aut_R(V)$ and $\Aut_S(V')$
coincides with an isomorphism
$\chi(\cdot)(\theta_1(\cdot)+\theta_2(\cdot^{-1}))$,
where $\chi(\cdot)$ is a~group homomorphism $DE_R(V)\to C(\Aut_S(V'))$,
$\theta_1\colon e\Endom_R(V)\to f \Endom_S(V')$ is a~ring isomorphism,
$\theta_2\colon (1-e)\Endom_R(V)\to (1-f)\Endom_S(V')$
is a~ring antiisomorphism,
$e$,~$f$ are central idempotents of the rings
$\Endom_R(V_1)$ and $M_S(V_2)$, respectively.
Since the rings $R$~and~$S$ do not contain any central
idempotents which are not equal to
$0$~or~$1$, we have that the rings $\Endom_R(V)$ and $\Endom_S(V')$
also do not contain any central idempotents which are not equal
to $0$~and~$1$, i.e., either $e=f=1$, or $e=f=0$.

1. If $e=f=1$, then $\varphi(\cdot)$ on $DE_R(V)$ coincides
with an isomorphism of the rings $\Endom_R(V)$ and $\Endom_S(V')$
of the form $\chi(\cdot)\theta_1(\cdot)$, i.e.,
the rings $\Endom_R(V)$ and $\Endom_S(V')$ are isomorphic.

2. If $e=f=0$, then $\varphi$ on $DE_R(V)$ coincides with an
antiisomorphism $\chi(\cdot)\theta_2(\cdot^{-1})$, i.e.,
the rings $\Endom_R(V)$ and $\Endom_S(V')^{\mathrm{op}}$ are isomorphic.

Suppose that we have this case. Consider in $\Endom_R(V)$
a~system of commuting conjugate orthogonal
idempotents with the condition
$$
\sum_{i\in I} e_{ii}\sim 1.
$$
This expression means that for every element~$a$ and every $i\in I$
there exist $i_1,\dots,i_n\in \varkappa$ such that
$$
\biggl(\,\sum_{j=1}^n e_{i_j i_j}\biggr) ae_{ii}=
ae_{ii}\biggl(\,\sum_{j=1}^n e_{i_ji_j}\biggr)=ae_{ii}.
$$

Now, as above, introduce a~\emph{system of matrix units}~$e_{ij}$
($i,j\in \varkappa$) by the condition
$$
e_{ij}e_{kl}=\delta_{jk} e_{il}.
$$

It is clear that such system $\{ e_{ij}\}$ in $\Endom_R(V)$ corresponds
to a~system $\{ f_{ij}\}$ in $\Endom_S(V')$, defined by
the condition
$$
f_{ij} f_{kl}=\delta_{il} f_{kj}.
$$

In $\Endom_R(V)$ there exists an element
$$
x\sim\sum_{i\in I} e_{1i},
$$
but in $\Endom_S(V_2)$ a~corresponding element
$$
y\sim\sum_{i\in I} f_{1i}
$$
can not exist.

We show this.

Let $W_i$ be the carrier of an idempotent $f_{ii}$. Then
$$
f_{1i}(W_1)=f_{ii} f_{1i}(W_1) \Rightarrow f_{1i}(W_1)\subset W_i.
$$
Further,
$$
W_1=f_{11}(W_1)=f_{i1}f_{1i}(W_1),
$$
i.e., $f_{ij}$ maps $W_j$ to $W_i$. Existence
of the element $f\sim\sum\limits_{i\in I} f_{1i}$ would mean
that $f$ maps some vector~$w$ from~$W_1$ to the sum
of an infinite number of vectors $w_j\in W_j$, but this is impossible.

Therefore, the condition
$$
\Endom_R(V)\cong \Endom_S(V')^{\mathrm{op}}
$$
is impossible.

The inverse implication is evident.
\end{proof}

\subsection[Elementary Equivalence of Automorphism Groups and
Endomorphism Rings of Modules of Infinite Ranks]{Elementary
Equivalence of Automorphism Groups and
Endomorphism Rings of \mbox{Modules} of Infinite Ranks}\label{ss5.2}
\begin{lemma}\label{mih-l-3}
For every ultrafilter $D$
$$
\prod_D \Endom_R(V)\cong \Endom_{\prod\limits_D}(V).
$$
\end{lemma}

\begin{proof}
By the definition of ultraproduct, every element $\prod\limits_D \Endom_R(V)$
is a~mapping (more precisely,
its equivalence class) $f\colon I\to \Endom_R(V)$,
i.e., a~set of pairs $\langle i,A\rangle$,
where $i\in I$, $A\in \Endom_R(V)$,
$\forall i\in I\losp \exists ! A\in \Endom_R(V)\losp
(\langle i,A\rangle\in f)$.
Every element of $A\in \Endom_R(V)$ is a~mapping
$a\colon \varkappa\times \varkappa \to R$ such that for every
$\alpha\in \varkappa$ there exists
only a~finite number of $\beta_j\in \varkappa$ such that
$a(\langle \alpha,\beta_j\rangle) \ne 0$, i.e.,
every element of $A\in \Endom_R(V)$ is a~set
of ordered triplets $\langle \alpha,\beta,r\rangle$, where
$\alpha,\beta\in \varkappa$, $r\in R$,
$\forall \alpha\losp \forall \beta\losp \exists! r\in R\losp
(\langle \alpha,\beta,r\rangle \in A)$.
Therefore, every element of the ultraproduct $\prod\limits_D \Endom_R(V)$
is a~set~$f$ of ordered quadruplets $\langle i,\alpha,\beta,r\rangle$
with $i\in I$, $\alpha,\beta\in \varkappa$, $r\in R$ and with the condition
$\forall i,\alpha,\beta\losp \exists ! r\losp
(\langle i,\alpha,\beta,r\rangle\in f)$.
In other words, it is a~mapping
$f\colon I\times \varkappa\times \varkappa\to R$
with the only condition that
for every $i\in I$ and $\alpha\in \varkappa$
there exist only a~finite number of
$\beta_j\in \varkappa$ such that $f(i,\alpha,\beta_j)\ne 0$.

Two such mappings $f,g\colon I\times \varkappa\times \varkappa\to R$
are equal if and only if
$$
\{ i\in I\mid \forall \alpha,\beta\in \varkappa\losp
(f(i,\alpha,\beta)=g(i,\alpha,\beta))\} \in D.
$$
For three mappings $f,g,h\colon I\times \varkappa\times \varkappa\to R$
we have $h=f+g$ if and only if
$$
\{ i\in I\mid \forall\alpha,\beta\in \varkappa\losp
(h(i,\alpha,\beta)=
f(i,\alpha,\beta)+g(i,\alpha,\beta))\}\in D.
$$
Similarly, for three mappings
$f,g,h\colon I\times \varkappa\times \varkappa\to R$
we have $h=fg$ if and only if
$$
\biggl\{ i\in I \biggm| \forall \alpha,\beta\in \varkappa\losp
\biggl(h(i,\alpha,\beta)=\sum_{\gamma\in \varkappa}
f(i,\alpha,\gamma)\cdot g(i,\gamma,\beta)\biggr)\biggr\}\in D.
$$
It is clear that we can write the sign of sum in this expression
because only a~finite number of elements of this sum
are nonzero.

Now consider the ring $\Endom_{\prod\limits_D}(V)$.
Completely the same arguments lead us to the fact that
the elements of this ring are mappings
$f\colon \varkappa\times \varkappa\times I\to R$ with the same condition
of finiteness and the same identity, sum, and product.
Therefore, the obtained isomorphism is natural
(it is the natural mapping
$I\times (\varkappa\times \varkappa)\to
(\varkappa \times \varkappa)\times I$).
\end{proof}

The proof of the following theorem is similar to the proof of
Theorem~4 of the paper~\cite{2}.

\begin{theorem}\label{4_2}\label{th3_5.2}
Suppose that rings $R$ and $S$ contain $1/2$ and do not contain any
central idempotents which are not equal to $1$~and~$0$.
Suppose that $V$ and~$V'$ are
free modules of infinite ranks over the rings $R$~and~$S$, respectively.
Then the groups $\Aut_R(V)$ and $\Aut_S(V')$
are elementarily equivalent if and only if
the rings $\Endom_R(V)$ and $\Endom_S(V')$ are elementarily equivalent.
\end{theorem}

\begin{proof}
Let the rings $\Endom_R(V)$ and $\Endom_S(V')$ be elementarily equivalent.

Consider an arbitrary sentence~$\varphi$ of the first order language of
the group theory. With the help of the sentence~$\varphi$ we construct
a~sentence~$\varphi'$ of the first order language of the ring theory
in the following way: every symbol-string of the form
$\forall x\losp (\dots)$
belonging to the sentence~$\varphi$ will be replaced by the symbol-string
$\forall x\losp (\exists x'\losp (xx'=x'x=1)\Rightarrow (\dots)$, and
every symbol-string of the form $\exists x\losp (\dots)$
will be replaced by the symbol-string
$\exists x\losp (\exists x'\losp (xx'=x'x=1)\logic\land(\dots))$.
It is clear that if the sentence~$\varphi$ holds in the group $\Aut_R(V)$,
then the sentence~$\varphi'$ holds in the ring $\Endom_R(V)$, and, therefore,
since the rings $\Endom_R(V)$ and $\Endom_S(V')$ are elementarily equivalent,
we have that it holds in the ring $\Endom_S(V')$. Consequently,
the sentence~$\varphi$ holds in the group $\Aut_S(V')$.
Now we see that the groups $\Aut_R(V)$ and $\Aut_S(V')$
are elementarily equivalent.

Let the operation~$*$ applied to some ring $A$ ($A^*$) be
taking the group of invertible elements of this ring.
It is clear that for every ultrafilter~$D$
$\smash[b]{\prod\limits_D \Aut_R(V)=
\prod\limits_D (\Endom_R(V))^* \cong
\Bigl(\,\prod\limits_D \Endom_R(V)\Bigr)^*}$,
i.e., that the operations~$*$ and $\prod\limits_D$ are permutable.

Let now the groups $\Aut_R(V)$ and $\Aut_S(V')$ be elementarily
equivalent. Then, by Theorem~\ref{ult-t-2} in Sec.~1.4, there
exist ultrapowers $G= \prod\limits_D \Aut_R(V)$ and
$G'=\prod\limits_{D}\Aut_S(V')$ of these groups
such that $G\cong G'$. Therefore,
$\Bigl(\,\prod\limits_D \Endom_R(V)\Bigr)^*
\cong \Bigl(\,\prod\limits_D \Endom_S(V')\Bigr)^*$,
and, by Lemma~\ref{mih-l-3},
$\Aut_{\prod\limits_D R}(V)\cong \Aut_{\prod\limits_D S}(V')$.
By Theorem~\ref{mih-t-2} from the previous subsection,
in this case
$\Endom_{\prod\limits_D R}(V)\cong \Endom_{\prod\limits_D S}(V')$.
Consequently, by Proposition~\ref{ult-t-1} in Sec.~1.4,
$\Endom_R(V)\equiv \Endom_S(V')$.
The theorem is proved.
\end{proof}

Therefore, in the case where we have associative rings with~$1/2$
which do not contain any central idempotents not equal to $0$~and~$1$
we can replace the question on elementary equivalence
of automorphism groups by the question on elementary equivalence of
endomorphism rings.

\subsection{The Main Theorem}\label{ss5.3}
In this section, we assume that a~cardinal number~$\varkappa_1$
is such that there exists a~maximal ideal of the ring~$R_1$
generated by at most~$\varkappa_1$ elements.

From Theorem \ref{3osn} in Sec.~3 and Theorem~\ref{4_2}
we easily obtain Theorem~\ref{4osn}.

\begin{theorem}\label{4osn}\label{th4_5.3}
Suppose that rings $R_1$~and~$R_2$ contain~$1/2$ and do not contain
any central idempotents which are not equal to $1$~or~$0$. Let
$V_1$~and~$V_2$ be free modules of infinite ranks
$\varkappa_1$~and~$\varkappa_2$ over the rings $R_1$~and~$R_2$,
respectively, and let
$\psi \in \Th_2^{\varkappa_1}(\langle \varkappa_1,R_1\rangle)$
be such that
$\psi\notin \Th_2^{\varkappa_1} (\langle \varkappa_1, R'\rangle )$
for any ring~$R'$ such that $R'$ is similar to~$R_1$ and
$\Th_2^{\varkappa_1}(\langle \varkappa_1,R_1\rangle )\ne
\Th_2^{\varkappa_1}(\varkappa_1, R'\rangle)$.
Then the groups $\Aut_{R_1}(V_1)$ and $\Aut_{R_2}(V_2)$ are elementarily
equivalent if and only if there exists a~ring~$S$
similar to the ring~$R_2$ and such that the theories
$\Th_2^{\varkappa_1}(\langle \varkappa_1,R_1\rangle)$ and
$\Th_2^{\varkappa_2}(\langle \varkappa_2,S\rangle)$ coincide.
\end{theorem}

\begin{corollary}\label{co1_5.3}
For free modules $V_1$~and~$V_2$ of infinite ranks
$\varkappa_1$~and~$\varkappa_2$ over skewfields
\textup{(}integral domains, commutative or local
rings without central idempotents not equal to $1$~or~$0$\textup{)}
$F_1$~and~$F_2$ with~$1/2$, respectively,
the groups $\Aut_{F_1}(V_1)$ and $\Aut_{F_2}(V_2)$
are elementarily equivalent if and only if the theories
$\Th_2^{\varkappa_1}(\langle \varkappa_1, F_1\rangle)$
and $\Th_2^{\varkappa_2}(\langle \varkappa_2, F_2\rangle)$
coincide.
\end{corollary}

\begin{corollary}\label{co2_5.3}
For free modules $V_1$~and~$V_2$ of infinite ranks
$\varkappa_1$~and~$\varkappa_2$ over Artinian
rings $R_1$~and~$R_2$
with~$1/2$ without central idempotents not equal
to $1$~or~$0$, respectively,
the groups $\Aut_{R_1}(V_1)$ and $\Aut_{R_2}(V_2)$
are elementarily equivalent
if and only if there exist rings $S_1$~and~$S_2$ such that
the ring~$R_1$ is similar to the ring~$S_1$, the ring~$R_2$
is similar to the ring~$S_2$, and the theories
$\Th_2^{\varkappa_1}(\langle \varkappa_1, S_1\rangle)$
and $\Th_2^{\varkappa_2}(\langle \varkappa_2, S_2\rangle)$
coincide.
\end{corollary}

\end{document}